\documentclass{article}
\pdfoutput=1
\usepackage{amsmath} 
\usepackage{amsfonts}
\usepackage{amssymb}
\usepackage{geometry}
\usepackage{makeidx}
\usepackage{color}
\begin{document}

\title {On logarithmic integrals, harmonic sums and variations}
\author{Ming Hao Zhao}
\date{}

\maketitle
\begin{abstract}
Based on various non-MZV approaches we evaluate certain logarithmic integrals and harmonic sums. More specifically, 85 LIs, 89 ESs, 263 PLIs, 28 non-alt QESs, 39 GESs, 26 BESs, 14 IBESs with weight $\leq5$, 193 QLIs, 172 QPLIs, 83 QESs, 7 QBESs, 3 IQBESs with weight $\leq4$. With help of MZV theory we evaluate other 10 QESs with weight $\leq4$. High weight hypergeometric series, nonhomogeneous integrals and other results are also established.
\end{abstract}

\bibliographystyle{plain}

\tableofcontents

\section{Preliminaries}

\subsection{Special functions}

Below are special functions used. We assume readers are aware of their basic properties \cite{ref13, ref14}.\\

Gamma function: The analytic continuation of $\Gamma(s)=\int_0^{\infty} t^{s-1}e^{-t}dt$\\

Euler Beta function: The analytic continuation of $B(a,b)=\int_0^1 t^{a-1}(1-t)^{b-1}dt= \frac{\Gamma (a) \Gamma (b)}{\Gamma (a+b)}$\\

Polygamma function: The analytic continuation of $\psi ^{(n)}(x)=\frac{\partial ^{n+1}\log (\Gamma (x))}{\partial x^{n+1}}$\\

Polylogarithm: The analytic continuation of $\text{Li}_n(z)=\sum _{k=1}^{\infty } \frac{z^k}{k^n}$\\

Riemann Zeta function: The analytic continuation of $\zeta (s)=\sum _{k=1}^{\infty } \frac{1}{k^s}$\\

Dirichlet Eta function: The analytic continuation of $\eta (s)=(1-2^{1-s})\zeta(s) =\sum _{k=1}^{\infty } \frac{(-1)^{k-1}}{k^s}$\\

Dirichlet Beta function: The analytic continuation of $\beta (s)=\sum _{k=0}^{\infty } \frac{(-1)^k}{(2 k+1)^s}$\\

Catalan constant and generalization: $C=\beta(2), \frac{\psi ^{(3)}\left(\frac{1}{4}\right)-\psi ^{(3)}\left(\frac{3}{4}\right)}{1536}=\beta(4)$\\

Euler-Mascheroni constant: $\gamma=-\psi^{(0)}(1)$\ \ \ \ \  Pochhammer symbol: $(a)_n=\frac{\Gamma (a+n)}{\Gamma (a)}$\\

Generalized hypergeometric function: The analytic continuation of $$\, _pF_q(a_1,\cdots,a_p;b_1,\cdots,b_q;z)=\sum _{n=0}^{\infty } \frac{ (a_1)_n \cdots (a_p)_n z^n}{(b_1)_n \cdots (b_q)_n n!}$$

Clausen function: $\text{Cl}_{2k}(x)=\sum _{n=1}^{\infty } \frac{\sin (n x)}{n^{2k}}, \text{Cl}_{2k+1}(x)=\sum _{n=1}^{\infty } \frac{\cos (n x)}{n^{2k+1}}$\\

Inverse tangent integral: $\text{Ti}_k(x)=\int_0^x \frac{\text{Ti}_{k-1}(t)}{t}dt, \text{Ti}_1(x)=\tan^{-1}(x)$\\

Elliptic K/E integral: $K(x)=\int_0^{\frac{\pi }{2}} \frac{1}{\sqrt{1-x \sin ^2(t)}} \, dt, E(x)=\int_0^{\frac{\pi }{2}} \sqrt{1-x \sin ^2(t)} \, dt$\\

Legendre P polynomial: $P_n(x)=\frac{1}{(2 n)\text{!!}}\frac{\partial ^n\left(x^2-1\right)^n}{\partial x^n}$\\

\subsection{Definition of subjects}
Clarification: Throughout this paper, LI will be the abbreviation of log integral, PLI polylog integral, QLI quadratic log integral, QPLI quadratic polylog integral, ES Euler sum, QES quadratic Euler sum, MZV multiple Zeta value. Minor abbreviations will be mentioned later.\\

\subsubsection{LI and PLI}

Notations:

$$f(0;x)=1-x, f(1;x)=x, f(2;x)=1+x$$

$$l(0;x)=\log (1-x), l(1;x)=\log (x), l(2;x)=\log (1+x), l(3;x)=\log(1+x^2), l(4;x)=\tan^{-1}(x)$$

$$L(n,1;x)=\text{Li}_n(x), L(n,2;x)=\text{Li}_n(-x), L(n,3;x)=\text{Li}_n(\frac{1-x}{2}),n\geq2$$

$$L(n,4;x)=\text{Li}_n(\frac{1+x}{2}), L(n,5;x)=\text{Li}_n(\frac{1-x}{1+x}), L(n,6;x)=\text{Li}_n(\frac{x-1}{1+x}),n\geq2$$\\

Definition:

$$LI(a(0),a(1),a(2);p)=\int_0^1 \frac{\prod _{m=0}^2 l(m;x)^{a(m)}}{f(p;x)} \, dx$$

$$PLI(a(0),a(1),a(2);b(2,1),\cdots,b(2,6);\cdots;b(N,1),\cdots,b(N,6);p)$$ $$=\int_0^1 \frac{\prod _{m=0}^2 l(m;x)^{a(m)} \prod _{k=1}^6 \prod _{n=2}^N L(n,k;x)^{b(n,k)}}{f(p;x)} \, dx$$\\

Abbreviation: In later discussion of PLIs, $N$ and $b(n,k)$ will be small and many of them equal to $0$. In this case, we omit those zero parameters and simply repeat the word '$nk$' $b(n,k)$ times for abbreviation. For instance we have \\
$$LI(1,2,3;0)=\int_0^1 \frac{\log (1-x) \log ^2(x) \log ^3(x+1)}{1-x} \, dx$$
$$PLI(0,0,1;24,24;32;2)=\int_0^1 \frac{\text{Li}_2\left(\frac{x+1}{2}\right)^2 \text{Li}_3(-x) \log (x+1)}{x+1} \, dx$$\\ 

Weight: We define the weight $W$ of LI, PLI as $$1+\sum _{m=0}^2 a(m),\ \ 1+\sum _{k=1}^6 \sum _{n=2}^N n b(n,k)+\sum _{m=0}^2 a(m)$$ respectively, so $LI(1,2,3;0)$ is of weight 7 and $PLI(0,0,1;24,24;32;2)$ weight 9.\\

For $W\leq 5$ there are 85 convergent LIs, since a simple counting shows the total number of weight $W$ convergent LIs is $\frac{1}{2} \left(3 W^2+W-2\right)$. For $W\leq 5$, after excluding those already contained in LIs, there are 263 convergent PLIs.\\

\subsubsection{QLI and QPLI}

Notations:
$$g(1;x)=1-x, g(2;x)=x, g(3;x)=1+x, g(4;x)=1+x^2, g(5;x)=\frac{1+x^2}{x}$$

$$k(1;x)=\log (1-x), k(2;x)=\log (x), k(3;x)=\log (1+x), k(4;x)=\log(1+x^2), k(5;x)=\tan^{-1}(x)$$

$$L(n,7;x)=\text{Li}_n(-x^2), L(n,8;x)=\text{Ti}_n(x)$$

Definition:
$$QLI(A(1),A(2),A(3),A(4),A(5);p)=\int_0^1 \frac{\prod _{m=1}^5 k(m;x)^{A(m)}}{g(p;x)} \, dx$$

{\small $$QPLI3(a;p)=\int_0^1 \frac{L(2,a;x)}{g(p;x)}dx, QPLI4(a;p)=\int_0^1 \frac{L(3,a;x)}{g(p;x)}dx, QPLI4(a;b;p)=\int_0^1 \frac{L(2,a;x)k(b;x)}{g(p;x)}dx$$}

Abbreviation: In later discussion of QLIs, $A(m)$ will be small and many of them equal to $0$. In this case, we omit those zero parameters and simply repeat the word '$m$' $A(m)$ times for abbreviation. For instance we have 

$$QLI(355;5)=\int_0^1 \frac{x \log (x+1) \tan ^{-1}(x)^2}{x^2+1} \, dx,\ \ QLI(145;4)=\int_0^1 \frac{\log (1-x) \log \left(x^2+1\right) \tan ^{-1}(x)}{x^2+1} \, dx$$

$$QPLI4(7;3)=\int_0^1 \frac{\text{Li}_3(-x^2)}{x+1}dx, QPLI4(8;2;1)=\int_0^1 \frac{\text{Ti}_2(x) \log(x)}{1-x}dx$$\\

Weight: We define the weight $W$ of QLI as $1+\sum _{m=1}^5 A(m)$. $QPLI3(\cdots), QPLI4(\cdots)$ are of weight 3, 4 respectively. All 4 examples above are of weight 4.\\

For $W\leq 4$, after excluding those already contained in LIs, there are 193 convergent QLIs. For $W\leq 4$, after excluding those already contained in LI/PLI/QLIs, there are 172 convergent QPLIs.\\

\subsubsection{ES and QES}
Notations:

$$H_n^{(k)}=\sum _{j=1}^n \frac{1}{j^k}, \widetilde{H_n^{(k)}}=\sum _{j=1}^n \frac{(-1)^{j-1}}{j^k}, H_n=H_n^{(1)}, \widetilde{H_n}=\widetilde{H_n^{(1)}}$$

$$f_{1k}(n)=H_n^{(k)}, f_{2k}(n)=\widetilde{H_n^{(k)}},  f_{3k}(n)=H_{2n}^{(k)}$$\\

Note that $\widetilde{H_{2n}^{(k)}}$ won't appear since it equals $H_{2n}^{(k)}-\frac{1}{2^{k-1}}H_{n}^{(k)}$.\\

Definition:

$$ES(a(1),\cdots,a(M),-b(1),\cdots,-b(N);\pm p)= \sum _{n=1}^{\infty } \frac{(\pm 1)^{n-1}  }{n^p}\prod _{k=1}^M H_n^{(a(k))} \prod _{j=1}^N\widetilde{H_n^{(b(j))}}$$

$$QES1(m_1 k_1, \cdots, m_s k_s; \pm p)=\sum_{n=1}^\infty \frac{(\pm1)^{n-1}}{n^p} \prod_{r=1}^s f_{m_r k_r}(n)$$

$$QES2(m_1 k_1, \cdots, m_s k_s; \pm p)=\sum_{n=1}^\infty \frac{(\pm1)^{n-1}}{(2n+1)^p} \prod_{r=1}^s f_{m_r k_r}(n)$$

We may assume $a(i), b(j), m_r k_r$ are non-decreasing. If there exists no $H_n^{(k)}$ (resp. $f_{2k_r}(n)$) and the outer index $+p$, we call an ES (resp. QES) non-alternating, otherwise alternating. \\

Weight: We define the weight $W$ of ES as $\sum _{k=1}^M a(k)+\sum _{j=1}^N b(j)+p$ and QES $\sum_{r=1}^s k_r +p$. For instance 
$$ES(2,-3,-3;-4)=\sum _{n=1}^{\infty } \frac{(-1)^{n-1} H_n^{(2)} \left(\widetilde{H_n^{(3)}}\right)^2}{n^4},  QES(11,12,12,32;3)=\sum _{n=1}^{\infty } \frac{H_n \left(H_n^{(2)}\right){}^2 H_{2 n}^{(2)}}{n^3}$$
are weight 12 alternating ES and weight 10 non-alternating QES, respectively.\\
 
For $W\leq 5$ there are 89 convergent ESs. For $W\leq 4$, after excluding those already contained in ESs, there are 93 convergent QESs. For $W=5$, after excluding those already contained in ESs, there are 28 convergent non-alternating QESs.\\

\subsubsection{Fibonacci basis}
This part is about the structure of closed-forms.\\

Weight of constants: We assign weight 0 to rational numbers, 1 to $\pi, \log(2)$, 2 to Catalan $C$, 4 to $\psi^{(3)}(\frac{1}{4})-\psi^{(3)}(\frac{3}{4})$, $n$ to $\zeta(n),\text{Li}_n(\frac{1}{2}), \Re/\Im \text{Li}_n(1+i)$. The weight given to product of these constants is defined as the sum of each components' weight, for instance $2\pi ^2 \text{Li}_4\left(\frac{1}{2}\right) \zeta (3) \log (2)$ has weight 10.\\

Linear relations: A linear relation is an equality of form $$\text{'rational combination of LI/PLI/QLI/QPLI/ES=certain constants'}$$
If all items in LHS have the same weight, we call it a homogeneous relation; otherwise we call it non-homogeneous.\\

Fibonacci basis: According to \cite{ref1}, all 85 LIs with weight $W\leq5$ are generated by certain $t(W)$ weight $W$ constants over $\mathbb{Q}$, here $t(1)=1, t(2)=2, t(3)=3, t(4)=5, t(5)=8$. We define the sets of these homogeneous constants generating LIs as Fibonacci bases due to consistency between $t(n)$ and Fibonacci numbers $F_n$. This term arises from \cite{ref6} but differs a little. See also \cite{ref1} for higher weight LIs.\\

We denote $A_W$ the weight $W$ Fibonacci basis containing $t(W)$ full-simplified polylog constants. By full-simplified we mean $\zeta (2 n) ,n>0$ and $\text{Li}_n\left(\frac{1}{2}\right) ,n<4$ won't appear, due to simplification formulas \cite{ref13, ref14} $$\text{Li}_2\left(\frac{1}{2}\right)=\frac{\pi ^2}{12}-\frac{\log ^2(2)}{2}, \text{Li}_3\left(\frac{1}{2}\right)=\frac{7 \zeta (3)}{8}+\frac{\log ^3(2)}{6}-\frac{1}{12} \pi ^2 \log (2), \zeta (2 n)=\frac{2^{2 n-1} \pi ^{2 n} \left| B_{2 n}\right| }{(2 n)!}$$  where $B_n$ are Bernoulli numbers. These simplifications keeps the weight unchanged. We have \cite{ref1}

$$A_1=\left\{\log (2)\right\},A_2=\left\{\pi ^2,\log ^2(2)\right\},A_3=\left\{\pi ^2 \log (2),\log ^3(2),\zeta (3)\right\}$$ 
$$A_4=\left\{\pi ^4,\pi ^2 \log ^2(2),\log ^4(2),\zeta (3) \log (2),\text{Li}_4\left(\frac{1}{2}\right)\right\}$$
$$A_5=\left\{\pi ^4 \log (2),\pi ^2 \log ^3(2),\log ^5(2),\pi ^2 \zeta (3),\zeta (3) \log ^2(2),\zeta (5),\text{Li}_4\left(\frac{1}{2}\right) \log (2),\text{Li}_5\left(\frac{1}{2}\right)\right\}$$\\
And that all LIs with $W\leq 5$ are generated by $A_W$ over $\mathbb{Q}$. A possible weight 6 representation is:

\begin{equation}
\begin{aligned}
A_6=\left\{\pi ^6,\pi ^4 \log ^2(2),\pi ^2 \log ^4(2),\log ^6(2),\pi ^2 \zeta (3) \log (2),\zeta (3) \log ^3(2),\zeta (5) \log (2),\right.\\ \left.\zeta (3)^2,\text{Li}_4\left(\frac{1}{2}\right) \log ^2(2),
\pi ^2 \text{Li}_4\left(\frac{1}{2}\right), \text{Li}_5\left(\frac{1}{2}\right) \log (2),\text{Li}_6\left(\frac{1}{2}\right),S\right\}\ \ \ \ \ \ \ \ \
\end{aligned}
\notag
\end{equation}\\
Where the constant $S=ES(1;-5)=\sum _{k=1}^{\infty } \frac{(-1)^{k-1} }{k^5}H_k$ is very likely to be irreducible. Higher weights' representations are omitted.\\

Extended basis: We define extended Fibonacci bases for $W\leq 4$ as\\
$$B_2=A_2 \cup \{C,\pi  \log (2)\}$$

$$B_3=A_3 \cup \left\{C \log (2),\pi  C,\pi  \log ^2(2), \pi ^3,\Im(\text{Li}_3(1+i))\right\}$$

\begin{equation}
\begin{aligned}
B_4=A_4\cup \left\{C^2, \pi C\log(2), C \log ^2(2),\pi ^2 C,\pi  \log ^3(2),\pi ^3 \log (2),\log (2) \Im(\text{Li}_3(1+i)),\right. \\ \left.\pi  \zeta (3), \pi \Im(\text{Li}_3(1+i)),\Im(\text{Li}_4(1+i)),\psi ^{(3)}\left(\frac{1}{4}\right)-\psi ^{(3)}\left(\frac{3}{4}\right)\right\}\ \ \ \ \ \ \ \ \ \ \ \ \ \ \
\end{aligned}
\notag
\end{equation}\\
Note that in these 3 cases each $B_W$ contains $2^W$ constants.\\

\subsection{Main results}

In section $3\sim6$ the following will be proved:\\

\noindent \textbf{Theorem 1.} All 89 ESs and 263 PLIs with weight $W\leq5$ are generated by $A_W$ over $\mathbb{Q}$.\\

\noindent \textbf{Theorem 2.} All 193 QLIs and 172 QPLIs with weight $W\leq 4$ are generated by $B_W$ over $\mathbb{Q}$.\\

Theorem 1, a natural consequence of MZV theory (see for instance \cite{ref18}), will be proved by a systematic non-MZV approach. Theorem 2, intimately connected with level 4 MZVs, will be established by generalizing the approach above.\\

See tables of section 9 for all their closed-forms and 8, 10 for applications.\\

\section{LI}

LIs with weight $W\leq 5$ are completely and explicitly solved in \cite{ref1}. We summarize and add some new results, generating linear relations between LIs.\\

\subsection{Brute force}

It is known that $\int \frac{\log (a+x) \log (b+x)}{c+x} \, dx$, $\int \frac{\log ^2(a+x) \log (b+x)}{a+x} \, dx$, $\int \frac{\log (a+x) \log ^2(b+x)}{a+x} \, dx$, $\int \frac{\log^n(a+x)}{b+x}$ have polylog primitives. By Mathematica we have:\\

{\small $$\text{F1:} \int \frac{\log ^3(a+x)}{b+x} \, dx=6 \text{Li}_4\left(\frac{a+x}{a-b}\right)+3 \log ^2(a+x) \text{Li}_2\left(\frac{a+x}{a-b}\right)$$ $$-6 \log (a+x) \text{Li}_3\left(\frac{a+x}{a-b}\right)+\log ^3(a+x) \log \left(-\frac{b+x}{a-b}\right)$$}

{\small $$\text{F2:} \int \frac{3 \log ^2(a+x) \log (b+x)}{a+x} \, dx=-6 \text{Li}_4\left(\frac{a+x}{a-b}\right)-3 \log ^2(a+x) \text{Li}_2\left(\frac{a+x}{a-b}\right)$$ $$+6 \log (a+x) \text{Li}_3\left(\frac{a+x}{a-b}\right)-\log ^3(a+x) \log \left(-\frac{b+x}{a-b}\right)+\log ^3(a+x) \log (b+x)$$}

{\scriptsize $$\text{F3:} \int \frac{\log (a+x) \log (b+x)}{c+x} \, dx=\text{Li}_3\left(\frac{(a-c) (b+x)}{(b-c) (a+x)}\right)+\left(\text{Li}_2\left(\frac{b+x}{a+x}\right)-\text{Li}_2\left(\frac{(a-c) (b+x)}{(b-c) (a+x)}\right)\right) \log \left(\frac{(a-c) (b+x)}{(a+x) (b-c)}\right)$$ $$+\text{Li}_2\left(\frac{a+x}{a-c}\right) \left(\log (b+x)-\log \left(\frac{(a-c) (b+x)}{(a+x) (b-c)}\right)\right)+\text{Li}_2\left(\frac{b+x}{b-c}\right) \left(\log \left(\frac{(a-c) (b+x)}{(a+x) (b-c)}\right)+\log (a+x)\right)$$ $$+\frac{1}{2} \left(-\log \left(\frac{(a-b) (c+x)}{(a+x) (c-b)}\right)+\log \left(\frac{a-b}{a+x}\right)+\log \left(\frac{c+x}{c-a}\right)\right) \log ^2\left(\frac{(a-c) (b+x)}{(a+x) (b-c)}\right)-\text{Li}_3\left(\frac{b+x}{a+x}\right)$$ $$-\text{Li}_3\left(\frac{a+x}{a-c}\right)-\text{Li}_3\left(\frac{b+x}{b-c}\right)+\log \left(\frac{b+x}{b-c}\right) \left(\log \left(-\frac{c+x}{b-c}\right)-\log \left(\frac{c+x}{c-a}\right)\right) \log \left(\frac{(a-c) (b+x)}{(a+x) (b-c)}\right)$$ $$+\frac{1}{2} \log \left(\frac{b+x}{b-c}\right) \left(\log \left(\frac{b+x}{b-c}\right)-2 \log (a+x)\right) \left(\log \left(\frac{c+x}{c-a}\right)-\log \left(-\frac{c+x}{b-c}\right)\right)+\log (a+x) \log (b+x) \log \left(\frac{c+x}{c-a}\right)$$}

{\tiny$$\text{ F4:} \int \frac{12 \log ^2(a+x) \log (b+x)}{b+x} \, dx=-24 \left(-\text{Li}_4\left(\frac{a+x}{a-b}\right)+\text{Li}_4\left(-\frac{b+x}{a-b}\right)+\text{Li}_4\left(\frac{b+x}{a+x}\right)\right)-12 \text{Li}_2\left(\frac{b+x}{a+x}\right) \log ^2\left(-\frac{b+x}{a+x}\right)$$ $$+12 \text{Li}_2\left(\frac{a+x}{a-b}\right) \left(\log \left(-\frac{b+x}{a-b}\right)-\log \left(-\frac{b+x}{a+x}\right)\right)^2+24 \text{Li}_2\left(\frac{a+x}{a-b}\right) \log \left(\frac{a+x}{a-b}\right) \left(\log (b+x)-\log \left(-\frac{b+x}{a-b}\right)\right)$$ $$+12 \text{Li}_2\left(-\frac{b+x}{a-b}\right) \left(2 \log \left(\frac{a+x}{a-b}\right) \left(\log (b+x)-\log \left(-\frac{b+x}{a-b}\right)\right)+\log \left(-\frac{b+x}{a-b}\right) \left(\log \left(-\frac{b+x}{a-b}\right)-2 \log \left(-\frac{b+x}{a+x}\right)\right)-2 \log (a+x) \log (b+x)\right)$$ $$+24 \text{Li}_3\left(\frac{a+x}{a-b}\right) \left(\log \left(-\frac{b+x}{a+x}\right)-\log (b+x)\right)+24 \text{Li}_3\left(-\frac{b+x}{a-b}\right) \left(\log \left(-\frac{b+x}{a+x}\right)+\log (a+x)\right)+24 \text{Li}_3\left(\frac{b+x}{a+x}\right) \log \left(-\frac{b+x}{a+x}\right)$$ $$+\log ^4\left(-\frac{b+x}{a-b}\right)+\log ^4\left(-\frac{b+x}{a+x}\right)-4 \left(\log \left(-\frac{b+x}{a-b}\right)+\log \left(\frac{a-b}{a+x}\right)\right) \log ^3\left(-\frac{b+x}{a+x}\right)+6 \log ^2\left(-\frac{b+x}{a+x}\right) \log ^2\left(-\frac{b+x}{a-b}\right)$$ $$-4 \left(3 \log \left(\frac{a+x}{a-b}\right)+\log \left(-\frac{b+x}{a-b}\right)\right) \log \left(-\frac{b+x}{a+x}\right) \log ^2\left(-\frac{b+x}{a-b}\right)+6 \log ^2(a+x) \log ^2(b+x)$$ $$+6 \log ^2\left(\frac{a+x}{a-b}\right) \left(\log (b+x)-\log \left(-\frac{b+x}{a-b}\right)\right) \left(3 \log \left(-\frac{b+x}{a-b}\right)+\log (b+x)\right)-4 \log \left(\frac{a+x}{a-b}\right) \left(3 \log (a+x) \log ^2(b+x)-2 \log ^3\left(-\frac{b+x}{a-b}\right)\right)$$}

Here logs and polylogs are principal, following the convention of Mathematica. Based on this, for each LI we may check by Mathematica whether its primitive is computable. If it is, we get the closed-form directly via Newton-Leibniz formula.\\

\subsection{General formulas}

\noindent \textbf{Proposition 1.} The following formulas hold:\\
\begin{equation}
\begin{aligned}
\int_0^1 x^m \log^n(x) \, dx=\frac{(-1)^n n!}{(m+1)^{n+1}}
\end{aligned}
\tag{1}
\end{equation}

\begin{equation}
\begin{aligned}
LI(0,0,n;2)=\frac{\log ^{n+1}(2)}{n+1}
\end{aligned}
\tag{2}
\end{equation}

\begin{equation}
\begin{aligned}
LI(0,n,0;0)=LI(n,0,0;1)=(-1)^n n! \zeta (n+1)
\end{aligned}
\tag{3}
\end{equation}

\begin{equation}
\begin{aligned}
LI(0,n,0;2)=\left(1-2^{-n}\right) (-1)^n n! \zeta (n+1)
\end{aligned}
\tag{4}
\end{equation}

\begin{equation}
\begin{aligned}
LI(n,0,0;2)=(-1)^n n! \text{Li}_{n+1}\left(\frac{1}{2}\right)
\end{aligned}
\tag{5}
\end{equation}

\begin{equation}
\begin{aligned}
LI(1,n,0;1)=LI(n,1,0,0)=(-1)^{n-1} n! \zeta (n+2)
\end{aligned}
\tag{6}
\end{equation}

\begin{equation}
\begin{aligned}
LI(0,n,1;1)=\left(1-2^{-(n+1)}\right) (-1)^n n! \zeta (n+2)
\end{aligned}
\tag{7}
\end{equation}

\begin{equation}
\begin{aligned}
LI(0,0,n;1)=-n! \sum _{k=0}^{n-1} \frac{\log ^k(2) \text{Li}_{-k+n+1}\left(\frac{1}{2}\right)}{k!}+n! \zeta (n+1)-\frac{n \log ^{n+1}(2)}{n+1}
\end{aligned}
\tag{8}
\end{equation}

\begin{equation}
\begin{aligned}
LI(0,1,n;2)=n! \sum _{k=0}^n \frac{\log ^k(2) \text{Li}_{-k+n+2}\left(\frac{1}{2}\right)}{k!}-n! \zeta (n+2)+\frac{\log ^{n+2}(2)}{n+2}
\end{aligned}
\tag{9}
\end{equation}

\begin{equation}
\begin{aligned}
LI(1,0,n;2)=n! \sum _{k=2}^{n+2} \frac{(-1)^{k-1} \zeta (k) \log ^{-k+n+2}(2)}{(-k+n+2)!}+(-1)^n n! \text{Li}_{n+2}\left(\frac{1}{2}\right)+\frac{\log ^{n+2}(2)}{n+1}
\end{aligned}
\tag{10}
\end{equation}

\begin{equation}
\begin{aligned}
LI(1,n,0;0)=LI(n,1,0;1)=\frac{1}{2} (-1)^{n-1} n! \left((n+1) \zeta (n+2)-\sum _{k=2}^n \zeta (k) \zeta (-k+n+2)\right)
\end{aligned}
\tag{11}
\end{equation}

\begin{equation}
\begin{aligned}
LI(0,n,1;0)= \frac{1}{2}(-1)^n n!  \sum _{k=2}^n \left(1-2^{1-k}\right) \left(1-2^{k-n-1}\right) \zeta (k) \zeta (-k+n+2) \\
-(-1)^n n! \left(\frac{1}{2} (n+1) \zeta (n+2)-2 \left(1-2^{-(n+1)}\right) \log (2) \zeta (n+1)\right)\ \ \ \ \
\end{aligned}
\tag{12}
\end{equation}

\begin{equation}
\begin{aligned}
LI(0,2n-1,1;2)=-(2 n-1)! \sum _{k=1}^{n-1} \left(1-2^{1-2 k}\right) \zeta (2 k) \zeta (-2 k+2 n+1) \\
-(2 n-1)!\left(2^{-(2 n+1)} \zeta (2 n+1)-(n-1) \left(1-2^{-2 n}\right) \zeta (2 n+1)\right)\ \ \ \ 
\end{aligned}
\tag{13}
\end{equation}

\begin{equation}
\begin{aligned}
LI(1,2n-1,0;2)=-(2 n-1)! \sum _{k=1}^{n-1} \left(1-2^{2 k-2 n}\right) \zeta (2 k) \zeta (-2 k+2 n+1)\ \ \ \ \ \ \ \  \\
+(2 n-1)!\left(\left((n+1) \left(1-2^{-2 n}\right)+2^{-(2 n+1)}\right) \zeta (2 n+1)-2 \left(1-2^{-2 n}\right) \log (2) \zeta (2 n)\right)
\end{aligned}
\tag{14}
\end{equation}

\begin{equation}
\begin{aligned}
LI(2,n,0;1)=(-1)^n n! \left((n+2) \zeta (n+3)-\sum _{k=2}^{n+1} \zeta (k) \zeta (-k+n+3)\right)
\end{aligned}
\tag{15}
\end{equation}

\begin{equation}
\begin{aligned}
LI(0,2n,2;1)=2 (2 n)! \sum _{k=1}^n \left(1-2^{1-2 k}\right) \zeta (2 k) \zeta (-2 k+2 n+3)\\
-2 (2 n)! \left(\left(1-2^{-(2 n+2)}\right) n \zeta (2 n+3)-2^{-(2 n+3)} \zeta (2 n+3)\right)\ \
\end{aligned}
\tag{16}
\end{equation}

\begin{equation}
\begin{aligned}
LI(1,2n,1;1)=-(2 n)! 2^{-(2 n+1)} \sum _{k=1}^n \left(1-2^{-2 k+2 n+2}\right) \zeta (2 k) \zeta (-2 k+2 n+3)\\
+(2 n)!\left((2 n+3) 2^{-(2 n+3)}-1\right) \zeta (2 n+3)\ \ \ \ \ \ \ \ \ \ \ \ \ \ \ \ \ \ \ \ \
\end{aligned}
\tag{17}
\end{equation}\\

\noindent Proof. (1): $x\to e^{-t}, \Gamma$. (2): Trivial. (3)-(5):  Series expansion of $\frac{1}{1-ax}$, then use (1). (6)-(7): Series expansion of $\log(1-ax)$. (8): Nielsen-Ramanujan integral. $t\to \frac{1}{1+x}$, then use $\int_{\frac{1}{2}}^1 x^m \log ^n(x) \, dx=\frac{\partial ^n}{\partial m^n}\frac{1-2^{-(m+1)}}{m+1}$. (9)-(10): Use (8), integrate by parts. (11)-(12): Series expansion of $\frac{\log(1-x)}{1-x}$, use (1) then ES results in \cite{ref11}. (13)-(17) Series expansion and \cite{ref11}. $\hfill\square$\\

For each LI, check whether they are directly solved using general formulas above.\\

Clarification: Throughout this paper, legitimacy of all operations of exchanging orders of multiple summation and integration can be justified by convergence theorems.\\

\subsection{Integration by parts}

Take $F(x)=\prod_{m=0}^2 l(m;x)^{a(m)}$ in N-L formula $\int_0^1 F'(x) \, dx=F(1)-F(0)$ to obtain linear a linear relation between LIs. Now consider all possible $F$. As one can see, this is essentially the process of partial integration.\\

Clarification: Though we don't carry out partial integration from the view of manipulating single LIs here, this will be used frequently later. Throughout this paper, by 'integration by parts lifting (up) $u'$ in integral $\int u' v$, we mean transforming it into $uv-\int u v'$.\\

\subsection{Fractional transformation}

Take $F(x)=\frac{\prod _{m=0}^2 l(m;x)^{a(m)}}{f(p;x)} (p\not=1)$ in the identity $\int_0^1 F(x) \, dx=\int_0^1 \frac{2 F\left(\frac{1-x}{x+1}\right)}{(x+1)^2} \, dx$. Since $$\log (x)\to \log (1-x)-\log (x+1),\log (1-x)\to \log (x)-\log (x+1)+\log (2),\log (x+1)\to \log (2)-\log (x+1)$$ under the substitution $x\to \frac{1-x}{x+1}$, we may obtain a linear relation after expansion and possible partial fraction decomposition on $\frac{1}{x(x+1)}$ for $p=0$. Now consider all possible $F$.\\

\subsection{Beta derivatives}

\noindent \textbf{Proposition 2.} The following formulas hold:\\
\begin{equation}
\begin{aligned}
LI(n,m,0;1)=LI(m,n,0;0)=\underset{\{a,b\}\to \{0,1\}}{\text{lim}}\frac{\partial ^{m+n}B(a,b)}{\partial a^m\, \partial b^n}
\end{aligned}
\tag{1}
\end{equation}

\begin{equation}
\begin{aligned}
\sum _{k=0}^n \binom{n}{k}LI(k,m,n-k;1)=2^{-(m+1)} \underset{\{a,b\}\to \{0,1\}}{\text{lim}}\frac{\partial ^{m+n}B(a,b)}{\partial a^m\, \partial b^n}
\end{aligned}
\tag{2}
\end{equation}

\begin{equation}
\begin{aligned}
\sum _{k=0}^n \binom{n}{k}LI(k,m,n-k;0)=2^{-(m+1)} \left(\underset{\{a,b\}\to \{1,0\}}{\text{lim}}\frac{\partial ^{m+n}B(a,b)}{\partial a^m\, \partial b^n}+\underset{\{a,b\}\to \left\{\frac{1}{2},0\right\}}{\text{lim}}\frac{\partial ^{m+n}B(a,b)}{\partial a^m\, \partial b^n}\right)
\end{aligned}
\tag{3}
\end{equation}

\begin{equation}
\begin{aligned}
\sum _{k=0}^n \binom{n}{k}LI(k,m,n-k;2)=2^{-(m+1)} \left(\underset{\{a,b\}\to \{\frac{1}{2},0\}}{\text{lim}}\frac{\partial ^{m+n}B(a,b)}{\partial a^m\, \partial b^n}-\underset{\{a,b\}\to \left\{1,0\right\}}{\text{lim}}\frac{\partial ^{m+n}B(a,b)}{\partial a^m\, \partial b^n}\right)
\end{aligned}
\tag{4}
\end{equation}

\begin{equation}
\begin{aligned}
LI(0,m,1;0)=\left(2^{-(m+1)}-1\right) \underset{\{a,b\}\to \{1,0\}}{\text{lim}}\frac{\partial ^{m+1}B(a,b)}{\partial a^m\, \partial b^1}+2^{-(m+1)} \underset{\{a,b\}\to \left\{\frac{1}{2},0\right\}}{\text{lim}}\frac{\partial ^{m+1}B(a,b)}{\partial a^m\, \partial b^1}
\end{aligned}
\tag{5}
\end{equation}

\begin{equation}
\begin{aligned}
\sum _{k=0}^n (-1)^k \binom{n}{k} LI(k,n-k,0;2)=\underset{s\to 1}{\text{lim}}\frac{\partial ^n}{\partial s^n}\frac{\left(2^{1-s}-1\right) \pi }{\sin (\pi  s)}
\end{aligned}
\tag{6}
\end{equation}

{\footnotesize\begin{equation}
\begin{aligned}
\sum _{j=0}^n \sum _{k=0}^{n+1} \binom{n}{j} \binom{n+1}{k} (-\log (2))^{j+k} LI(-k+n+1,0,n-j;2)=\frac{1}{2}\underset{\{x,y\}\to \left\{0,1\right\}} {\text{lim}}\frac{\partial ^{2n+1}B(a,b)}{\partial a^n\, \partial b^{n+1}}-\frac{\log ^{2 n+2}(2)}{2 n+2} 
\end{aligned}
\tag{7}
\end{equation}}

\begin{equation}
\begin{aligned}
\sum _{k=0}^n (-1)^k  \binom{n}{k} LI(0,n-k,k;2)=(-1)^n\sum _{l=0}^n  l! \binom{n}{l}\text{Li}_{l+1}\left(\frac{1}{2}\right)  \log ^{n-l}(2)
\end{aligned}
\tag{8}
\end{equation}

\begin{equation}
\begin{aligned}
(-1)^k \sum _{j=0}^m (-1)^j \binom{m}{j} LI(0,j+k,-j+m+n;1)+\sum _{j=0}^n (-1)^j \binom{n}{j} LI(0,j+k,-j+m+n;1)\\=(-1)^{m+n} \sum _{j=0}^k (-1)^k \binom{k}{j} \left(\underset{\{x,y\}\to \left\{0,1\right\}} {\text{lim}}\frac{\partial ^{k+m+n}B(a,b)}{\partial a^{j+m}\, \partial b^{-j+k+n}}+\underset{\{x,y\}\to \left\{1,0\right\}} {\text{lim}}\frac{\partial ^{k+m+n}B(a,b)}{\partial a^{j+m}\, \partial b^{-j+k+n}}\right)\ \
\end{aligned}
\tag{9}
\end{equation}
\\

\noindent Proof. (1): Differentiate $B(a,b)$. (2)-(4): consider $\int_0^1 \frac{\log ^m(x) \log ^n\left(1-x^2\right)}{f(p;x)} \, dx$. Using $\frac{1}{1-x}=\frac{x+1}{1-x^2},\frac{1}{x+1}=\frac{1-x}{1-x^2}$ and that $x^2\to t$ give RHS. Alternatively $\log(1-x^2)=\log(1-x)+\log(1+x)$ and binomial theorem for expansion give LHS.  (5): Multiply both sides of $\frac{\log (x+1)}{1-x}=\frac{x \log \left(1-x^2\right)}{1-x^2}+\frac{\log \left(1-x^2\right)}{1-x^2}-\frac{\log (1-x)}{1-x}$ with $\log^m(x)$ then integrate, $t\to x^2$. (6): Consider $F(x)=\frac{\log ^n\left(\frac{x}{1-x}\right)}{x+1}$, substituting $t\to \frac{x}{1-x}$ and differentiating $\int_0^\infty \frac{x^{s-1}}{(x+1) (2 x+1)} \, dx=\frac{\pi  \left(2^{1-s}-1\right)}{\sin (\pi  s)}$ yield RHS, and binomial theorem yields LHS. (7): $\int_0^{\frac{1}{2}} \frac{\log ^n(1-x) \log ^{n+1}(x)}{1-x} \, dx$ is evaluated in \cite{ref15}. Let $x\to \frac{1-u}{2}$, use binomial theorem twice. (8): Similar to (6), subsitute $t\to \frac{x}{x+1}$ in integral of $F_1(x)=\frac{\log ^n\left(\frac{x}{x+1}\right)}{x+1}$, series expansion. (9): Consider $\int_0^{\infty } \frac{\log ^k(x) \log ^m(x+1) \log ^n\left(\frac{1}{x}+1\right)}{x} \, dx, m, n>0$. Substitution $x\to\frac{1}{t}$ for $(1,\infty)$ part and binomial theorem yields LHS, while $x\to \frac{1-t}{t}$ for the whole integral and partial fractions give RHS. $\hfill\square$\\

Partial derivatives in above expressions can be evaluated by Mathematica (or manually by computing Laurent series of polygamma and plugging in their special values). Now consider all possible relations generated from these formulas.\\

\subsection{Contour integration}

Integrate $$F(z)=\frac{\log^k(1-z) \log^m(z) \log^n(z+1)}{z (1-z)}$$ (with all logs principal) along the large upper semicircular contour, by Cauchy theorem we have $\int_{-\infty+0i}^{\infty+0i} F(z)dz=0$. Split the integral into 4 parts by branch points $0, \pm 1$, use $x\to\pm t, x\to \pm \frac{1}{t}$ to map them to $(0,1)$, after expansion we obtain a linear relation (see \cite{ref1} and sections below for details and examples). Now consider all possible $F$.\\

Clarification: Here and below, by 'principal' we follow the convention of Mathematica, that is, the argument of log/polylogs vary in $(-\pi, \pi)$.\\

\subsection{Double integration}

In \cite{ref1} a symmetric double log integral is investigated, which essentially yields:\\
\begin{equation}
\begin{aligned}
LI(2n,0,1;1)=(2 n)! \left(\sum _{k=0}^{n-1} (-1)^k \text{Li}_{k+1}\left(\frac{1}{2}\right) \text{Li}_{-k+2 n+1}\left(\frac{1}{2}\right)+\frac{1}{2} (-1)^n \text{Li}_{n+1} \left(\frac{1}{2}\right)^2+\text{Li}_{2 n+2}\left(\frac{1}{2}\right)\right)
\end{aligned}
\notag
\end{equation}\\

\subsection{Hypergeometric identities}

Differentiating hypergeometric identities can be used to generate the most nontrivial high weight relations, see \cite{ref1}.\\

\subsection{Solution}

All linear relations involving LIs with $W\leq5$ generated by subsection 2-1$\sim$2-6 already yields:\\

\noindent \textbf {Lemma 1.} All LIs with weight $W\leq5$ can be generated by $A_W$ over $\mathbb{Q}$.\\

In fact $W=6$ case is also solvable using 2-1$\sim$2-6 only, but we won't record their closed-forms below due to great complexity (for their values we refer readers to \cite{ref1}). $W=6$ results will be frequently used later so it's better to be aware of the general structure.\\

Clarification: In section $N$, all subjects computed in section $1,\cdots, N-1$ are supposed to be known. As an example, LI values are frequently used below with the statement '$\cdots$ modulo LIs'.\\

\section{ES}
Since the theory of Euler sums is well developed already, this section won't be a thorough summary, but a simple review of ES evaluation and solution of several ESs.\\

Excluding the trivial $\sum _{n=1}^{\infty } \frac{\left(\pm1\right)^{n-1}}{n^p}$, there lefts 89 ESs with weight $W\leq5$ (2, 9, 24, 54 for weight 2, 3, 4, 5 respectively) .\\

\subsection{General formulas}

According to \cite{ref11}, we have\\

\noindent \textbf {Lemma 2.} $ES(\pm1;p), ES(\pm1;-p)$ ($p$ even), $ES(p;p), ES(p;q)$ ($p+q$ odd or smaller than 7), $ES(1,1;p)$ ($p$ odd or smaller than 5) are reducible to Zeta values.\\

See \cite{ref2, ref7} for alternative proofs.\\

\subsection{Symmetric relations}

By inclusion-exclusion (abbr. I-E) principle we have:\\
$$\sum _{0<j,k,l}^{}=\sum _{0<j,k\leq l}^{}+\sum _{0<j,l\leq k}^{}+\sum _{0<k,l\leq j}^{}-\sum _{0<j\leq k=l}^{}-\sum _{0<k\leq j=l}^{}-\sum _{0<l\leq j=k}^{}+\sum _{0<j=k=l}^{}$$

Apply this identity as an operator on all triple summands $\frac{(\pm 1)^{j-1} (\pm 1)^{k-1} (\pm 1)^{l-1}}{j^a k^b l^c}$ to obtain linear relations. For instance, if we choose 3 '$\pm$'s to be '+' and assuming convergence one obtain\\
{\small $$\zeta (a) \zeta (b) \zeta (c)=ES(a,b;c)+ES(a,c;b)+ES(b,c;a)-ES(a;b+c)-ES(b;a+c)-ES(c;a+b)+ES(a+b+c)$$}

Consider all other choices and $n$-ple generalizations to obtain relations, see for instance \cite{ref21, ref22}. See also \cite {ref9} for similar contents on MZVs.\\

\subsection{Partial fractions}

Consider Cauchy product of the double sum $\sum _{j=1}^{\infty } \sum _{k=1}^{\infty } \frac{1}{j^a k^b}$. By partial fraction decomposition we know that all sums of the form $\sum _{k=1}^{n-1} \frac{1}{k^a (n-k)^b}$ are expressible in a combination of homogeneous $\frac{H_n^{(k)}}{n^{a+b-k}}, k=1,\cdots,a+b-1$, we found that $\zeta(a) \zeta(b)$ is expressed in a corresponding combination of ESs, which is a linear relation. We may also consider $\frac{(\pm 1)^{j-1} (\pm 1)^{k-1}}{j^a k^b}$ and generalizations to triple sums. See \cite{ref9} for more.\\

\subsection{Contour integration}

See \cite{ref11} again for generating ES relations using contour integration and polygamma kernels.\\

\subsection{Solution}
\subsubsection{Preparation}
In \cite{ref11, ref19, ref21, ref22} all but 1 weight 3 ES, 2 weight 4 ESs and 7 weight 5 ESs are solved. 3 weight $<5$ ESs are evaluated trivially. Now we will sketch a proof for the remaining 7, that is\\
$$ES(-1,-1,-1;-2), ES(-1,-1,-1,-1;-1), ES(1,1,-2;-1), ES(1,1,-1,-1;-1)$$
$$ES(1,-1,-1;-2), ES(1,-1,-2;-1), ES(1,-1,-1,-1;-1)$$\\

All LIs, low weight ESs and other 47 weight 5 ESs are supposed to be known. Also, several GES (Geometric Euler Sums, i.e. sums involving $\frac{1}{2^n}$ like $\sum _{n=1}^{\infty } \frac{H_n}{n^3 2^n}$) computed in \cite{ref19} will be applied too.\\

\noindent \textbf{Proposition 3.} (1)-(15) hold where $\left| x\right| <1$ and naturally extends to $|z|\leq1$ assuming convergence:\\
\begin{equation}
\begin{aligned}
\frac{1}{1\mp x}=\sum _{n=0}^{\infty } \left(\pm x\right)^n, -\log (1\mp x)=\sum _{n=1}^{\infty } \frac{\left(\pm x\right)^n}{n}
\end{aligned}
\tag{1}
\end{equation}

\begin{equation}
\begin{aligned}
-\frac{\log (1\mp x)}{1\mp x}=\sum _{n=1}^{\infty } H_n \left(\pm x\right)^n, \frac{\log (1\pm x)}{1\mp x}=\sum _{n=1}^{\infty } \widetilde{H_n^{}} \left(\pm x\right)^n
\end{aligned}
\tag{2}
\end{equation}

\begin{equation}
\begin{aligned}
\log ^2(1\pm x)=\sum _{n=1}^{\infty } 2 (-1)^{n-1} \left(\frac{1}{n^2}-\frac{H_n}{n}\right) \left(\pm x\right)^n
\end{aligned}
\tag{3}
\end{equation}

\begin{equation}
\begin{aligned}
\log (1-x) \log (x+1)=-\sum _{n=1}^{\infty } \left(\frac{\widetilde{H_{2 n}^{}}}{n}+\frac{1}{2 n^2}\right) x^{2 n}
\end{aligned}
\tag{4}
\end{equation}

\begin{equation}
\begin{aligned}
\sum _{n=1}^{\infty } \frac{\left(\pm x\right)^n}{n^k}=\text{Li}_k(\pm x),  \sum _{n=1}^{\infty }  H_n^{(k)}\left(\pm x\right)^n=\frac{\text{Li}_k(\pm x)}{1\mp x},  \sum _{n=1}^{\infty }\widetilde{H_n^{(k)}} \left(\pm x\right)^n =-\frac{\text{Li}_k(\mp x)}{1\mp x}
\end{aligned}
\tag{5}
\end{equation}

\begin{equation}
\begin{aligned}
\sum _{n=1}^{\infty } \frac{H_n x^n}{n}=\text{Li}_2(x)+\frac{1}{2} \log ^2(1-x)
\end{aligned}
\tag{6}
\end{equation}

\begin{equation}
\begin{aligned}
\sum _{n=1}^{\infty } \frac{\widetilde{H_{n}^{}} x^n}{n}=\text{Li}_2\left(\frac{1-x}{2}\right)-\text{Li}_2(-x)-\text{Li}_2\left(\frac{1}{2}\right)-\log (2) \log (1-x)
\end{aligned}
\tag{7}
\end{equation}

\begin{equation}
\begin{aligned}
\sum _{n=1}^{\infty } \left(H_n\right){}^2 x^n=\frac{\text{Li}_2(x)}{1-x}+\frac{\log ^2(1-x)}{1-x}
\end{aligned}
\tag{8}
\end{equation}

\begin{equation}
\begin{aligned}
\sum _{n=1}^{\infty }\left(\widetilde{H_{n}^{}}\right)^2 x^n=\frac{2 \left(\text{Li}_2\left(\frac{1-x}{2}\right)-\text{Li}_2\left(\frac{1}{2}\right)+\log \left(\frac{x+1}{2}\right) \log (1-x)\right)+\text{Li}_2(x)}{1-x}
\end{aligned}
\tag{9}
\end{equation}

\begin{equation}
\begin{aligned}
\sum _{n=1}^{\infty } H_n \widetilde{H_{n}^{}} x^n=\frac{\text{Li}_2\left(\frac{1-x}{2}\right)-\text{Li}_2(-x)-\text{Li}_2\left(\frac{1}{2}\right)-\frac{1}{2} \log ^2(x+1)-\log (2) \log (1-x)}{1-x}
\end{aligned}
\tag{10}
\end{equation}

\begin{equation}
\begin{aligned}
\sum _{n=1}^{\infty } H_n H_n^{(2)} x^n=\frac{\text{Li}_3(1-x)+\text{Li}_3(x)+\frac{1}{2} \log (x)\log ^2(1-x) -\frac{1}{6} \pi ^2 \log (1-x)-\zeta (3)}{1-x}
\end{aligned}
\tag{11}
\end{equation}

{\small\begin{equation}
\begin{aligned}
\sum _{n=1}^{\infty } \left(H_n\right){}^3 x^n=\frac{3 \text{Li}_3(1-x)+\text{Li}_3(x)-\log ^3(1-x)+\frac{3}{2} \log (x) \log ^2(1-x)-\frac{1}{2} \pi ^2 \log (1-x)-3 \zeta (3)}{1-x}
\end{aligned}
\tag{12}
\end{equation}}

\begin{equation}
\begin{aligned}
\sum _{n=1}^{\infty } \frac{\left(H_n\right){}^2 x^n}{n}=\text{Li}_3(x)-\text{Li}_2(x) \log (1-x)-\frac{1}{3} \log ^3(1-x)
\end{aligned}
\tag{13}
\end{equation}

{\footnotesize\begin{equation}
\begin{aligned}
\sum _{n=1}^{\infty } \frac{H_n^{(2)} x^n}{n}=2 \text{Li}_3(1-x)+\text{Li}_3(x)-2 \text{Li}_2(1-x) \log (1-x)-\text{Li}_2(x) \log (1-x)-\log (x) \log ^2(1-x)-2 \zeta (3)
\end{aligned}
\tag{14}
\end{equation}}
\begin{equation}
\begin{aligned}
\sum _{n=1}^{\infty } \frac{H_n x^n}{n^2}=-\text{Li}_3(1-x)+\text{Li}_3(x)+\text{Li}_2(1-x) \log (1-x)+\frac{1}{2}\log (x) \log ^2(1-x) +\zeta (3)
\end{aligned}
\tag{15}
\end{equation}\\

\noindent Proof. (1): Trivial. (2)-(5): Cauchy product, use (1). (6)-(7):  Applying $L(f)=\int_0^x \frac{f(t)}{t} \, dt$ on certain $f$ in (5). (8): Denote LHS as $F(x)$, suppose it is of form $F(x)=\frac{G(x)}{1-x}$, then by Cauchy we have $G(x)=\sum _{n=1}^{\infty } \left(\left(H_n\right){}^2-\left(H_{n-1}\right){}^2\right) x^n$. Using $H_{n-1}=H_n-\frac{1}{n}$ and (5), (6) to arrive at RHS. (9)-(12): Calculating differences like (8). (13)-(15): Make use of $L(f)$ on (5), (6), (8). Mathematica is used to calculate primitives. $\hfill\square$\\\\

\noindent \textbf{Proposition 4.} (1)-(10) hold where $\left| x\right| <1$ and naturally extends to $|z|\leq1$ assuming convergence:\\

{\small $$(1)\ \ \sum _{n=1}^{\infty } \frac{\left(H_n\right){}^2 x^n}{n^2}=\frac{\text{Li}_2(x){}^2}{2}-2 \text{Li}_4(1-x)+\text{Li}_4(x)$$ $$-\text{Li}_2(1-x) \log ^2(1-x)+2 \text{Li}_3(1-x) \log (1-x)-\frac{1}{3} \log (x) \log ^3(1-x)+\frac{\pi ^4}{45}$$}

{\scriptsize $$(2)\ \ \sum _{n=1}^{\infty } \frac{H_n^{(2)} x^n}{n^2}=-\frac{1}{6} \log ^4\left(\frac{1}{x}-1\right)+\frac{2}{3} \log (1-x) \log ^3\left(\frac{1}{x}-1\right)+\frac{2}{3} \log \left(\frac{1}{x}\right) \log ^3\left(\frac{1}{x}-1\right)$$ $$-\log ^2(1-x) \log ^2\left(\frac{1}{x}-1\right)+2 \text{Li}_2\left(\frac{x-1}{x}\right) \log ^2\left(\frac{1}{x}-1\right)-2 \text{Li}_2(x) \log ^2\left(\frac{1}{x}-1\right)+\frac{2}{3} \log ^3(1-x) \log \left(\frac{1}{x}-1\right)$$ $$+2 \log ^2(1-x) \log (x) \log \left(\frac{1}{x}-1\right)+4 \log (1-x) \text{Li}_2(1-x) \log \left(\frac{1}{x}-1\right)+4 \log (1-x) \text{Li}_2(x) \log \left(\frac{1}{x}-1\right)$$ $$-4 \text{Li}_3(1-x) \log \left(\frac{1}{x}-1\right)-4 \text{Li}_3\left(\frac{x-1}{x}\right) \log \left(\frac{1}{x}-1\right)-4 \text{Li}_3(x) \log \left(\frac{1}{x}-1\right)+\frac{\pi ^4}{18}-\frac{1}{12} \log ^4\left(\frac{x}{1-x}\right)$$ $$-\frac{1}{6} \log ^4(1-x)+\frac{1}{3} \log \left(\frac{1}{1-x}\right) \log ^3\left(\frac{x}{1-x}\right)+\frac{1}{3} \log \left(\frac{x}{1-x}\right) \log ^3(x)-\frac{2}{3} \log (1-x) \log ^3(x)-\frac{1}{2} \log ^2(x) \log ^2\left(\frac{x}{1-x}\right)$$ $$+\frac{3}{2} \log ^2(1-x) \log ^2(x)+\log \left(\frac{x}{1-x}\right) \log (1-x) \log ^2(x)+\frac{\text{Li}_2(x){}^2}{2}+\frac{1}{3} \log ^3\left(\frac{x}{1-x}\right) \log (x)-\frac{4}{3} \log ^3(1-x) \log (x)$$ $$+\log ^2\left(\frac{x}{1-x}\right) \text{Li}_2\left(-\frac{x}{1-x}\right)-\log ^2\left(\frac{x}{1-x}\right) \text{Li}_2(1-x)-2 \log ^2(1-x) \text{Li}_2(1-x)-\log ^2(x) \text{Li}_2(1-x)$$ $$+2 \log \left(\frac{x}{1-x}\right) \log (x) \text{Li}_2(1-x)+2 \log (1-x) \log (x) \text{Li}_2(1-x)-2 \log ^2(1-x) \text{Li}_2(x)-\log ^2(x) \text{Li}_2(x)$$ $$+2 \log \left(\frac{x}{1-x}\right) \log (x) \text{Li}_2(x)+2 \log (1-x) \log (x) \text{Li}_2(x)-2 \log \left(\frac{x}{1-x}\right) \text{Li}_3\left(-\frac{x}{1-x}\right)-2 \log \left(\frac{x}{1-x}\right) \text{Li}_3(1-x)$$ $$-2 \log \left(\frac{x}{1-x}\right) \text{Li}_3(x)+2 \log (1-x) \text{Li}_3(x)+2 \text{Li}_4\left(-\frac{x}{1-x}\right)+2 \text{Li}_4(1-x)+4 \text{Li}_4\left(\frac{x-1}{x}\right)-\text{Li}_4(x)-2 \log (x) \zeta (3)-\frac{\log ^4(x)}{12}$$\\

$$(3)\ \ \sum _{n=1}^{\infty } \frac{H_n x^n}{n^3}=\frac{1}{12} \log ^4\left(\frac{1}{x}-1\right)-\frac{1}{3} \log (1-x) \log ^3\left(\frac{1}{x}-1\right)-\frac{1}{3} \log \left(\frac{1}{x}\right) \log ^3\left(\frac{1}{x}-1\right)$$ $$+\frac{1}{2} \log ^2(1-x) \log ^2\left(\frac{1}{x}-1\right)-\text{Li}_2\left(\frac{x-1}{x}\right) \log ^2\left(\frac{1}{x}-1\right)+\text{Li}_2(x) \log ^2\left(\frac{1}{x}-1\right)-\frac{1}{3} \log ^3(1-x) \log \left(\frac{1}{x}-1\right)$$ $$-\log ^2(1-x) \log (x) \log \left(\frac{1}{x}-1\right)-2 \log (1-x) \text{Li}_2(1-x) \log \left(\frac{1}{x}-1\right)-2 \log (1-x) \text{Li}_2(x) \log \left(\frac{1}{x}-1\right)$$ $$+2 \text{Li}_3(1-x) \log \left(\frac{1}{x}-1\right)+2 \text{Li}_3\left(\frac{x-1}{x}\right) \log \left(\frac{1}{x}-1\right)+2 \text{Li}_3(x) \log \left(\frac{1}{x}-1\right)+\frac{1}{24} \log ^4\left(\frac{x}{1-x}\right)+\frac{1}{12} \log ^4(1-x)$$ $$+\frac{\log ^4(x)}{24}-\frac{1}{6} \log \left(\frac{1}{1-x}\right) \log ^3\left(\frac{x}{1-x}\right)-\frac{1}{6} \log (x) \log ^3\left(\frac{x}{1-x}\right)+\frac{1}{3} \log (1-x) \log ^3(x)+\frac{1}{4} \log ^2\left(\frac{x}{1-x}\right) \log ^2(x)$$ $$-\frac{3}{4} \log ^2(1-x) \log ^2(x)-\frac{1}{6} \log ^3(x) \log \left(\frac{x}{1-x}\right)-\frac{1}{2} \log (1-x) \log ^2(x) \log \left(\frac{x}{1-x}\right)+\frac{2}{3} \log ^3(1-x) \log (x)$$ $$+\frac{1}{2} \log ^2\left(\frac{x}{1-x}\right) \text{Li}_2(1-x)+\log ^2(1-x) \text{Li}_2(1-x)+\frac{1}{2} \log ^2(x) \text{Li}_2(1-x)-\log (x) \log \left(\frac{x}{1-x}\right) \text{Li}_2(1-x)$$ $$-\log (1-x) \log (x) \text{Li}_2(1-x)-\frac{1}{2} \log ^2\left(\frac{x}{1-x}\right) \text{Li}_2\left(\frac{x}{x-1}\right)+\log ^2(1-x) \text{Li}_2(x)+\frac{1}{2} \log ^2(x) \text{Li}_2(x)$$ $$-\log (x) \log \left(\frac{x}{1-x}\right) \text{Li}_2(x)-\log (1-x) \log (x) \text{Li}_2(x)+\log \left(\frac{x}{1-x}\right) \text{Li}_3(1-x)+\log \left(\frac{x}{1-x}\right) \text{Li}_3\left(\frac{x}{x-1}\right)$$ $$+\log \left(\frac{x}{1-x}\right) \text{Li}_3(x)-\log (1-x) \text{Li}_3(x)-\text{Li}_4(1-x)-\text{Li}_4\left(\frac{x}{x-1}\right)-2 \text{Li}_4\left(\frac{x-1}{x}\right)+2 \text{Li}_4(x)+\log (x) \zeta (3)-\frac{\pi ^4}{36}$$\\

$$(4)\ \ \sum _{n=1}^{\infty } \frac{H_n H_n^{(2)} x^n}{n}=\frac{\text{Li}_2(1-x){}^2}{2}-\frac{\text{Li}_2(x){}^2}{2}+\frac{\pi ^2 \text{Li}_2(x)}{6}+\text{Li}_4(1-x)-\text{Li}_4\left(\frac{x}{x-1}\right)$$ $$+\frac{1}{2} \text{Li}_2(1-x) \log ^2\left(\frac{x}{1-x}\right)-\frac{1}{2} \text{Li}_2\left(\frac{x}{x-1}\right) \log ^2\left(\frac{x}{1-x}\right)+\frac{1}{2} \text{Li}_2(1-x) \log ^2(1-x)+\frac{1}{2} \text{Li}_2(1-x) \log ^2(x)$$ $$+\frac{1}{2} \text{Li}_2(x) \log ^2(x)-\text{Li}_2(1-x) \log (x) \log \left(\frac{x}{1-x}\right)-\text{Li}_2(x) \log (x) \log \left(\frac{x}{1-x}\right)+\text{Li}_3(1-x) \log \left(\frac{x}{1-x}\right)$$ $$+\text{Li}_3\left(\frac{x}{x-1}\right) \log \left(\frac{x}{1-x}\right)+\text{Li}_3(x) \log \left(\frac{x}{1-x}\right)-\text{Li}_2(x) \log (1-x) \log (x)-\text{Li}_3(1-x) \log (1-x)+\zeta (3) \log (1-x)$$ $$-\zeta (3) \log (x)+\frac{1}{24} \log ^4\left(\frac{x}{1-x}\right)+\frac{\log ^4(x)}{24}-\frac{1}{6} \log \left(\frac{1}{1-x}\right) \log ^3\left(\frac{x}{1-x}\right)-\frac{1}{6} \log (x) \log ^3\left(\frac{x}{1-x}\right)-\frac{1}{6} \log ^3(x) \log \left(\frac{x}{1-x}\right)$$ $$+\frac{1}{3} \log (1-x) \log ^3(x)+\frac{1}{4} \log ^2(x) \log ^2\left(\frac{x}{1-x}\right)-\frac{1}{2} \log (1-x) \log ^2(x) \log \left(\frac{x}{1-x}\right)+\frac{1}{12} \pi ^2 \log ^2(1-x)-\frac{1}{4} \log ^2(1-x) \log ^2(x)-\frac{\pi ^4}{40}$$}

{\small $$(5)\ \ \sum _{n=1}^{\infty } \frac{H_n^{(3)} x^n}{n}=-\frac{\text{Li}_2(x){}^2}{2}+\text{Li}_4(x)-\text{Li}_3(x) \log (1-x)$$}

{\tiny$$(6)\ \ \sum _{n=1}^{\infty } \frac{\left(H_n\right){}^3 x^n}{n}=\frac{3 \text{Li}_2(1-x){}^2}{2}-\frac{\text{Li}_2(x){}^2}{2}+\frac{\pi ^2 \text{Li}_2(x)}{2}-3 \text{Li}_4(1-x)-3 \text{Li}_4\left(\frac{x}{x-1}\right)-2 \text{Li}_4(x)+\frac{3}{2} \text{Li}_2(1-x) \log ^2\left(\frac{x}{1-x}\right)$$ $$-\frac{3}{2} \text{Li}_2\left(\frac{x}{x-1}\right) \log ^2\left(\frac{x}{1-x}\right)-\frac{3}{2} \text{Li}_2(1-x) \log ^2(1-x)+\frac{3}{2} \text{Li}_2(1-x) \log ^2(x)+\frac{3}{2} \text{Li}_2(x) \log ^2(x)-3 \text{Li}_2(1-x) \log (x) \log \left(\frac{x}{1-x}\right)$$ $$-3 \text{Li}_2(x) \log (x) \log \left(\frac{x}{1-x}\right)+3 \text{Li}_3(1-x) \log \left(\frac{x}{1-x}\right)+3 \text{Li}_3\left(\frac{x}{x-1}\right) \log \left(\frac{x}{1-x}\right)+3 \text{Li}_3(x) \log \left(\frac{x}{1-x}\right)-3 \text{Li}_2(x) \log (1-x) \log (x)$$ $$+3 \text{Li}_3(1-x) \log (1-x)+2 \text{Li}_3(x) \log (1-x)+3 \zeta (3) \log (1-x)-3 \zeta (3) \log (x)+\frac{1}{8} \log ^4\left(\frac{x}{1-x}\right)+\frac{1}{4} \log ^4(1-x)+\frac{\log ^4(x)}{8}$$ $$-\frac{1}{2} \log \left(\frac{1}{1-x}\right) \log ^3\left(\frac{x}{1-x}\right)-\frac{1}{2} \log (x) \log ^3\left(\frac{x}{1-x}\right)-\frac{1}{2} \log ^3(x) \log \left(\frac{x}{1-x}\right)+\log (1-x) \log ^3(x)-\log ^3(1-x) \log (x)$$ $$+\frac{3}{4} \log ^2(x) \log ^2\left(\frac{x}{1-x}\right)-\frac{3}{2} \log (1-x) \log ^2(x) \log \left(\frac{x}{1-x}\right)+\frac{1}{4} \pi ^2 \log ^2(1-x)-\frac{3}{4} \log ^2(1-x) \log ^2(x)-\frac{\pi ^4}{120}$$

$$(7)\ \ \sum _{n=1}^{\infty } \frac{\widetilde{H_n} x^n}{n^2}=-\text{Li}_3\left(\frac{x+1}{2}\right)-\text{Li}_3(-x)-\text{Li}_3(x)-\text{Li}_3\left(\frac{x}{x+1}\right)+\text{Li}_3\left(\frac{2 x}{x+1}\right)+\log (2) \text{Li}_2(x)+\text{Li}_2\left(\frac{x+1}{2}\right) \left(\log (x)-\log \left(\frac{2 x}{x+1}\right)\right)$$ $$+\text{Li}_2\left(\frac{1}{2}-\frac{x}{2}\right) \log (x)+\text{Li}_2(x) \left(\log \left(\frac{x}{x+1}\right)+\log (x+1)\right)+\left(\text{Li}_2\left(\frac{x}{x+1}\right)-\text{Li}_2\left(\frac{2 x}{x+1}\right)\right) \log \left(\frac{2 x}{x+1}\right)$$ $$+\frac{1}{2} \log ^2(2) \log (x)+\frac{1}{2} \log ^2\left(\frac{2 x}{x+1}\right) \left(\log \left(\frac{1-x}{2}\right)+\log \left(\frac{1}{x+1}\right)-\log \left(-\frac{x-1}{x+1}\right)\right)+\log (2) \log (x) \log \left(\frac{2 x}{x+1}\right)$$ $$-\frac{1}{2} \log (2) \log (x) (\log (x)-2 \log (x+1)+\log (4))-\frac{1}{12} \pi ^2 \log (x)+\log \left(\frac{1-x}{2}\right) \log \left(\frac{x+1}{2}\right) \log (x)+\frac{7 \zeta (3)}{8}+\frac{\log ^3(2)}{6}$$

$$(8)\ \ \sum _{n=1}^{\infty } \frac{H_n \widetilde{H_n} x^n}{n}=-\text{Li}_3\left(\frac{1-x}{2}\right)-2 \text{Li}_3\left(\frac{x+1}{2}\right)+\text{Li}_3\left(\frac{x+1}{1-x}\right)-\text{Li}_3(1-x)-\text{Li}_3\left(\frac{x+1}{x-1}\right)-\text{Li}_3(-x)-\text{Li}_3(x)-\text{Li}_3\left(\frac{x}{x+1}\right)$$ $$+\text{Li}_3\left(\frac{2 x}{x+1}\right)+\frac{1}{2} \log (2) \left(2 \text{Li}_2(x)+\log ^2(1-x)\right)+\text{Li}_2\left(\frac{1-x}{2}\right) \log (1-x)-\text{Li}_2\left(\frac{1-x}{2}\right) (\log (1-x)-\log (x))$$ $$+\text{Li}_2\left(\frac{x+1}{2}\right) \left(\log (x)-\log \left(\frac{2 x}{x+1}\right)\right)+\text{Li}_2\left(\frac{x+1}{2}\right) \log (x+1)+\text{Li}_2(1-x) \left(\log (x+1)-\log \left(\frac{x+1}{1-x}\right)\right)$$ $$+\left(\text{Li}_2\left(\frac{x+1}{x-1}\right)-\text{Li}_2\left(\frac{x+1}{1-x}\right)\right) \log \left(\frac{x+1}{1-x}\right)+\text{Li}_2(-x) (\log (1-x)-\log (x))+\text{Li}_2(-x) \log (x)+\text{Li}_2(x) \left(\log \left(\frac{x}{x+1}\right)+\log (x+1)\right)$$ $$+\left(\text{Li}_2\left(\frac{x}{x+1}\right)-\text{Li}_2\left(\frac{2 x}{x+1}\right)\right) \log \left(\frac{2 x}{x+1}\right)+\text{Li}_2(x+1) \left(\log \left(\frac{x+1}{1-x}\right)+\log (1-x)\right)-\text{Li}_2(x+1) \log (x+1)$$ $$+\frac{1}{2} \log \left(\frac{1-x}{2}\right) \log ^2(x+1)-\frac{1}{2} \log (-x) \log ^2(x+1)+\frac{1}{12} \left(\pi ^2-6 \log ^2(2)\right) \log (1-x)+\frac{1}{2} \log ^2\left(\frac{x+1}{1-x}\right) \left(-\log \left(\frac{2 x}{x-1}\right)+\log \left(-\frac{2}{x-1}\right)+\log (x)\right)$$ $$-\frac{1}{12} \left(\pi ^2-6 \log ^2(2)\right) \log (x)+\frac{1}{2} \log ^2\left(\frac{2 x}{x+1}\right) \left(\log \left(\frac{1-x}{2}\right)+\log \left(\frac{1}{x+1}\right)-\log \left(-\frac{x-1}{x+1}\right)\right)+\log (2) \log (x) \log \left(\frac{2 x}{x+1}\right)$$ $$-\frac{1}{2} \log (2) \log (x) (\log (x)-2 \log (x+1)+\log (4))+\log \left(\frac{1-x}{2}\right) \log \left(\frac{x+1}{2}\right) \log (x)+\log \left(\frac{x+1}{1-x}\right) (\log (-x)-\log (x)) \log (x+1)$$ $$+\log (1-x) \log (x) \log (x+1)+\frac{1}{2} (\log (x)-\log (-x)) \log (x+1) (\log (x+1)-2 \log (1-x))+\frac{15 \zeta (3)}{8}+\frac{\log ^3(2)}{2}-\frac{1}{6} \pi ^2 \log (2)$$

$$(9)\ \ \sum _{n=1}^{\infty } \frac{\left(\widetilde{H_n}\right)^2 x^n}{n}=-4 \text{Li}_3\left(\frac{1-x}{2}\right)-2 \text{Li}_3\left(\frac{x+1}{2}\right)+2 \text{Li}_3\left(\frac{x+1}{1-x}\right)-2 \text{Li}_3\left(\frac{x+1}{x-1}\right)-\text{Li}_3(x)-2 \text{Li}_3\left(\frac{x}{x+1}\right)+2 \text{Li}_3\left(\frac{2 x}{x+1}\right)$$ $$-2 \text{Li}_3(x+1)+\log (2) \left(2 \text{Li}_2(x)+\log ^2(1-x)\right)+4 \text{Li}_2\left(\frac{1-x}{2}\right) \log (1-x)-2 \text{Li}_2\left(\frac{1-x}{2}\right) (\log (1-x)-\log (x))+2 \text{Li}_2\left(\frac{x+1}{2}\right) \left(\log (x)-\log \left(\frac{2 x}{x+1}\right)\right)$$ $$+2 \text{Li}_2(1-x) \left(\log (x+1)-\log \left(\frac{x+1}{1-x}\right)\right)-2 \text{Li}_2(1-x) \log (1-x)+2 \left(\text{Li}_2\left(\frac{x+1}{x-1}\right)-\text{Li}_2\left(\frac{x+1}{1-x}\right)\right) \log \left(\frac{x+1}{1-x}\right)-\text{Li}_2(x) (\log (1-x)-\log (x))$$ $$-\text{Li}_2(x) \log (x)+2 \text{Li}_2(x) \left(\log \left(\frac{x}{x+1}\right)+\log (x+1)\right)+2 \left(\text{Li}_2\left(\frac{x}{x+1}\right)-\text{Li}_2\left(\frac{2 x}{x+1}\right)\right) \log \left(\frac{2 x}{x+1}\right)+2 \text{Li}_2(x+1) \left(\log \left(\frac{x+1}{1-x}\right)+\log (1-x)\right)$$ $$-\log (2) \log ^2(1-x)+\frac{1}{6} \left(\pi ^2-6 \log ^2(2)\right) \log (1-x)+\log ^2\left(\frac{x+1}{1-x}\right) \left(-\log \left(\frac{2 x}{x-1}\right)+\log \left(-\frac{2}{x-1}\right)+\log (x)\right)-\log ^2(1-x) \log (x)$$ $$-\frac{1}{6} \left(\pi ^2-6 \log ^2(2)\right) \log (x)+\log ^2\left(\frac{2 x}{x+1}\right) \left(\log \left(\frac{1-x}{2}\right)+\log \left(\frac{1}{x+1}\right)-\log \left(-\frac{x-1}{x+1}\right)\right)+2 \log (2) \log (x) \log \left(\frac{2 x}{x+1}\right)$$ $$-\log (2) \log (x) (\log (x)-2 \log (x+1)+\log (4))+2 \log \left(\frac{1-x}{2}\right) \log \left(\frac{x+1}{2}\right) \log (x)+2 \log \left(\frac{x+1}{1-x}\right) (\log (-x)-\log (x)) \log (x+1)$$ $$+2 \log (1-x) \log (x) \log (x+1)+(\log (x)-\log (-x)) \log (x+1) (\log (x+1)-2 \log (1-x))+\frac{15 \zeta (3)}{4}+\log ^3(2)-\frac{1}{3} \pi ^2 \log (2)$$}

(10) The following four generating functions are polylog expressible, using formula (1)-(6): $$\sum _{n=1}^{\infty } \left(H_n\right){}^4 x^n,\sum _{n=1}^{\infty } \left(H_n\right){}^2 H_n^{(2)} x^n,\sum _{n=1}^{\infty } \left(H_n^{(2)}\right){}^2 x^n,\sum _{n=1}^{\infty } H_n H_n^{(3)} x^n$$\\

\noindent Proof. (1)-(9): Again, make use of $L(f)$ in appropriate functions in Proposition 3, then Mathematica. They are verified via differentiation. (10): Calculate differences. The results will be combinations of (1)-(6) thus extremely lengthy, so we omit them. $\hfill\square$\\\\

\noindent \textbf{Proposition 5.} The following formulas hold:\\
\begin{equation}
\begin{aligned}
\int_0^1 x^{n-1} \log^k(x) \, dx=\frac{(-1)^k k!}{n^{k+1}}
\end{aligned}
\tag{1}
\end{equation}

\begin{equation}
\begin{aligned}
\int_0^1 x^{n-1} \log (1-x) \, dx=-\sum _{j=1}^{\infty } \frac{1}{j (j+n)}=-\frac{H_n}{n}
\end{aligned}
\tag{2}
\end{equation}

\begin{equation}
\begin{aligned}
\int_0^1 x^{n-1} \log ^2(1-x) \, dx=\sum _{j=1}^{\infty } \frac{1}{j (j+n)^2}=\frac{\left(H_n\right){}^2+H_n^{(2)}}{n}
\end{aligned}
\tag{3}
\end{equation}

\begin{equation}
\begin{aligned}
\int_0^1 x^{n-1} \log (x) \log (1-x) \, dx=\frac{H_n}{n^2}+\frac{H_n^{(2)}}{n}-\frac{\pi ^2}{6 n}
\end{aligned}
\tag{4}
\end{equation}

\begin{equation}
\begin{aligned}
\int_0^1 x^{n-1} \log ^3(1-x) \, dx=-\frac{3 H_n H_n^{(2)}+\left(H_n\right){}^3+2 H_n^{(3)}}{n}
\end{aligned}
\tag{5}
\end{equation}
\begin{equation}
\begin{aligned}
-\int_0^1 x^{n-1} \log (x+1) \, dx=\sum _{j=1}^{\infty } \frac{(-1)^j}{j (j+n)}=\frac{(-1)^{n-1} \widetilde{H_{n}^{}}-(-1)^{n-1} \log (2)-\log (2)}{n}
\end{aligned}
\tag{6}
\end{equation}

\begin{equation}
\begin{aligned}
\sum _{j=1}^{\infty } \frac{(-1)^j}{j (j+n)^2}=\frac{(-1)^{n-1} \widetilde{H_{n}^{\left(2 \right)}}}{n}+\frac{(-1)^{n-1} \widetilde{H_{n}^{}}}{n^2}-\frac{(-1)^{n-1} \log (2)}{n^2}-\frac{\log (2)}{n^2}-\frac{\pi ^2 (-1)^{n-1}}{12 n}
\end{aligned}
\tag{7}
\end{equation}\\

\noindent Proof. (1): Trivial. (2)-(5): Differentiate Beta functions, then transform polygamma constants into harmonic ones (for instance $\psi ^{(2)}(n+1)=2 H_n^{(3)}-2 \zeta (3),\psi ^{(3)}(n+1)=\frac{\pi ^4}{15}-6 H_n^{(4)}$, etc). (6)-(7): Invoke partial fractions and that $\eta (1)=\log (2),\eta (2)=\frac{\pi ^2}{12}$. $\hfill\square$\\

\subsubsection{Confrontation}
Now we compute 7 remaining weight 5 Euler sums:\\

\noindent \textbf{ Proposition 6.} The following formulas hold:\\

$$(1)\ ES(-1,-1,-1;-2)$$$$=-18 \text{Li}_5\left(\frac{1}{2}\right)-\frac{\pi ^2 \zeta (3)}{8}+\frac{1229 \zeta (5)}{64}+\frac{21}{8} \zeta (3) \log ^2(2)+\frac{3 \log ^5(2)}{20}+\frac{3}{4} \pi ^2 \log ^3(2)-\frac{29}{160} \pi ^4 \log (2)$$

$$(2)\ ES(-1,-1,-1,-1;-1)$$$$=-48 \text{Li}_5\left(\frac{1}{2}\right)-4 \text{Li}_4\left(\frac{1}{2}\right) \log (2)+\frac{5 \pi ^2 \zeta (3)}{48}+\frac{733 \zeta (5)}{16}+\frac{13 \log ^5(2)}{30}+\frac{3}{2} \pi ^2 \log ^3(2)-\frac{71}{180} \pi ^4 \log (2)$$

$$(3)\ ES(1,1,-2;-1)$$$$=6 \text{Li}_5\left(\frac{1}{2}\right)-\frac{\pi ^2 \zeta (3)}{48}-\frac{93 \zeta (5)}{64}-\frac{1}{20} \log ^5(2)+\frac{1}{36} \pi ^2 \log ^3(2)-\frac{13 \pi ^4 \log (2)}{1440}$$

$$(4)\ ES(1,1,-1,-1;-1)$$$$=12 \text{Li}_5\left(\frac{1}{2}\right)+\frac{\pi ^2 \zeta (3)}{16}-\frac{201 \zeta (5)}{16}-\frac{3}{8} \zeta (3) \log ^2(2)+\frac{\log ^5(2)}{10}+\frac{1}{12} \pi ^2 \log ^3(2)+\frac{97}{720} \pi ^4 \log (2)$$

$$(5)\ ES(1,-1,-1;-2)$$$$=10 \text{Li}_5\left(\frac{1}{2}\right)-2 \text{Li}_4\left(\frac{1}{2}\right) \log (2)+\frac{\pi ^2 \zeta (3)}{96}-10 \zeta (5)-\frac{7}{8} \zeta (3) \log ^2(2)-\frac{1}{6} \log ^5(2)+\frac{2}{9} \pi ^2 \log ^3(2)+\frac{1}{10} \pi ^4 \log (2)$$

$$(6)\ ES(1,-1,-2;-1)$$$$=4 \text{Li}_4\left(\frac{1}{2}\right) \log (2)+\frac{7}{16} \zeta (3) \log ^2(2)+\frac{\log ^5(2)}{6}-\frac{5}{72} \pi ^2 \log ^3(2)-\frac{1}{480} \pi ^4 \log (2$$

$$(7)\ ES(1,-1,-1,-1;-1)$$$$=2 \text{Li}_5\left(\frac{1}{2}\right)+\frac{13 \pi ^2 \zeta (3)}{96}-\frac{83 \zeta (5)}{32}+\frac{9}{16} \zeta (3) \log ^2(2)-\frac{13}{60}  \log ^5(2)+\frac{41}{72} \pi ^2 \log ^3(2)-\frac{1}{180} \pi ^4 \log (2)$$\\

\noindent Proof. (1)-(2): Exploit symmetry. Take (2) as an example, we have
$$ES(-1,-1,-1,-1;-1)=\sum _{0<j,k,l,m\leq n}^{} \frac{(-1)^{j-1} (-1)^{k-1} (-1)^{l-1} (-1)^{m-1} (-1)^{n-1}}{j k l m n}$$
Consider the analogous sum without restrictions, by I-E principle
$$\log ^5(2)=5 ES(-1,-1,-1,-1;-1)-10ES(-1,-1,-1;2)+10ES(-1,-1;-3)-5 ES(-1;4)+\eta (5)$$
Since all ESs but the desired one are known, (2) is established. (1) is similar.\\

(3): Multiply both sides of Prop. 5(7) with $\left(H_n\right){}^2$ then sum w.r.t $n$, we have:\\
\begin{equation}
\begin{aligned}
-\frac{\pi ^2}{12}ES(1,1;-1)-\log (2)ES(1,1;2)-\log (2)ES(1,1;-2)\ \ \\+ES(1,1,-1;-2)+ES(1,1,-2;-1)=\sum _{n=1}^{\infty } \sum _{j=1}^{\infty } \frac{(-1)^j \left(H_n\right){}^2}{j (j+n)^2}
\end{aligned}
\notag
\end{equation}\\
Using formulas above we have: \\
\begin{equation}
\begin{aligned}
\sum _{n=1}^{\infty } \sum _{j=1}^{\infty } \frac{(-1)^j \left(H_n\right){}^2}{j (j+n)^2}=\int_0^1 \frac{\log (x)}{x} \left(\sum _{j=1}^{\infty } \frac{(-1)^{j-1} x^j}{j}\right) \left(\sum _{n=1}^{\infty } \left(H_n\right){}^2 x^n\right) \, dx\\ =\int_0^1 \left(\frac{1}{x}+\frac{1}{1-x}\right) \log (x) \log (x+1) \left(\text{Li}_2(x)+\log ^2(1-x)\right) \, dx\ \ \ \ \ \ \
\end{aligned}
\notag
\end{equation}\\
Therefore we only need to evaluate $PLI(0,1,1;21;0), PLI(0,1,1;21;1)$. Indeed:\\
\begin{equation}
\begin{aligned}
PLI(0,1,1;21;1)=\sum _{n=1}^{\infty } \sum _{j=1}^{\infty } \frac{(-1)^{j-1} \int_0^1 \log (x) x^{j+n-1} \, dx}{j n^2}=\sum _{n=1}^{\infty } \sum _{j=1}^{\infty } \frac{(-1)^j}{j n^2 (j+n)^2}
\end{aligned}
\notag
\end{equation}
\begin{equation}
\begin{aligned}
PLI(0,1,1;21;0)=\cdots=-\sum _{n=1}^{\infty } \sum _{j=1}^{\infty } \frac{\widetilde {H_j}}{n^2 (j+n+1)^2}
\end{aligned}
\notag
\end{equation}\\
For the first one, sum w.r.t either $n$ or $j$ first and the PLI reduces to known ESs. For the second, sum w.r.t $n$ yields\\
\begin{equation}
\begin{aligned}
-\sum _{n=1}^{\infty } \sum _{j=1}^{\infty } \frac{\widetilde {H_j}}{n^2 (j+n+1)^2}=-\sum _{j=1}^{\infty } \left(\widetilde {H_j}-\frac{(-1)^{j-1}}{j}\right) \left(-\frac{2 H_j}{j^3}-\frac{H_j^{(2)}}{j^2}+\frac{\pi ^2}{3 j^2}\right)
\end{aligned}
\notag
\end{equation}\\
All ES components on RHS are also known so (3) is proved.\\

(4): Multiply both sides of Prop. 5(3) with $\left(-1\right)^{n-1} {\widetilde {H_n}}^2$ and sum $n$, the ES boils down to $PLI(2,0,0;22;1/2), PLI(2,0,0;24;1/2)$ via Prop. 3(9). For $PLI(2,0,0;22;1/2)$, use Prop. 3(5) for expansion to reduce them to known ESs. For $PLI(2,0,0;24;1/2)$, use that:\\
$$\text{Li}_2\left(\frac{1-x}{2}\right)+\text{Li}_2\left(\frac{x+1}{2}\right)=\log (2) \log (1-x)+\log (2) \log (x+1)-\log (1-x) \log (x+1)+\frac{\pi ^2}{6}-\log ^2(2)$$\\
to reduce them to $PLI(2,0,0;23;1/2)$ modulo LIs. For the final 2 PLIs, by substitution $t\to1-x$, single/double expansion for $\frac{\text{Li}_2\left(\frac{x}{2}\right)}{1-\frac{x}{2}}/\frac{\text{Li}_2\left(\frac{x}{2}\right)}{1-x}$ respectively and ES/GES results \cite{ref19} they are solved either so (4) is proved.\\

(5): Same as (4), except that we start from Prop. 5(4). One may verify that all ESs appeared during reduction are known and all PLIs can be reduced to known ESs.\\

(6): Multiply both sides of Prop. 5(7) with $H_n \widetilde {H_n}$, remaining steps are alike to (3)-(5), in which the value of (5) is used.\\

(7): Denote $f(x)=\sum _{n=1}^{\infty } \left(\widetilde{H_n}\right){}^3 x^n$, by calculating differences of the generating function we have $f(-x)=\frac{-3 P(x)-3 Q(x)-\text{Li}_3(x)}{x+1}$, where $P(x)=\sum _{n=1}^{\infty } \frac{\left(\widetilde{H_n}\right){}^2 x^n}{n}, Q(x)=\sum _{n=1}^{\infty } \frac{(-1)^n \widetilde{H_n} x^n}{n^2}$. Thus\\
\begin{equation}
\begin{aligned}
ES(1,-1,-1,-1;-1)=\int_0^1 \frac{f(-x) \log (1-x)}{x} \, dx=\int_0^1 \frac{\log (1-x) (-3 P(x)-\text{Li}_3(x)-3 Q(x))}{x (x+1)} \, dx
\end{aligned}
\notag
\end{equation}\\
Integrate by parts using\\
 $$\int \frac{\log (1-x)}{x (x+1)} \, dx=u(x)=\frac{\pi ^2}{6}-\text{Li}_2\left(\frac{1-x}{2}\right)-\text{Li}_2(x)-\log (1-x) \log \left(\frac{x+1}{2}\right)$$\\
 We only have to compute\\
\begin{equation}
\begin{aligned}
\int_0^1 \frac{u(x) }{x} \left(3 \sum _{n=1}^{\infty } \frac{(-1)^n}{n} \widetilde{H_n} x^n+3 \sum _{n=1}^{\infty } \left(\widetilde{H_n}\right){}^2 x^n+\text{Li}_2(x)\right)\, dx
\end{aligned}
\notag
\end{equation}\\
Which, by Proposition 3, equals to\\
\begin{equation}
\begin{aligned}
3\int_0^1 \frac{u(x) }{x (1-x)} \left(2 \text{Li}_2\left(\frac{1-x}{2}\right)+\text{Li}_2(x)-2 \text{Li}_2\left(\frac{1}{2}\right)+2 \log \left(\frac{x+1}{2}\right) \log (1-x)\right)\, dx 
\end{aligned}
\notag
\end{equation}
\begin{equation}
\begin{aligned}
 +3\int_0^1 \frac{u(x) }{x}\left(\text{Li}_2\left(\frac{x+1}{2}\right)-\text{Li}_2(x)-\text{Li}_2\left(\frac{1}{2}\right)-\log (2) \log (x+1)\right) \, dx  +\int_0^1 \frac{\text{Li}_2(x) u(x)}{x} \, dx
\end{aligned}
\notag
\end{equation}\\

Expand all these expressions. Some parts of the integrand are divergent so we put them together again, for instance $\int_0^1 \frac{\text{Li}_2(x) \left(\text{Li}_2\left(\frac{x+1}{2}\right)-\text{Li}_2\left(\frac{1}{2}\right)\right)}{x} \, dx$. Apply integration by parts on all these composite integrals (in the example lift up a $\text{Li}_3(x)$) to convert them into convergent PLIs. After this we reduce it to plenty of LIs and PLIs. For each of them, integrate by parts and use techniques in (3)-(6) over and over again (or more generally, use the expansion algorithm in section 4), all PLIs can be reduced to lower weight ESs and 53 weight 5 ESs already known, hence (7) is proved.$\hfill\square$\\

Thus all 89 ESs with $W\leq5$ are evaluated.\\

\section{PLI}

Due to existence of a polylog term, a PLI's weight $W\geq3$. There are 11, 55, 197 PLIs for weight 3, 4, 5 respectively.\\

\subsection{Brute force}
Some of the PLIs can be calculated by brute force, for instance $\frac{\text{Li}_2\left(\frac{x+1}{2}\right) \log (1-x)}{x+1}$ has primitive $\log (2) \text{Li}_3\left(\frac{x+1}{2}\right)-\frac{1}{2} \text{Li}_2\left(\frac{x+1}{2}\right)^2$. For each PLI, check by Mathematica whether the primitive is solvable and evaluate it if so.\\

\subsection{General formulas}

See \cite{ref12, ref15} for proofs of the following:\\

\noindent \textbf{Lemma 3.} $\int_0^1 \frac{\text{Li}_m(x) \text{Li}_n(x)}{x} \, dx$, $\int_0^1 \frac{\log ^m(x) \text{Li}_n(x)}{1-x} \, dx$ ($m+n$ even) and $\int_0^1 \frac{\text{Li}_{n+1}(x) \log ^n(x)}{x+1} \, dx$ are reducible to Zeta values.\\

Thus PLIs of these types are solved directly.\\

\subsection{Integration by parts}

In spirit of subsection 2-3, generalize $F$ to $\prod _{m=0}^2 l(m;x)^{a(m)} \prod _{k=1}^6 \prod _{n=2}^N L(n,k;x)^{b(n,k)}$. For those involving $L(n,k;x), k\in\{3,5,6\}$, after partial fractions another integration by parts lifting $\log(1-x)$ is needed to obtain a regular relation modulo LIs. Now consider all possible $F$.\\

\subsection{Fractional transformation}

Similar to subsection 2-4, apply $x\to\frac{1-x}{x+1}$ in integral of $F=\frac{\prod _{m=0}^2 l(m;x)^{a(m)} \prod _{k=1}^6 \prod _{n=2}^N L(n,k;x)^{b(n,k)}}{f(p;x)}$, where $b(n,3)=b(n,4)=0$. For $p=1$, after partial fractions another integration by parts lifting $\log(1-x)$ is needed to obtain a regular relation modulo LIs. Now consider all possible $F$.\\

\subsection{Complementary integrals}

Denote $P(x)=\prod _{n=2}^N \text{Li}_n(x){}^{a(n)}$, by equally splitting $(0,1)$ and substituting $t=\frac{1\pm x}{2}$ we have:

{\scriptsize $$\int_0^1 \frac{P\left(\frac{1-x}{2}\right) \log ^m\left(\frac{1-x}{2}\right) \log ^n\left(\frac{x+1}{2}\right)}{1\mp x}+\frac{P\left(\frac{x+1}{2}\right) \log ^m\left(\frac{x+1}{2}\right) \log ^n\left(\frac{1-x}{2}\right)}{1\pm x}dx=\int_0^1 \frac{P(x) \log ^m(x) \log ^n(1-x)}{f_{\pm}(x)} \, dx$$}\\
Here $f_{+}(x)=x, f_{-}(x)=1-x$. Expand LHS to get a linear relation. Now consider all possible $P$.\\

\subsection{Polylog identities}

\noindent \textbf{Proposition 7.} The following formulas hold on $(0,1)$ and naturally extend to the plane, where all logs/polylogs are principal:\\
\begin{equation}
\begin{aligned}
\text{Li}_n(-x)+\text{Li}_n(x)=\frac{\text{Li}_n\left(x^2\right)}{2^{n-1}}
\end{aligned}
\tag{1}
\end{equation}

\begin{equation}
\begin{aligned}
\text{Li}_2(1-x)+\text{Li}_2(x)=\frac{\pi ^2}{6}-\log (x) \log (1-x)
\end{aligned}
\tag{2}
\end{equation}

\begin{equation}
\begin{aligned}
\text{Li}_2\left(\frac{x}{x-1}\right)+\text{Li}_2(x)=-\frac{1}{2} \log ^2(1-x)
\end{aligned}
\tag{3}
\end{equation}

\begin{equation}
\begin{aligned}
\text{Li}_2\left(-\frac{1}{x}\right)+\text{Li}_2(-x)=-\frac{\pi ^2}{6}-\frac{1}{2} \log ^2(x)
\end{aligned}
\tag{4}
\end{equation}

{\footnotesize\begin{equation}
\begin{aligned}
\text{Li}_2\left(\frac{1-x}{2}\right)+\text{Li}_2\left(\frac{x+1}{2}\right)=\log (2) \log (x+1)+\log (2)\log (1-x) -\log (1-x) \log (x+1)+\frac{\pi ^2}{6}-\log ^2(2)
\end{aligned}
\tag{5}
\end{equation}}

\begin{equation}
\begin{aligned}
\text{Li}_2\left(\frac{1-x}{2}\right)+\text{Li}_2\left(\frac{x-1}{x+1}\right)=-\frac{1}{2} \log ^2(x+1)+\log (2) \log (x+1)-\frac{1}{2} \log ^2(2)
\end{aligned}
\tag{6}
\end{equation}

\begin{equation}
\begin{aligned}
\text{Li}_3(-x)-\text{Li}_3\left(-\frac{1}{x}\right)=-\frac{1}{6} \log ^3(x)-\frac{1}{6} \pi ^2 \log (x)
\end{aligned}
\tag{7}
\end{equation}

{\small\begin{equation}
\begin{aligned}
\text{Li}_3\left(\frac{x}{x-1}\right)+\text{Li}_3(1-x)+\text{Li}_3(x)=\frac{1}{6} \log ^3(1-x)-\frac{1}{2} \log (x) \log ^2(1-x)+\frac{1}{6} \pi ^2 \log (1-x)+\zeta (3)
\end{aligned}
\tag{8}
\end{equation}}

\begin{equation}
\begin{aligned}
\text{Li}_3\left(\frac{1-x}{2}\right)+\text{Li}_3\left(\frac{x+1}{2}\right)+\text{Li}_3\left(\frac{x-1}{x+1}\right)=\frac{1}{6} \log ^3(x+1)-\frac{1}{2} \log ^2(2) \log (1-x)\\-\frac{1}{2} \log ^2(x+1) \log (1-x)+\left(\frac{\pi ^2}{6}-\frac{\log ^2(2)}{2}\right) \log (x+1)\ \ \ \ \ \ \ \ \ \ \ \ \ \ \ \ \ \  \\+\log (2) \log (1-x) \log (x+1)+\zeta (3)+\frac{\log ^3(2)}{3}-\frac{1}{6} \pi ^2 \log (2)\ \ \ \ \ \ \ \ \ \ \ \ \ \ \ 
\end{aligned}
\tag{9}
\end{equation}\\

\noindent Proof. (1)-(4), (7)-(8): Classical formulas, see \cite{ref14}. (5), (6), (9): Substitute $t\to \frac{1-x}{2}$ in (2), (3), (8) respectively.  $\hfill\square$\\

Now, multiply both sides of (1) with $F(x)=\frac{\log ^m\left(1-x^2\right) \log^k(x)}{x}$, integrate on $(0,1)$. Substitute $t\to x^2$ on RHS, expand $(\log (1-x)+\log (x+1))^m$ on LHS, we obtain a relation. Another choice is $F(x)=\frac{x \log ^m\left(1-x^2\right) \log^k(x)}{1-x^2}$, where we also use partial fractions to decompose LHS. Now consider all possible $F$.\\

Moreover, multiply both sides of (9) with $F(x)=\frac{\prod _{m=0}^2 l(m;x)^{a(m)}}{f(p;x)}$, it generates a relation between 3 PLIs since all LIs on RHS are known. For (5), (6) the process is similar, except that we may broaden the range of $F$ to $\frac{\prod _{m=0}^2 l(m;x)^{a(m)}}{f(p;x)}$ and $\frac{L(2,k;x)}{f(p;x)}$. Now consider all possible $F$ for 3 formulas to engender linear relations.\\

Clarification: If we derive all closed-forms of PLIs involving $L(2,2;x)$ hence $L(2,6;x)$ (via fractional transformation), those involving $L(2,3/4;x)$ will be natural corollaries thus redundant in the systematic solution. However, as we will see in the next subsection, solving all $L(2,2;x)$ type PLIs directly is unachievable. As a matter of fact, some of the independent relations only generate from manipulating $L(2,3/4;x)$ type, say subsection 4-7 (reflection of $\text{Li}_n\left(\frac{1-x}{2}\right)$), 6-10-4, hence this part is necessary.\\

\subsection{Series expansion}

Weight 5 PLIs are sorted into 4 classes, that is, (2,1,1), (3,1), (4), (2,2). The index indicates the quantity as well as order of polylog terms, for instance, $PLI(1,1,0;22;2)$ belongs to (2,1,1) class for it contains one $\text{Li}_2$ term and two log terms in numerator. Note that class (1,1,1,1,1) is excluded for they are simply LIs.\\

From now on we only discuss about weight 5 PLIs for lower weights are relevantly trivial. Moreover, we restrict the use of this method to (2,1,1), (3,1) and (4) class PLIs. Furthermore, we require the only polylog term  $L(n,k;x)$ ($n=2, 3, 4$) contained in the PLI to satisfy $1\leq k\leq3$. Denote the set of these PLIs as $Q$.\\

We only present the general procedure of converting the PLI into combination of ESs. In this subsection 'ES' represents all ordinary ESs and all GESs.\\

Step 1. Reflection. If the only polylog term involved is $\text{Li}_n\left(\frac{1-x}{2}\right)$, substitute $t=1-x$ in the integral. If the only polylog term involved is $\text{Li}_n(-x)$, keep it unchanged. If the only polylog term involved is $\text{Li}_n(x)$, keep it unchanged, or reduce the integral to the corresponding one involving $\text{Li}_n(1-x)$ via (3.4.2) modulo known LIs, then substitute $t=1-x$. After this operation, the argument in polylog term is either $x$ or $\frac{x}{2}$. In this step we have 2 choices for case $\text{Li}_n(x)$ and 1 choice for other 2 cases.\\

Step 2. Determine the invariant part of integrand. Since we will use (2.5.16)-(2.5.19) to convert the integral to ESs, we choose a part of the integrand to be invariant. This part should be of form $F_1(x)=\log ^m(x) \log ^n(1-x)$, with $m$ equals to the degree of $\log(x)$ in the original integrand, $n$ do not exceed the degree of $\log(1-x)$ in the original integrand. In this step we have at most 3 choices on determining $n$, namely $n=0, 1, 2$.\\

Step 3. Multiple expansion. Denote $F=F_1 F_2$, where $F$ is the whole integrand after possible manipulations in step 1, $F_1$ the invariant part determined in step 2. We have $F_2$ a combination of products of terms belong to the set:
$$\left\{\frac{1}{x},\frac{1}{1-x},\frac{1}{x+1},\frac{1}{1-\frac{x}{2}},\log (1-x),\log (x+1),\log \left(1-\frac{x}{2}\right),\text{Li}_n(x),\text{Li}_n\left(\frac{x}{2}\right)\right\}$$
Where by combination we mean the integral is composed of several integrals with different weights (see examples below). Separate $F_2$, expand each term other than $\frac{1}{x}$ into elementary power series, exchange the order of evaluation, we get a at most quadruple sum with the summand an integral consistent with Prop. 5. In this step we have only 1 choice for the process of integration is definite.\\

Step 4. Converting into multiple harmonic series. This is straightforward if we use Prov. 5 to evaluate the integral in rational/harmonic terms. In this step we have only 1 choice for the process of integration is definite.\\

Step 5. Cauchy simplification. After step 4 we arrive at a multiple sum w.r.t at most 4 summation indexes. If there's only 1 index, it's already a combination of ESs hence the procedure is finished. If there are 2 indexes, go to the next step. If there are more than 2 indexes, consider 2 manipulations: (1) Sum w.r.t all indexes with rational summands to reduce the weight. (2) Consider all index pairs $(j,k)$, check whether their Cauchy product is expressible in harmonic terms with a new index $m$, if so, replace $(j,k)$ with $m$ in the summation to reduce the weight. Carry out (1)(2) until it's not reducible any more. If the resulting sum is double, go to the next step, otherwise discard this sum and stop the procedure. In this step we have a great deal of choices to simplify the sum, but much of them will be discarded.\\

Step 6. Reduction to ESs. If we go through step 5, we will have a double series in hand. Denote the index of summation as $(j,k)$. If the summand is rational functions of $j, k$, sum w.r.t either of the indexes first to finish the procedure. If the summand is rational to $j$ but harmonic to $k$, sum w.r.t $j$ also finishes. If we have all harmonic terms with index $j+k$ and other parts rational, invoke Cauchy product again to finish the procedure. For other cases (e.g. the sum involves $H_k H_{j+k}$), discard it. In this step we have 1 or 2 choices for using iterated summation or Cauchy simplification.\\

Fix a weight 5 PLI that belongs to $Q$, consider the tree of choices that we would make through step 1$\sim$6 of the procedure above. If there exist a branch of this tree that is not ended with 'discarded' in step 5$\sim$6, we have reduced the PLI to ESs successfully. Since all weight 5 ESs are already evaluated, we readily obtained the closed-form of this PLI. If all branches are discarded, we get nothing new. Now carry out the procedure for every PLI in $Q$ to obtain independent relations giving PLI closed-forms directly. We show some typical examples:\\

A successful example: Consider $PLI(1,0,1;21;2)=\int_0^1 \frac{\text{Li}_2(x) \log (1-x) \log (x+1)}{x+1} \, dx$.\\

Step 1: Modulo known LIs, we only need to evaluate $\int_0^1 \frac{\text{Li}_2(1-x) \log (1-x) \log (x+1)}{x+1} \, $, so let $t\to 1-x$.\\

Step 2: The problem boils down to evaluating $\int_0^1 \log(2) \frac{\text{Li}_2(x) \log (x)}{1-\frac{x}{2}} \, dx$ and $\int_0^1 \frac{\text{Li}_2(x) \log (x) \log \left(1-\frac{x}{2}\right)}{1-\frac{x}{2}} \, dx$. Set $\log(x)$ to be invariant in both 2 integrals. Making use of double expansion the weight 4 integral is easily solved. Now we deal with the latter one.\\

Step 3: Expand $\text{Li}_2(x),\log \left(1-\frac{x}{2}\right),\frac{1}{1-\frac{x}{2}}$ with 3 independent index $j, k, l$.\\

Step 4: Make use of Prop. 5(1).\\

Step 5: Consider Cauchy product of the summation w.r.t $k, l$. The convolutional term is harmonic hence the triple sum is reduced to double (or we may expand $\frac{\log \left(1-\frac{x}{2}\right)}{1-\frac{x}{2}}$ directly in step 3 to obtain the double sum).\\

Step 6: Denote the Cauchy index replacing $k,l$ as $m$. The double sum we have in hand is $\sum _{j=1}^{\infty } \sum _{m=1}^{\infty } \frac{H_m}{j^2 2^m (j+m+1)^2}$. Sum w.r.t $j$ first (this is nothing but the elementary sum occurred in proof of Prop. 6(3)), we reduced this sum to a combination of three GESs already known. Hence
\begin{equation}
\begin{aligned}
PLI(1,0,1;21;2)=-2 \text{Li}_5\left(\frac{1}{2}\right)-\text{Li}_4\left(\frac{1}{2}\right) \log (2)+\frac{23 \pi ^2 \zeta (3)}{96}-\frac{15 \zeta (5)}{64}\\-\frac{21}{16} \zeta (3) \log ^2(2)-\frac{1}{40} \log ^5(2)+\frac{5}{72} \pi ^2 \log ^3(2)-\frac{11}{720} \pi ^4 \log (2)\ \ \ \ \ \ \ \
\end{aligned}
\notag
\end{equation}\\

Unsuccessful examples: Readers may verify that $PLI(1,0,1;22;2), PLI(2,0,0;22;2)$ are unsolvable using this procedure. More specifically, it is impossible to make choices in the procedure above, to avoid the multiple sum containing summand of form $H_k H_{j+k}$ or its alternating analogue, hence in step 6 they must be discarded. In fact these two integrals are evaluated nontrivially as we put all relations generated by this section together. This reveals the necessity of a systematic approach providing linear relations.\\

\subsection{Solution}

All linear relations generated by section 3 and 4 forms a equation system, from which all ESs and PLIs with weight $\leq5$ are explicitly solved thus Theorem 1 is proved.\\

Note: Based on Prop. 7(1)(2)(3) and all LI/PLIs computed, we may add the following \\
$$\text{Li}_2\left(x^2\right),\text{Li}_2(1-x),\text{Li}_2\left(\frac{x}{x-1}\right),\text{Li}_2\left(\frac{x}{x+1}\right)$$\\
to the set of 6 polylog terms with new PLIs remain in $A_W$ for corresponding $W$. As an example, see \cite{ref15} for evaluation of $\int_0^1 \frac{\text{Li}_2\left(\frac{x}{x-1}\right) \text{Li}_2\left(\frac{x}{x+1}\right)}{x} \, dx$. These new PLIs will be of great use in section of QPLI below.\\

In fact, there exists a more essential way of proving Theorem 1 using MZV decomposition, which is implicitly outlined in \cite{ref4, ref5, ref6, ref18}. It is of great complexity compared to our methods above. Since we focus on elementary proofs, its explanation will be omitted. Moreover, we establish all QLI/QPLI results, and all but 10 QES results using MZV-free methods below. Only in the last part of solution, $10$ QESs will be solved by applying MZV theory.\\

\section{QLI}

There are 14, 50, 129 QLIs for weight 2, 3, 4 respectively. \\

\subsection{Logsine integrals and contour deformation}

\subsubsection{LSI, PLSI, ITF}
We introduce logsine integrals (abbr. LSIs) as follow:

$$LSI(a,b,c;d)=\int_0^{\frac{\pi }{2}} x^a  \log ^b(2 \sin (x)) \log ^c(2 \cos (x))\cot ^d(x) \, dx$$\\
Here $a, b, c\geq0, d=0, \pm1$. The weight of an LSI is defined as $W=a+b+c+1$.\\

\noindent \textbf{Proposition 8.} All convergent LSIs with weight $W \leq 5$ are generated by $A_W \cup \pi A_{W-1}$ over $\mathbb{Q}$. Here $\pi S=\{\pi  s|s \in S \}, A_0=\{1\}$.\\

\noindent Proof. Consider $F_0(x)=\frac{l (0;x)^k l(1;x)^m l(2;x)^n }{f(p;x)}$ with all logs principal, integrate along the contour encircling the upper unit disc. Map the $(-1,0)$ part back to $(0,1)$ via $t \to -x$, parametrize the semicircular part by $t\to e^{2 i x}$. Since $\log (-x)=\log (x)+i \pi$ on upper $(-1,0)$ and $\log (1+e^{2 i x})=\log (2 \cos (x))+i x,\log (1-e^{2 i x})=\log (2 \sin (x))+i \left(x-\frac{\pi }{2}\right)$ on $(0,\frac{\pi}{2})$, we have that:

$$\int_0^{\frac{\pi}{2}} F(x) G_l(x) \, dx=LI(k,m,n;l)-\sum_{j=0}^m \binom{m}{j}(i \pi)^j LI(n,m-j,k;2-l)$$\\
Where $F(x)=(2 i x)^m \left(\log (2 \sin (x))+i \left(x-\frac{\pi }{2}\right)\right)^k (\log (2 \cos (x))+i x)^n, G_0(x)=\cot (x)+i, G_1(x)=-2 i$ and $G_2(x)=\tan (x)-i$. Taking $\Re/\Im$ yields a LI-LSI relation. Now consider all $F_0$ to build up the equation system, from which all LSIs desired are solved modulo known LIs. $\hfill\square$\\

Clarification: There's no complete table for LSIs in this article. Readers may compute them by using Prop. 8 themselves. Several of LSIs will be used later in section 8. Also note that this solution naturally extends to higher weights.\\

Though Prop. 8 is already a complete solution, we may also apply Fourier expansion.  Since \\
$$\sum _{k=1}^{\infty } \frac{(-1)^{k-1} \cos (2 k x)}{k}=\log (2 \cos (x)), -\sum _{k=1}^{\infty } \frac{\cos (2 k x)}{k}=\log (2 \sin (x))$$\\
 hold on $(0,\frac{\pi}{2})$ \cite{ref13}, we see immediately that $\int_0^{\frac{\pi }{2}} x^n \log (2 \sin (x)) \, dx$ or $\int_0^{\frac{\pi }{2}} x^n \log (2 \cos (x)) \, dx$ are given closed-forms using elementary results of $\int_0^{\frac{\pi }{2}} x^m e^{i n x} \, dx$ and Fubini theorem. This is naturally extended to the integral on $(0,\frac{\pi}{4})$, see subsection 5-12 for applications.\\

Another way of solving LSIs is to differentiate the following formula \cite{ref13} w.r.t $v, a$:\\
$$\int_0^{\frac{\pi }{2}} \cos (a x) \cos ^{v-1}(x) \, dx=\frac{\pi }{2^v v B\left(\frac{a+v+1}{2} ,\frac{-a+v+1}{2} \right)}$$\\
By which all LSIs of form $\int_0^{\frac{\pi }{2}} x^{2 n} \log ^m(2 \cos (x)) \, dx, \int_0^{\frac{\pi }{2}} x^{2 n+1} \tan (x) \log ^m(2 \cos (x)) \, dx$ are solved. Similarly, differentiate the following formula w.r.t $p$ (for its proof, substitute $\tan(x)\to t$ then refer to Mellin transform subsection below):\\
$$\int_0^{\frac{\pi }{2}} x \tan ^p(x) \, dx=\frac{1}{4} \pi  \csc \left(\frac{\pi  p}{2}\right) \left(\psi ^{(0)}\left(1-\frac{p}{2}\right)-2 \psi ^{(0)}(1-p)-\gamma \right)$$\\
We obtain closed-forms of $\int_0^{\frac{\pi }{2}} x \log ^n(\tan (x)) \, dx$ and $\int_0^{\frac{\pi }{2}} x \cot (x) \log ^n(\tan (x)) \, dx$, which yield LSI relations after expansion.\\

We define poly-logsine integrals (abbr. PLSIs) as an integral on $(0,\frac{\pi}{2})$, with its integrand $F G$ with $F=x^a  \log ^b(2 \sin (x)) \log ^c(2 \cos (x))\cot ^d(x)$, $G$ a product whose terms belong to the set \\
$$\left\{\text{Cl}_2(2 x), \text{Cl}_2(4 x),\text{Li}_2\left(-\tan ^2(x)\right), \text{Ti}_2(\tan (x))\right\}$$\\
  They arise in parametrizations of $L(2,k;x),k=1,2,5,6$ via substitution $x\to e^{2 i t}$ in PLIs \cite{ref14}: \\
$$\text{Li}_2\left(x\right)=\text{Li}_2\left(e^{2 i t}\right)=i \text{Cl}_2(2 t)+t^2-\pi  t+\frac{\pi ^2}{6}$$

$$\text{Li}_2\left(-x\right)=\text{Li}_2\left(-e^{2 i t}\right)=i \left(\frac{1}{2}\text{Cl}_2(4 t)-\text{Cl}_2(2 t)\right)+t^2-\frac{\pi ^2}{12}$$

$$\text{Li}_2\left(\pm \frac{1-x}{x+1}\right)=\text{Li}_2\left(\mp i \tan(t)\right)=\frac{1}{4}\text{Li}_2\left(- \tan^2(t)\right)\mp i \text{Ti}_2 \left(\tan(t)\right)$$\\
Deforming contours like before, PLSIs can be evaluated.\\

Furthermore, by $x\to \tan(t)$ and repeated integration by parts, it's possible to convert the inverse tangent family (abbr. ITF) \\
$$\int_0^{\infty } \frac{\tan ^{-1}(x)^k \log ^m\left(x^2+1\right)}{x^l \left(x^2+1\right)^n} \, dx$$\\
 to combination of LSIs hence obtain their closed-forms. Differentiating $\int_0^{\infty } \frac{x^{a-1}}{(x+1)^{a+b}} \, dx=B(a,b)$ or integrating $\frac{\log ^m(1+i z)}{z^n}$ along the upper plane also generate relations between them.\\

\subsubsection{QLSI and QLI}
Quadratic logsine integrals (abbr. QLSIs) are defined as:\\
$$QLSI(a,b,c)=\int_0^{\frac{\pi }{4}} x^a  \log ^b(2 \sin (x)) \log ^c(2 \cos (x)) dx$$\\
The name 'logsine' originate from \cite{ref10}. QLSIs are involved in solution of certain QLIs. For instance, by $x\to \tan(t)$ we have \\
$$\int_0^1 \frac{\log ^m(x) \log ^n\left(x^2+1\right)}{x^2+1} \, dx=(-2)^n \int_0^{\frac{\pi }{4}} \log ^m(\tan (t)) \log ^n(\cos (t)) \, dt$$\\
In our case $m+n\leq 3$ so that brute force may be used. Expand RHS using $\log(\tan(t))=\log(\sin(t))-\log(\cos(t))$, deform contour via $z\to e^{2it}$ like before, QLSIs are converted into complex log integrals $\int_1^i \cdots dz$, which are completely solvable by Mathematica (using $F1\sim F4$ in subsection 2-1). This method is also applicable for another type of QLI: $\int_0^1 \frac{\log^m(x) \tan^{-1}(x)^n}{x^2+1}dx$.  There's also no complete table for QLSIs in this article.\\

\subsection{Inverse tangent family}
Start from a integral in ITF, map $(1, \infty)$ part to (0,1) via $t\to \frac{1}{x}$, after simplification on logs and $\tan^{-1}(x)$ (note $\tan ^{-1}\left(\frac{1}{x}\right)+\tan ^{-1}(x)=\frac{\pi }{2}$) we obtain a QLI relation. Now consider all ITF integrals.\\

\subsection {Polylog special values}
Apparently $\text{Li}_n(\pm 1), \text{Li}_n(\pm i)$ can be represented by $\zeta, \eta, \beta, \psi^{(n)}$. Moreover we have (proved by plugging in special values in formulas of  Prop. 7 and \cite{ref14}):\\
$$\text{Li}_2(1+i)=\frac{\pi ^2}{16}+i \left(C+\frac{1}{4} \pi  \log (2)\right), \text{Li}_2\left(\frac{1}{2}+\frac{i}{2}\right)=i \left(C-\frac{1}{8} \pi  \log (2)\right)+\frac{5 \pi ^2}{96}-\frac{1}{8} \log ^2(2)$$

$$\Re(\text{Li}_3(1+i))=\frac{35 \zeta (3)}{64}+\frac{1}{32} \pi ^2 \log (2), \Re\left(\text{Li}_3\left(\frac{1+i}{2}\right)\right)=\frac{35 \zeta (3)}{64}+\frac{\log ^3(2)}{48}-\frac{5}{192} \pi ^2 \log (2)$$ 

$$\Im\left(\text{Li}_3\left(\frac{1+i}{2}\right)\right)+\Im(\text{Li}_3(1+i))=\frac{7 \pi ^3}{128}+\frac{3}{32} \pi  \log ^2(2)$$

$$\Im(\text{Li}_4(1+i))-\Im\left(\text{Li}_4\left(\frac{1+i}{2}\right)\right)=\frac{1}{64} \pi  \log ^3(2)+\frac{7}{256} \pi ^3 \log (2)$$\\

Below is the only nontrivial pair of special values.\\

\noindent \textbf{Proposition 9.} The following formulas hold: \\
$$\Re(\text{Li}_4(1+i))=-\frac{5 \text{Li}_4\left(\frac{1}{2}\right)}{16}+\frac{97 \pi ^4}{9216}-\frac{5 \log ^4(2)}{384}+\frac{1}{48} \pi ^2 \log ^2(2)$$

$$\Re\left(\text{Li}_4\left(\frac{1}{2}+\frac{i}{2}\right)\right)=\sum _{k=1}^{\infty } \frac{\sum _{l=0}^{\left[\frac{k}{2}\right]} (-1)^l \binom{k}{2 l}}{k^4 2^k}=\frac{5 \text{Li}_4\left(\frac{1}{2}\right)}{16}+\frac{343 \pi ^4}{92160}+\frac{\log ^4(2)}{96}-\frac{5}{768} \pi ^2 \log ^2(2)$$\\
In which $\left[a\right]$ denotes the integer part of $a$.\\

\noindent Proof: Two formulas are equivalent due to inversion formula of $\text{Li}_4(z)$. Consider evaluating $QLI(225;4)$ in 2 ways: \\
$$(1)\ \ QLI(225;4)\to \int_0^{\frac{\pi }{4}} t \log ^2(\tan (t)) \, dt\to \text{QLSIs}$$\\
 Where QLSIs are evaluated by brute force in the subsection 5-1-2. \\
 $$(2)\ \ QLI(225;4)\to \cdots+\sum _{n=1}^{\infty } \frac{(-1)^n H_{2 n}}{n^3}\to \cdots \Re\left(\sum _{n=1}^{\infty } \frac{i^n H_n}{n^3}\right)$$ \\
 Where we have used double expansion like in subsection 4-7, then a simplified version of Prop. 4(3). Equating 2 forms generated by 2 methods (use Mathematica for simplicity) gives the closed-from of $\Re(\text{Li}_4(1+i))$ and completes the proof.  $\hfill\square$\\
 
\subsection{Brute force}
\subsubsection{Class 1}
By writing $\log \left(x^2+1\right), 2i \tan ^{-1}(x)$ as $\log (1+i x)\pm\log(1-i x)$, all QLIs with $W\leq 3$ are decomposed into several integrals of form $\int_0^1 \frac{\log (a+x) \log (b+x)}{c+x} \, dx$, which are trivial due $F3$, subsection 2-1.\\

\subsubsection{Class 2}
Moreover, due to $F1, F2, F4$, subsection 2-1, all integrals of form \\
$$\int \frac{\log ^2(A(x)) \log (B(x))}{A(x)} \, dx, \int \frac{\log ^2(A(x)) \log (B(x))}{B(x)} \, dx, \int \frac{\log ^n(C(x))}{g(p;x)} \, dx$$\\
 (here $A(x), B(x)\in Z=\{1,x,1\pm x,1\pm i x\}$, $C(x)=\frac{P(x)}{Q(x)}$ with $P(x), Q(x) \in Z$, $p\not=1$) can be calculated by brute force as well. For former two we simply use $F2, F4$, while the latter boils down to transformation $z\to C(x)$, partial fraction decomposition outside the log (since the inverse $D=C^{-1}$ is also a fractional transformation, the outside $R(z)=\frac{D'(z)}{g(p;D(z))}$ is rational), and usage of $F1$. Consider all possible $A, B, C, p$, expand the expression, take $\Re/\Im$ to get QLI relations. Also some of the QLIs are solved directly in this way, say $QLI(555;3)$. One may try on these examples: \\
 $$\int_0^1 \frac{\log (1-i x) \log ^2(x+1)}{1-i x} \, dx,\int_0^1 \frac{\log ^2(1+i x) \log (1-x)}{1+i x} \, dx $$
 
  $$\int_0^1 \frac{ \log ^3\left(\frac{1-x}{1+i x}\right)}{x} \, dx,\int_0^1 \frac{x \log ^3\left(\frac{x}{1+i x}\right)}{x^2+1} \, dx, \int_0^1 \frac{\log ^3\left(\frac{1+i x}{1-i x}\right)}{x+1} \, dx$$\\
All can be done without much effort with the help of Mathematica and subsection 5-3.\\

\subsection {General formulas} 
Evidently the following are expressible in terms of $\zeta, \eta, \beta, \text{Li}_k$ for arbitrary $n$:\\
 $$\int_0^1 \frac{\log ^n(x)}{x^2+1} \, dx,\int_0^1 \frac{x \log ^n(x)}{x^2+1} \, dx,\int_0^1 \frac{x \log ^n\left(x^2+1\right)}{x^2+1} \, dx,\int_0^1 \frac{\tan ^{-1}(x)^n}{x^2+1} \, dx$$\\
By series expansion and calculating primitives. Based on subsection 2-2 and these four, we may also evaluate the following in a general sense, using $x\to x^2$ and integration by parts:\\
$$\int_0^1 \frac{\log ^n\left(x^2+1\right)}{x} \, dx, \int_0^1 \frac{x \log (x) \log ^n\left(x^2+1\right)}{x^2+1} \, dx, \int_0^1 \frac{\log \left(x^2+1\right) \log ^n(x)}{x} \, dx$$

$$\int_0^1 \frac{\tan ^{-1}(x) \log ^n(x)}{x} \, dx, \int_0^1 \frac{\log ^2\left(x^2+1\right) \log ^{2 n}(x)}{x} \, dx, \int_0^1 \frac{x \log \left(x^2+1\right) \log ^{2 n-1}(x)}{x^2+1} \, dx$$\\
One may write down exact formulas themselves. Now use them to solve QLIs directly.\\

\subsection {Feynman's trick on Mellin transform}

\noindent \textbf{Proposition 10.} Denote \\
$$f(s)=\frac{1}{4} \pi  \sec \left(\frac{\pi  s}{2}\right) \left(\psi ^{(0)}\left(1-\frac{s}{2}\right)-\psi ^{(0)}\left(\frac{1}{2}\right)\right)$$ 

$$g(s)=2 \beta \left(\frac{1-s}{2}\right) \csc \left(\frac{\pi  s}{2}\right)+2 \beta \left(\frac{2-s}{2}\right) \sec \left(\frac{\pi  s}{2}\right)-2 \pi  \csc (\pi  s)$$
Then we have\\
$$\int_0^1 \frac{\tan ^{-1}(x) \log ^{2 n+1}(x)}{x^2+1} \, dx=\frac{1}{2} \underset{s\to 0}{\text{lim}}\frac{\partial ^{2 n+1}f(s)}{\partial s^{2 n+1}}-\frac{1}{4} \pi  (2 n+1)! \beta (2 n+2)$$

$$\int_0^1 \frac{x \tan ^{-1}(x) \log ^{2 n}(x)}{x^2+1} \, dx=\frac{1}{2} (2 n)! \beta (2 n+2)+\frac{\pi  \left(1-2^{-2 n}\right) (2 n)! \zeta (2 n+1)}{2^{2 n+3}}-\frac{1}{2} \underset{s\to 0}{\text{lim}}\frac{\partial ^{2 n}f(s)}{\partial s^{2 n}}$$

$$\int_0^1 \frac{\tan ^{-1}(x) \log ^{2 n}(x)}{x+1} \, dx=\frac{1}{2} (2 n)! \beta (2 n+2)+\frac{1}{4} \pi  \left(1-2^{-2 n}\right) (2 n)! \zeta (2 n+1)-\frac{1}{16} \pi  \underset{s\to 0}{\text{lim}}\frac{\partial ^{2 n}g(s)}{\partial s^{2 n}}$$\\

\noindent Proof. Consider Mellin transform of $\frac{\tan ^{-1}(x)}{x+1},\frac{\tan ^{-1}(x)}{x^2+1}$. Take the latter as the example: let \\
$$M(s)=\int_0^{\infty } \frac{x^{s-1} \tan ^{-1}(x)}{x^2+1} \, dx, J(a,s)=\int_0^{\infty } \frac{x^{s-1} \tan ^{-1}(a x)}{x^2+1} \, dx$$ \\
then by Feynman's trick and proof of  Prop. 2(6) we get\\
$$\frac{\partial J}{\partial a}=\int_0^{\infty } \frac{x^s}{\left(x^2+1\right) \left(a^2 x^2+1\right)} \, dx=\frac{\pi  \left(1-a^{1-s}\right) \sec \left(\frac{\pi  s}{2}\right)}{2 \left(1-a^2\right)}$$ \\
Integrate w.r.t $a$ using the following (see \cite{ref13}) $$\int_0^1 \frac{x^u-x^v}{1-x} \, dx=\psi ^{(0)}(v+1)-\psi ^{(0)}(u+1)$$ we get $J$. Differentiate w.r.t $s$, substitute $t\to \frac{1}{x}$ for $(1,\infty)$ part on RHS and assign appropriate values of $s$ gives the first two formulas. The last one is similar, see \cite{ref15}. $\hfill\square$\\

Now use them to solve QLIs directly, in our case $QLI(225;5), QLI(225;3)$.\\

\subsection {Integration by parts} 
In spirit of subsection 2-3, generalize $F$ to $\prod _{m=1}^5 k(m;x)^{A(m)}$. Now consider all possible $F$.\\

\subsection {Fractional transformation} 
Similar to subsection 2-4, substitute $x\to \frac{1-x}{x+1}$ in integral of $F= \frac{\prod _{m=1}^5 k(m;x)^{A(m)}}{g(p;x)} $, where $p\not=2$. 
Using $$\tan ^{-1}(x)+\tan ^{-1}\left(\frac{1-x}{x+1}\right)=\frac{\pi }{4}, \log \left(\left(\frac{1-x}{x+1}\right)^2+1\right)=\log \left(\frac{2 \left(x^2+1\right)}{(x+1)^2}\right)$$
we may split the integrand to obtain QLI relations (after possible decomposition $\frac{1-x}{(1+x)(1+x^2)}\to \frac{1}{1+x}-\frac{x}{1+x^2}$). Now consider all possible $F$.\\

\subsection {Power substitution} 
For each LI, substitute $x\to t^2$, then expand $\log^n(1-t^2)$ via binomial formula (and possibly decompose $\frac{t}{1-t^2}$) to get QLI relation. Now consider all LIs.\\

\subsection {Weierstrass substitution} 
For each LSI, substitute $x\to2 \tan^{-1}(t)$. Since $(\sin(x), \cos(x))\to \left(\frac{2t}{1+t^2}, \frac{1-t^2}{1+t^2}\right)$, $dx\to 2\frac{dt}{1+t^2}$, we get a QLI relation after simplification. Now consider all LSIs.\\

\subsection {Parametric integration} 
This method only generates one relation. By (keyhole) contour integration we may show \\
$$\int_0^{\infty } \frac{\log ^2(x)}{\left(x^2+1\right) (a+x)} \, dx=\frac{\pi ^2 a}{8 \left(a^2+1\right)}-\frac{\log ^3(a)}{3 \left(a^2+1\right)}-\frac{\pi ^2 \log (a)}{3 \left(a^2+1\right)} (a>0)$$\\
 Integrate $a$ on (0,1), expand the RHS integrals, $\int_0^{\infty } \frac{\log (x+1) \log ^2(x)}{x^2+1} \, dx$ is computed. Now map the $(1, \infty)$ part to $(0,1)$.\\

\subsection {Fourier expansion} 
We may solve $\int_0^1 \frac{\tan ^{-1}(x)^n}{x} \, dx$ and $\int_0^1 \frac{x \tan ^{-1}(x)^n}{x^2+1} \, dx$ by Fourier expansion after $x\to \tan(t)$. The former becomes $\int_0^{\frac{\pi }{4}} t^n (\tan (t)+\cot (t)) \, dt$, then $\int_0^{\frac{\pi }{4}} t^{n-1} \log (\sin (t)) \, dt, \int_0^{\frac{\pi }{4}} t^{n-1} \log (\cos (t)) \, dt$ after integration by parts, which are completely solved using Fourier series in subsection 5-1-1. The latter boils down to $\int_0^{\frac{\pi }{4}} t^n \tan (t) \, dt$ either. See also \cite{ref15} for evaluation of $\int_0^1 \frac{\log \left(1-x^2\right) \tan ^{-1}(x)^2}{x^2+1} \, dx$ using Fourier expansion twice.\\

\subsection {Beta derivatives} 
\subsubsection{Class 1}
Similar to subsection 2-5, all integrals of form $\int_0^1\frac{\log ^m(x) \log ^n\left(1-x^4\right)}{g(p;x)}dx$ ($p=1, 2, 3, 4, 5$) are solvable. For instance, when $p=3$ we write $\frac{1}{x+1}=\frac{-x^3+x^2-x+1}{1-x^4}$, substitute $x^4\to t$ then split it into 4 integrals expressible via Beta derivatives. On the other hand, we may expand $\log^n(1-x^4)=(\log(1-x)+\log(x+1)+\log(x^2+1))^n$ to obtain a QLI relation. Now consider all $m, n, p$.\\

\subsubsection{Class 2}
Another variation is to multiply both sides of the equality \\
$$\frac{\log \left(x^2+1\right)}{1\mp x}=\frac{\left(x^3\pm x^2\pm x+1\right) \log \left(1-x^4\right)}{1-x^4}-\frac{(\pm x+1) \log \left(1-x^2\right)}{1-x^2}$$\\
 with $\log^m(x)$ and integrate. Therefore $\int_0^1 \frac{\log \left(x^2+1\right) \log ^m(x)}{1\mp x} \, dx$, in our case $QLI(224;1), QLI(224;3)$, are evaluated by Beta derivatives using $x^2\to t, x^4\to t$ on RHS.\\

\subsection {Contour integration} 
\subsubsection{Class 1}
Similar to subsection 2-6, consider integrating \\
$$F_1(z)=\frac{\log^k(1-z) \log^m(z)\log^n(z+1)}{1+z^2}\ \ \text{or}\ \  F_2(z)=\frac{\log^k(1-z) \log^m(z)\log^n(z+1)}{z(1+z^2)}$$\\
 along large upper semicircular contour with all logs principal. Residue theorem yields \\
 $$\int_{-\infty+0i}^{\infty+0i}F_k(z)dz=2\pi i \text{Res}(F_k, i)$$\\
  Split LHS into 4 parts by branch points $0, \pm 1$, use $x\to\pm t, x\to \pm \frac{1}{t}$ to map them to $(0,1)$, after expansion and taking $\Re/\Im$ we obtain a linear relation. Now consider all possible $F_1, F_2$.\\

\subsubsection{Class 2}
Alternatively, consider integrating \\
$$F(z)=\frac{\log^k(1-z) \log^m(z)\log^n(1+z) \log^s(1+iz) \log^t(1-iz)}{z(1+z)}$$\\
 along large quartercircular contour encircling the first quadrant, with all logs principal. Cauchy theorem yields $\int_0^{\infty}F(z)dz-\int_0^{i\infty}F(z)dz=0$. Split the integral into 4 parts by $0, 1, i$, use $x\to t, x\to \frac{1}{t}, x\to i t, x\to \frac{i}{t}$ to map them to $(0,1)$, after expansion and taking $\Re/\Im$ we obtain a linear relation. Now consider all possible $F$.\\

Take $F(z)=\frac{\log (1-z) \log (1+z) \log (1+i z)}{z (1+z)}$ as an example, after routine simplification we have

{\scriptsize $$\int_0^1 -(\frac{x}{x^2+1}-\frac{i}{x^2+1}) (\log (1-x)-\log (x)+i \pi ) \left(\frac{1}{4} \log ^2\left(x^2+1\right)-\log (x) \log \left(x^2+1\right)+\log ^2(x)+\tan ^{-1}(x)^2-\pi  \tan ^{-1}(x)+\frac{\pi ^2}{4}\right)$$ $$-(\frac{1}{x}-\frac{i}{x^2+1}-\frac{x}{x^2+1}) \log (1-x) \left(\frac{1}{4} \log ^2\left(x^2+1\right)+\tan ^{-1}(x)^2\right)+\left(\frac{1}{x}-\frac{1}{x+1}\right) \log (1-x) \log (x+1) \left(\frac{1}{2} \log \left(x^2+1\right)+i \tan ^{-1}(x)\right) $$ $$+\frac{1}{x+1}(\log (1-x)-\log (x)-i \pi ) (\log (x+1)-\log (x)) \left(\frac{1}{2} \log \left(x^2+1\right)-\log (x)-i \tan ^{-1}(x)+\frac{\pi  i}{2}\right)\, dx=0$$}\\
Now expand this expression and take $\Re/\Im$ to obtain 2 QLI relations.\\

Note: In fact we have a `Beta' use of contour integration, that is, integrate $\frac{\log ^m(z) \log ^n(1-z)}{z (1-z)}$ along $(0,\infty)\to(0,0)\to(1,0)\to(1,\infty)$. Since it only generates superfluous relations, we won't elaborate more here.\\

\subsection{Series expansion} 
The procedure of multiple expansion in subsection 4-7 can be naturally extended to QLIs, since $\log \left(x^2+1\right),\tan ^{-1}(x),\frac{1}{x^2+1}, \frac{x}{x^2+1}$, and products like $\log \left(x^2+1\right) \tan ^{-1}(x), \tan ^{-1}(x)^2, \frac{\tan ^{-1}(x)}{x^2+1}$ have elementary or (quadratic) harmonic power series. Finally, if not discarded, set $x\to i$ in all generating functions in Prop. 4, take $\Re/\Im$ to see whether the resultant sum is included.\\

Take $QLI(115;2)$ for an example. After the procedure we get quadratic sums $\sum _{k=0}^{\infty } \frac{(-1)^k H_{2 k+1}^{(2)}}{(2 k+1)^2}$, $\sum _{k=0}^{\infty } \frac{(-1)^k \left(H_{2 k+1}\right){}^2}{(2 k+1)^2}$, which are directly solved by setting $x\to i$ in Prop. 4(1)(2), taking imaginary parts and using Prop. 9 for simplification. Moreover, using Prop 4. and Prop. 5(3)(4), $QLI(125;2), QLI(114;2), QLI(124;2)$ can be evaluated similarly. \\

\subsection{Solution}
Solving the equation system of QLI relations generated by this section, all 193 QLIs with weight $\leq4$ are given closed-forms. Moreover, one may observe there are essentially two types of closed-forms in QLIs with fixed weight. Take $W=4$, then 68 QLIs lie in the subspace spanned by 8 constants in the extended basis $B_4$ over $\mathbb{Q}$, in which those five in $A_4$ are included. The remaining 61 QLIs are completely expressible by other 8 constants in $B_4$. This phenomenon of separation will extend to QPLIs.\\

\section{QPLI}
There are 20, 152 QPLIs with weight 3, 4 respectively.\\

\subsection{Brute force}
Some QPLIs can be calculated by brute force, for instance $\int_0^1 \frac{\text{Ti}_2(x) \tan ^{-1}(x)}{x} \, dx$. For each QPLI, check by Mathematica whether the premitive is solvable and evaluate it if so.\\

\subsection{General formulas}
See \cite{ref15} for proof of the following:\\

\noindent \textbf{Lemma 4.} $\int_0^1 \frac{\text{Ti}_2(x) \log ^{2 n-1}(x)}{x+1} \, dx,\int_0^1 \frac{\text{Li}_{n+1}(-x) \log ^n(x)}{x^2+1} \, dx$ and $\int_0^1 \frac{x \log ^n(x) \text{Li}_{n+1}(\pm x)}{x^2+1} \, dx$ are reducible to $\zeta, \beta$ values.\\

Thus QPLIs of these types are solved directly.\\

\subsection{Integration by parts}
Lifting up the denominator, all QPLIs with $W=3$ are trivial. Generally, in spirit of subsection 2-3, generalize $F$ to $\prod _{m=1}^5 k(m;x)^{A(m)} \prod _{k=1}^8 \prod _{n=2}^N L(n,k;x)^{b(n,k)}$ (in $W=4$ case it is of form $L(4,a;x), L(3,a;x)k(b;x), L(2,a;x)k(b;x)k(c;x)$ or $L(2,a;x)L(2,b;x)$). For those involving $L(n,a;x), a\in\{3,5,6\}$, after partial fractions another integration by parts lifting $\log(1-x)$ is needed to obtain a regular relation modulo QLIs. Now consider all possible $F$.\\

Take $F=\text{Li}_2\left(\frac{1-x}{2}\right) \text{Ti}_2(x)$ as an example, after the extra procedure we obtain a relation between $QPLI4(8;1;3), QPLI4(3;5;2)$.\\

\subsection{Fractional transformation}
Similar to subsection 2-4, apply $x\to\frac{1-x}{x+1}$ in integral of $F=\frac{\prod _{m=1}^5 k(m;x)^{A(m)} \prod _{k=1}^8 \prod _{n=2}^N L(n,k;x)^{b(n,k)}}{g(p;x)}$ (in $W=4$ case it is of form $\frac{L(2,a;x)k(b;x)}{g(p;x)}$ or $\frac{L(3,a;x)}{g(p;x)}, a\in \{1,2,5,6\}$). For $p=2$, after partial fractions, possibly another integration by parts lifting $\log(1-x)$ is needed to obtain a regular relation. Now consider all possible $F$.\\

\subsection{Polylog identities}
Like in subsection 4-6, multiply both sides of Prop. 7(1) with $F=\frac{\log^m(1-x^2)\log^n(1+x^2)\log^k(x)}{x}$ or $\frac{x\log^m(1-x^2)\log^n(1+x^2)\log^k(x)}{1\pm x^2}$, Prop. 7(9) with $F=\frac{1}{x^2+1}$ or $\frac{x}{x^2+1}$, Prop. 7(5)(6) with $F=\frac{k(m;x)}{g(p;x)}$, integrate on (0,1) to obtain relations modulo known QLIs. Now consider all 4 formulas (1)(5)(6)(9) and all corresponding possible $F$.\\

\subsection{Power substitution}
Similar to subsection 5-9, for each PLI containing only $\text{Li}_n(\pm x)$ and log terms, substitute $x\to t^2$ then use Prop. 7(1) for decomposition, we obtain a QPLI relation (for instance taking $F=\frac{\text{Li}_2(x) \log (1-x)}{x+1}$ yields a relation between $QPLI4(1/2;1/3;5)$). Now consider all suitable PLIs.\\

\subsection{Valean's master formula}
\subsubsection{Class 1}
In \cite{ref16}, by calculating Cauchy products of polylogarithms based on subsection 3-3, the following is proved whenever $|a|\leq1$:\\
$$\int_0^1\frac{\ln(t)\text{Li}_2(t)}{1-at}\ dt=\frac{\text{Li}_2^2(a)}{2a}+3\frac{\text{Li}_4(a)}{a}-2\zeta(2)\frac{\text{Li}_2(a)}{a}$$\\
It naturally extends to the plane. Letting $a\to i$ and taking $\Re/\Im$, $QPLI4(1;2;4/5)$ are computed.\\

\subsubsection{Class 2}
Moreover, by using $a\to \frac{i}{i-1}$, $\int_0^1 \frac{\text{Li}_2(1-x) \log (1-x)}{1-i x} \, dx=\frac{1}{1-i}\int_0^1 \frac{\text{Li}_2(x) \log (x)}{1-\frac{i}{i-1}x} \, dx$ is computable due to subsection 5-3. By Prop. 7(2) the corresponding $\int_0^1 \frac{\text{Li}_2(x) \log (1-x)}{1-i x} \, dx$ is also evaluated modulo QLIs, thus giving closed-forms of $QPLI4(1;1;4/5)$.\\

\subsection{Contour integration}
\subsubsection{Class 1}
Since we only discuss the weight 4 case, integrands will be restricted. Similar to subsection 5-14, consider integrating \\
$$F_1(z)=\frac{\text{Li}_2(\pm z)\log(1+a z)}{1+z^2}\ \ \text{or}\ \  F_2(z)=\frac{\text{Li}_2(\pm z)\log(1+a z)}{z(1+z^2)}$$\\
along large upper semicircular contour (here logs and polylogs are principal, $a\in\{\pm1, -i\}$), apply residue theorem, map each part to (0,1), take $\Re/\Im$ to obtain relations. Now consider all possible $F_1, F_2$.\\

For instance, take $F_2(z)=\frac{\text{Li}_2(z) \log (1-z)}{z \left(z^2+1\right)}$, after routine simplification we have that:\\

{\tiny $$\int_0^1 -\frac{z }{z^2+1}(\log (z+1)-\log (z)) \left(-\text{Li}_2(-z)-\frac{1}{2} \log ^2(z)-\frac{\pi ^2}{6}\right)+\frac{z}{z^2+1} (\log (1-z)-\log (z)-i \pi ) \left(-\text{Li}_2(z)-\frac{1}{2} \log ^2(z)-i \pi  \log (z)+\frac{\pi ^2}{3}\right)$$ $$-\left(\frac{1}{z}-\frac{z}{z^2+1}\right) \text{Li}_2(-z) \log (z+1)+\left(\frac{1}{z}-\frac{z}{z^2+1}\right) \text{Li}_2(z) \log (1-z) \, dz=2 \pi  i \text{Res}\left(F_2(z),i\right)=-\frac{1}{4} i \pi ^2 C+\frac{1}{2} \pi  C \log (2)+\frac{\pi ^4}{192}+\frac{1}{96} i \pi ^3 \log (2)$$}\\

\subsubsection{Class 2}
Alternatively, consider integrating\\
 $$F(z)=\frac{\text{Li}_2(az) \log(1+bz)}{z(1+z)}$$ \\
 along large quartercircular contour encircling the first quadrant (here logs, polylogs are principal, $a,b\in\{\pm1,\pm i\}$ ), apply Cauchy theorem, map each part to (0,1), take $\Re/\Im$ to obtain relations. Now consider all possible $F$.\\

For instance, take $F(z)=\frac{\text{Li}_2(z) \log (1-z)}{z (z+1)}$, after routine simplification we have that:\\

{\scriptsize$$\int_0^1 -\left(\frac{z}{z^2+1}-\frac{i}{z^2+1}\right) \left(\frac{1}{2} \log \left(z^2+1\right)-\log (z)+i \tan ^{-1}(z)-\frac{i \pi }{2}\right)\left(-\frac{\text{Li}_2\left(-z^2\right)}{4}+i \text{Ti}_2(z)-\frac{1}{2} \log ^2(z)-\frac{1}{2} \pi  i \log (z)-\frac{\pi ^2}{24}\right)$$ $$-\left(-\frac{z}{z^2+1}-\frac{i}{z^2+1}+\frac{1}{z}\right) \left(\frac{\text{Li}_2\left(-z^2\right)}{4}+i \text{Ti}_2(z)\right) \left(\frac{1}{2} \log \left(z^2+1\right)-i \tan ^{-1}(z)\right)$$ $$+\frac{1}{z+1}(\log (1-z)-\log (z)-i \pi ) \left(-\text{Li}_2(z)-\frac{1}{2} \log ^2(z)-i \pi  \log (z)+\frac{\pi ^2}{3}\right)+\left(\frac{1}{z}-\frac{1}{z+1}\right) \text{Li}_2(z) \log (1-z)\, dz=0$$}\\

Note: In simiplifications of relations above $\log (1+i z)=\frac{1}{2} \log \left(z^2+1\right)+i \tan ^{-1}(z)$, $\text{Li}_2(i z)=\frac{1}{4}\text{Li}_2\left(-z^2\right)+i \text{Ti}_2(z)$ are frequently used. Moreover, due to Prop. 7(4) and definition of principle branches, we have $\text{Li}_2(z)=-\text{Li}_2\left(\frac{1}{z}\right)-\frac{1}{2} \log ^2(z)+i \pi  \log (z)+\frac{\pi ^2}{3}$ on upper side of $(1,\infty)$, here $\text{Li}_2\left(\frac{1}{z}\right)=\sum _{k=1}^{\infty } \frac{1}{k^2 z^k}$. Also reminds that $\text{Ti}_2(z)-\text{Ti}_2\left(\frac{1}{z}\right)=\frac{1}{2} \pi  \log (z)$ for $z>0$.\\

\subsection{Series expansion}

\subsubsection{Class 1}
By Prop. 3(5), 5(2) we have\\
$$\int_0^1 \frac{\text{Li}_2(i x) \log (1-x)}{x (1-i x)} \, dx=-\sum _{n=1}^{\infty } \frac{i^n H_n H_n^{(2)}}{n}$$\\
The sum is solved by Prop. 4(4) and 9. Thus, after partial fraction decomposition, expanding $\text{Li}_2(ix)$ and taking real, imaginary parts on LHS, relations between $(QPLI4(7;1;4), QPLI4(8;1;5))$, $(QPLI4(7;1;5), QPLI4(8;1;4))$ are deduced, modulo known $QPLI4(7/8;1;2)$ (since by lifting up $\text{Li}_2(x)$ then $\log(x)$ they boils down to $QPLI4(1;4/5;2)$ then $QPLI4(1;2;4/5)$ modulo QLIs, which are solved in subsection 6-7).\\

\subsubsection{Class 2}
Moreover, using Prop. 3(5) again gives \\
$$\int_0^1 \frac{\text{Li}_2\left(-x^2\right) \log (1-x)}{x^2+1} \, dx=\sum _{k=0}^{\infty } \frac{(-1)^{k-1} H_{2 k+1} H_k^{(2)}}{2 k+1}$$\\
Substitute $H_k^{(2)}=-\int_0^1 \frac{\left(1-x^k\right) \log (x)}{1-x} \, dx$ in the sum, by Fubini we have it equals to
$$\int_0^1 \frac{\log (x) \left(\Im(f(i))-\frac{\Im\left(f\left(i \sqrt{x}\right)\right)}{\sqrt{x}}\right)}{1-x} \, dx=\int_0^1 \frac{4 x \log (x) \left(\Im(f(i))-\frac{\Im(f(i x))}{x}\right)}{1-x^2} \, dx$$
Where $f(x)=\sum _{k=1}^{\infty } \frac{H_k x^k}{k}=\text{Li}_2(x)+\frac{1}{2} \log ^2(1-x)$ by Prop. 3(6). Extracting imaginary parts using formulas in subsection above, it is furtherly transformed into
$$\int_0^1 \left(\frac{4 C x \log (x)}{1-x^2}-\frac{4 \text{Ti}_2(x) \log (x)}{1-x^2}-\frac{\pi  x \log (2) \log (x)}{2 \left(1-x^2\right)}+\frac{2 \log \left(x^2+1\right) \log (x) \tan ^{-1}(x)}{1-x^2}\right) \, dx$$
Modulo LI/QLIs and partial fractions it boils down to $QPLI4(8;2;1/3)$, then $QPLI4(1/2;5;2)$ and $QPLI4(1/2;2;4)$ by lifting up $\text{Li}_2(\pm x)$ then $\log(x)$. Since $QPLI4(1/2;2;4)$ are evaluated in subsection 6-2, 6-7, the original $QPLI4(7;1;4)$ is solved directly.\\

By taking real parts, $QPLI4(7;1;5)$ is dealt in a similar way.\\

\subsection{Double integration}
\subsubsection{Class 1}
The elementary identity\\
 $$-\int_0^1 \frac{x \log (y)}{1-x y} \, dy=\text{Li}_2(x)\ (|x|<1)$$ \\
 produces several linear relations. For instance, apply this to rewrite the expression of $QPLI4(1;1;5)$, after Fubini and calculating $\int_0^1 \frac{x^2 \log (1-x)}{\left(x^2+1\right) (1-x y)} \, dx$ by brute force (here and below, indefinite integration w.r.t $x$ are carried out by Mathematica, easily verified by differentiation), we have:\\
 $$\int_0^1 \frac{x \text{Li}_2(x) \log (1-x)}{x^2+1} \, dx=-\int_0^1\log(y)\left(\int_0^1\frac{x^2\log(1-x)}{(x^2+1)(1-xy)}dx\right)dy$$
 
 $$=-\int_0^1 \log (y) \left(\frac{C}{y^2+1}+\frac{\text{Li}_2\left(\frac{y}{y-1}\right)}{y \left(y^2+1\right)}+\frac{5 \pi ^2 y}{96 \left(y^2+1\right)}-\frac{y \log ^2(2)}{8 \left(y^2+1\right)}-\frac{\pi  \log (2)}{8 \left(y^2+1\right)}\right) \, dy$$\\
Modulo PLI/QLIs, partial fractions and Prop. 7(3) $\left(\text{Li}_2\left(\frac{y}{y-1}\right)\to \text{Li}_2(y)\right)$, the only one remains is $\int_0^1 \frac{y \text{Li}_2\left(y\right) \log (y)}{y^2+1} \, dy$, thus we get a relation between $QPLI4(1;1;5), QPLI4(1;2;5)$.\\

3 more relations between $(QPLI4(2;1;5), QPLI4(2;2;5))$, $(QPLI4(1;1;4), QPLI4(1;2;4))$ and $(QPLI4(2;1;4), QPLI4(2;2;4))$ can be obtained by double integration similarly, possibly with a substitution $x\to-x$ in the identity above and Prop. 7(3).\\

\subsubsection{Class 2}
Alternatively, apply this method to $QPLI4(1;3;5)$ yields:\\
$$\int_0^1\frac{x\text{Li}_2(x)\log(1+x)}{1+x^2}\ dx=-\int_0^1\log(y)\left(\text{log terms}+\frac{\text{Li}_2\left(\frac{y}{1+y}\right)}{y(1+y^2)}-\frac{\text{Li}_2\left(\frac{2y}{1+y}\right)}{y(1+y^2)}\right)\ dy$$\\
With the part involving log terms trivially evaluated via LIs. Modulo PLI/QLIs, partial fractions and Prop. 7(3) $\left(\text{Li}_2\left(\frac{y}{y+1}\right)\to \text{Li}_2(-y)\right)$ and $QPLI4(2;2;5)$ (evaluated in subsection 6-2), the $\int_0^1 \log(y)\frac{\text{Li}_2\left(\frac{y}{1+y}\right)}{y(1+y^2)}dy$ part is solved. $\int_0^1 \log(y) \frac{\text{Li}_2\left(\frac{2y}{1+y}\right)}{y}dy$ is also solved by lifting up $\frac{\log^2(y)}{2}$ so only $J=\int_0^1 y\log(y) \frac{\text{Li}_2\left(\frac{2y}{1+y}\right)}{1+y^2}dy$ remains. For $J$, integrate by parts twice (lifting firstly $\frac{1}{2}\log(1+y^2)$ then $\frac{1}{2}\text{Li}_2(-y^2)$) to eliminate $\text{Li}_2\left(\frac{2y}{1+y}\right)$ term, modulo QLIs it boils down to evaluation of\\
 $$\int_0^1\log\left(\frac{1-y}{1+y}\right)\frac{\text{Li}_2(-y^2)}{y}dy, \int_0^1\log\left(\frac{1-y}{1+y}\right)\frac{\text{Li}_2(-y^2)}{1+y}dy$$ \\
 The former two $\int_0^1\log\left(1\pm y\right)\frac{\text{Li}_2(-y^2)}{y}dy$ are also reduced to subsection 6-2 after integrating by parts twice (lifting firstly $-\text{Li}_2(\mp y)$ then $\log(y)$), moreover $\int_0^1\log\left(1+y\right)\frac{\text{Li}_2(-y^2)}{1+y}dy$ reduces to QLIs by lifting up $\frac{1}{2}\log^2(1+y)$. Thus the final relation is between the original integral and $\int_0^1\log\left(1-y\right)\frac{\text{Li}_2(-y^2)}{1+y}dy$, i.e. $QPLI4(1;3;5)$ and $QPLI4(7;1;3)$.\\

2 analogous relations between $(QPLI4(1;3;4), QPLI4(8;1;3))$, $(QPLI4(2;3;4), QPLI4(8;1;3))$ can be obtained by double integration and repeated integration by parts similarly.\\

\subsubsection{Class 3}
Furthermore, apply this method to one involving imaginary arguments yields:\\
$$\int_0^1 \frac{\text{Li}_2(x) \log (1+i x)}{x^2+1} \, dx=\int_0^1 -\frac{\log (y) }{y^2+1} \left(\text{quadratic log terms}+\text{Li}_2\left(\frac{y}{i+y}\right)-\text{Li}_2\left(\frac{(1+i) y}{i+y}\right)\right)\, dy$$\\
With the quadratic log terms of form $\log(1\pm y), \log(1\pm iy)$, all logs and polylogs are principal. Denote $f(y)=\text{Li}_2\left(\frac{y}{i+y}\right)-\text{Li}_2\left(\frac{(1+i) y}{i+y}\right)$, then $f'(y)=\left(\frac{i}{y^2+1}-\frac{y}{y^2+1}+\frac{1}{y}\right) \log (1-y)$ is of good form. Modulo known LI/QLIs, our only interest is $\int_0^1 \frac{f(y) \log (y)}{y^2+1} \, dy$. Apply partial integration twice (lifting $\tan ^{-1}(y)$ then $\text{Ti}_2(y)$), modulo QLIs again the integral is reduced to $\int_0^1 \text{Ti}_2(y) f'(y) \, dy$, which is a combination of $QPLI4(8;1;2/4/5)$. Now, modulo QLIs again $QPLI4(8;1;2)$ is equivalent to $QPLI4(1;2;4)$ (solved in subsection 6-7) by partial integration lifting $-\text{Li}_2(y)$ then $\log(y)$, while $QPLI4(8;1;4/5)$ are deduced by combining relations in subsection 6-9 (in which the original integral is not used so it's not a circular argument). Thus we've evaluated the original integral i.e. $QPLI4(1;4/5;4)$ directly.\\

2 relations between $(QPLI4(2;4;4), QPLI4(8;3;2/5))$, $(QPLI4(2;5;4), QPLI4(8;3;4))$ can be obtained similarly if we start from $\int_0^1 \frac{\text{Li}_2(-x) \log (1+i x)}{x^2+1} \, dx$ instead.\\

\subsubsection{Class 4}
Now we make use of the generalized \\
$$\text{Li}_3(x)=\int_0^1 \frac{x \log ^2(y)}{2 (1-x y)} \, dy \ (|x|<1)$$\\
 Consider $QPLI4(3;4)$, by Fubini theorem, rational integration w.r.t $x$ and reflection $y\to 1-y$:
$$\int_0^1 \frac{\text{Li}_3\left(\frac{1-x}{2}\right)}{x^2+1} \, dx=\int _0^1\int _0^1\frac{(1-x) \log ^2(y)}{2 \left(x^2+1\right) (2-(1-x) y)}dydx$$ $$=\int_0^1 \frac{(-\pi  y-4 \log (2-y)+\pi +\log (4)) \log ^2(y)}{8 \left(y^2-2 y+2\right)} \, dy=\int_0^1 \frac{\log ^2(1-y) (\pi  y-4 \log (y+1)+\log (4))}{8 \left(y^2+1\right)} \, dy$$\\
Plug in QLI values, the desired $QPLI4(3;4)$ is solved directly.\\

Note that $QPLI4(4;4)$ can be computed in the same way.\\

\subsubsection{Class 5}
Now we make use of a new kernel:\\
 $$\text{Li}_2\left(-x^2\right)=\int_0^1 \frac{4 x^2 y \log (y)}{x^2 y^2+1} \, dy \ (|x|<1)$$\\
Apply this to another involving imaginary arguments yields:\\
$$\int_0^1 \frac{\text{Li}_2\left(-x^2\right) \log (1+i x)}{x^2+1} \, dx=\int _0^1\int _0^1\frac{4 x^2 y \log (1+i x) \log (y)}{\left(x^2+1\right) \left(x^2 y^2+1\right)}dydx$$

{\small $$=\int_0^1 \frac{2 i \log (y)}{1-y^2} \left(\text{quadratic log terms}+\text{Li}_2\left(\frac{y}{y-1}\right)-\text{Li}_2\left(\frac{(1+i) y}{y-1}\right)-\text{Li}_2\left(\frac{y}{y+1}\right)+\text{Li}_2\left(\frac{(1+i) y}{y+1}\right)\right) \, dy$$}\\
By Prop. 7(3) and known LI/PLI/QLIs, the part involving quadratic log and $\text{Li}_2\left(\frac{y}{y\pm1}\right)$ are directly evaluated. Thus, up to a constant factor $-i$, the only one remains is $\int_0^1 \frac{2 h(y)\log (y)}{1-y^2} \, dy$, where\\
 $$h(y)=\text{Li}_2\left(\frac{(1+i) y}{y-1}\right)-\text{Li}_2\left(\frac{(1+i) y}{y+1}\right), h'(y)=\frac{\log \left(\frac{1-i y}{y+1}\right)}{y (y+1)}-\frac{\log \left(\frac{1+i y}{1-y}\right)}{(1-y) y}$$\\
Apply partial integration using:\\
$$K(y)=\int \frac{2 \log (y)}{1-y^2} \, dy=\text{Li}_2(-y)-\text{Li}_2(y)-\log (1-y) \log (y)+\log (y+1) \log (y)+\frac{\pi ^2}{4}$$\\
We have $\int_0^1 \frac{2 h(y)\log (y)}{1-y^2} \, dy$ boils down to $\int_0^1 K(y) h'(y) \, dy$, which is decomposed into combination of $QPLI4(1/2;\cdot;1/2/3)$ after separation, appropriate partial integration lifting $-\log(1-y)$, modulo LI/PLI/QLIs. Thus, we obtain 2 independent relations between the original $QPLI4(7;4/5;4)$ and these $QPLI4(1/2;\cdot;1/2/3)$. In fact, by combining all relations generated by subsection 6-1$\sim$6-10-4, all QPLIs of form $QPLI4(1/2;\cdot;1/2/3)$ are solved so $QPLI4(7;4/5;4)$ is evaluated.\\

Another relation between $QPLI4(7;3;4), QPLI4(2;5;5)$ and $QPLI4(2;4;4)$ can be obtained similarly if we consider applying this kernel on $\int_0^1 \frac{\text{Li}_2(-x^2) \log (1+ x)}{x^2+1} \, dx$ instead. Indeed, in this case the remaining $\int_0^1 K(y) h'(y) \, dy$ will be replaced by $\int_0^1 K(y) l'(y) \, dy$, where:\\
$$l(y)=\text{Li}_2\left(\frac{y}{-i+y}\right)-\text{Li}_2\left(\frac{2 y}{-i+y}\right)-\text{Li}_2\left(\frac{y}{i+y}\right)+\text{Li}_2\left(\frac{2 y}{i+y}\right)$$

$$l'(y)=-\frac{i \log \left(y^2+1\right)}{y^2+1}+\frac{2 i y \tan ^{-1}(y)}{y^2+1}-\frac{2 i \tan ^{-1}(y)}{y}$$\\
From which the mentioned relation is extracted.\\

\subsection{Solution}
One may verify that subsection 6-1$\sim$6-10 readily offers 171 independent relations, from which 160 QPLIs ($W\leq4$) out of the total 172 are solved, while the remaining 12 QPLIs are connected by 11 relations (generated by subsection 6-3, 6-4, 6-5, 6-10-3, 6-10-5). In other words, only one needs to be solved explicitly to complete the system. The 12 (up till now) unsolved are:\\
$$QPLI4(2/3/4/6;4;4), QPLI4(2/3/4/6;5;5)$$ 

$$QPLI4(8;3;5), QPLI4(8;4;3), QPLI4(7;3;4), QPLI4(7;5;3)$$\\
Apparently, $\int_0^1 \frac{\text{Li}_2(-y) \log (1-i y)}{1-i y} \, dy$ is evaluable, since by partial integration it boils down to integral $\int_0^1 \frac{\log ^2(1-i y) \log (y+1)}{y} \, dy$ then QLIs. On the other hand, taking its real part gives a relation between $QPLI4(2;4;4), QPLI4(2;5;5)$ (essentially this is subsection 6-3).\\

Now we evaluate $\int_0^1 \frac{\text{Li}_2(-y) \log (1+i y)}{1-i y} \, dy$ directly to obtain the last nontrivial relation. Apply double expansion using Prop. 3(2)(5), sum w.r.t $j$ in spirit of Prop. 5, after an index change we arrive at:\\
$$\int_0^1 \frac{\text{Li}_2(-x) \log (1+i x)}{1-i x} \, dx=\sum _{j=1}^{\infty } \sum _{k=1}^{\infty } \frac{(-1)^j i^k \widetilde{H_k}}{j^2 (j+k+1)}$$

$$=\sum _{k=1}^{\infty } i^k  \widetilde{H_k} \left(\frac{(-1)^k (\log (2)- \widetilde{H_{k+1}})}{(k+1)^2}-\frac{\pi ^2}{12 (k+1)}+\frac{\log (2)}{(k+1)^2}\right)$$

$$=\sum _{k=1}^{\infty } i^{k-1} \left(\widetilde{H_k}-\frac{(-1)^{k-1}}{k}\right) \left(\frac{(-1)^{k-1} (\log (2)-\widetilde{H_k})}{k^2}+\frac{\log (2)}{k^2}-\frac{\pi ^2}{12 k}\right)$$\\
Modulo $\zeta, \eta$ values and up to constant factors, we have to deal with 5 sums:\\
$$S_1=\sum _{n=1}^{\infty } \frac{i^n \widetilde{H_n}}{n}, S_{2/3}=\sum _{n=1}^{\infty } \frac{(\pm i)^n \widetilde{H_n}}{n^2}, S_4=\sum _{n=1}^{\infty } \frac{i^n \widetilde{H_n}}{n^3}, S_5=\sum _{n=1}^{\infty } \frac{i^n \left(\widetilde{H_n}\right)^2}{n^2}$$\\
Using Prop. 3(7), 4(7) and all closed-forms in subsection 5-3, plug in $x\to \pm i$, $S_1, S_2, S_3$ are solved. For $S_4$, by using $\widetilde{H_n}=\int_0^1 \frac{1-(-x)^n}{x+1} \, dx$ and partial integration we have\\
$$\sum _{n=1}^{\infty } \frac{i^n  \widetilde{H_n}}{n^3}=\int_0^1 \frac{\text{Li}_3(i)-\text{Li}_3(-i x)}{x+1} \, dx=\int_0^1 \frac{\log (x+1)}{x} \left(\frac{\text{Li}_2\left(-x^2\right)}{4}+i \text{Ti}_2(x)\right) \, dx+\frac{1}{16} i \pi ^3 \log (2)$$\\
Since $QPLI4(7/8;3;2)$ belong to 160 QPLIs already solved, $S_4$ is finished. To conquer $S_5$, recall Prop. 5(1), 3(9), exchange the order:\\
$$\sum _{n=1}^{\infty } \frac{i^n \left(\widetilde{H_n}\right)^2}{n^2}=\int_0^1 -\frac{\log (x) }{x}\sum _{n=1}^{\infty } \left( \widetilde{H_n}\right)^2 (-i x)^n \, dx$$ 

$$=\int_0^1 \left(-\frac{1}{x}+\frac{i}{1+i x}\right) \log (x) \left(\text{Li}_2(-i x)+2 \text{Li}_2\left(\frac{ i x+1}{2}\right)-2 \text{Li}_2\left(\frac{1}{2}\right)+A(x)\right) \, dx$$\\
Here $A(x)=2 \log (1-i x) \log (1+i x)-2 \log (2) \log (1+i x)$. Apparently $\int_0^1 \left(-\frac{1}{x}+\frac{i}{1+i x}\right) \log (x) A(x) dx$ trivially reduces to LI/QLIs. Moreover, by using partial integration lifting $\frac{1}{2}\log^2(x)$, polylogs in $\int_0^1 \frac{\left(\text{Li}_2(-i x)+2\text{Li}_2\left(\frac{ i x+1}{2}\right)-2 \text{Li}_2\left(\frac{1}{2}\right)\right) \log (x)}{x} \, dx$ part are eliminated hence also solved via QLIs. Furthermore, $\int_0^1 \frac{\text{Li}_2\left(\frac{1}{2}\right) \log (x)}{1+i x} \, dx, \int_0^1 \frac{\text{Li}_2(-i x) \log (x)}{1+i x} \, dx$ are trivial since $QPLI4(7/8;2;4/5)$ also belong to 160 known QPLIs. Therefore, up to constants the only one left is $\int_0^1 \frac{\text{Li}_2\left(\frac{i x+1}{2}\right) \log (x)}{1+i x} \, dx$. Like in subsection 6-10-1, with the aid of\\
$$\frac{\text{Li}_2\left(\frac{1}{2} (i x+1)\right)}{1+i x}=-\int_0^1 \frac{\log (y)}{2-(1+i x) y} \, dy \ (|x|<1)$$\\
We apply double integration\\
$$\int_0^1 \frac{\text{Li}_2\left(\frac{i x+1}{2}\right) \log (x)}{1+i x} \, dx=-\int _0^1\int _0^1\frac{\log (x) \log (y)}{2-(1+i x) y}dydx=-i\int_0^1 \frac{\text{Li}_2\left(\frac{i y}{2-y}\right) \log (y)}{y} \, dy$$\\
Lift up $\frac{1}{2}\log^2(y)$ again, let $y\to 1-y$ and simplify, we have it equals to \\
$$-\frac{i}{2}\int_0^1 \frac{\log ^2(1-y) \log \left(\frac{(1-i) (1+i y)}{y+1}\right)}{1-y} \, dy-\frac{i}{2}\int_0^1 \frac{\log ^2(1-y) \log \left(\frac{(1-i) (1+i y)}{y+1}\right)}{y+1} \, dy$$\\
 Obviously they reduce to LI/QLIs after separation and a possible partial integration lifting up $-\frac{1}{3}\log^3(1-y)$ for the former, so $S_5$ is computed.\\

Combining $S_1\sim S_5$, desired $\int_0^1 \frac{\text{Li}_2(-x) \log (1+i x)}{1-i x} \, dx$ is solved. Add it up to $\int_0^1 \frac{\text{Li}_2(-x) \log (1-i x)}{1-i x} \, dx$ above then take real part, $QPLI4(2;4;4)$ hence other 11 unknowns are calculated, finishing evaluation of all 172 QPLIs with weight $\leq4$. Thus, Theorem 2 is proved.\\

\section{QES}
We solve 83 QESs with $W\leq4$ (out of 93) and 28 non-alt QESs with $W=5$ using non-MZV approach. 10 tough alternating weight 4 QESs remain are deduced with help of MZV theory \cite{ref24}.\\

\subsection{Non-alternating case}
For each non-alternating QES, write $H_n=H_{2n}-\widetilde H_{2n}, \frac{1}{n^p}=\frac{2^p}{(2 n)^p}$, use\\
$$\sum _{n=1}^{\infty } f(2 n)=\frac{1}{2} \left(\sum _{n=1}^{\infty } f(n)-\sum _{n=1}^{\infty } (-1)^{n-1} f(n)\right)$$\\
and possible index change to decompose it into alternating ESs. For instance we have\\
$$QES(11,31,31;22)=\sum _{n=1}^{\infty } \frac{(H_{2n}-\widetilde H_{2n}) \left(H_{2 n}\right){}^2}{(2 n+1)^2}=\sum _{n=0}^{\infty } \frac{1-(-1)^{n-1}}{2}\frac{(H_{n}-\widetilde H_{n}) \left(H_{n}\right){}^2}{(n+1)^2}$$

$$=\sum _{n=0}^{\infty } \frac{1-(-1)^{n-1}}{2} \frac{\left(H_{n+1}-\widetilde H_{n+1}+\frac{(-1)^n}{n+1}-\frac{1}{n+1}\right) \left(H_{n+1}-\frac{1}{n+1}\right){}^2}{(n+1)^2}$$

$$=\sum _{n=1}^{\infty } \frac{1+(-1)^{n-1}}{2} \frac{\left(H_n-\widetilde H_n+\frac{(-1)^{n-1}}{n}-\frac{1}{n}\right) \left(H_n-\frac{1}{n}\right){}^2}{n^2}$$\\
Expand this expression. Since all ESs of weight 5 are recorded/computed in section 3, the QES is solved directly. Now do the same to all non-alt QESs.\\

\subsection{Generating functions}
For each formula in Prop. 3 or 4, let $z\to i$, take $\Re/\Im$ (like in subsection 5-15, 6-9-1, etc) then apply formulas of subsection 5-3 to simplify, we obtain QES relations after possible decomposition on $H_{2n+1}, \widetilde H_{2n+1}$. Now consider all formulas.\\

\subsection{QLI/QPLI expansion}
For each QLI/QPLI, carry out the procedure of subsection 4-7 to obtain QES relations. For instance, double expansion on $\text{Li}_2(x),\frac{\log \left(x^2+1\right)}{x^2+1}$ gives\\
$$QPLI4(1,4;4)=\sum _{n=1}^{\infty } (-1)^{n-1} H_n \left(-\frac{H_{2n}}{(2 n+1)^2}+\frac{\pi ^2}{6 (2 n+1)}-\frac{1}{(2 n+1)^3}\right)$$\\
Which directly leads to $QES(11,31;-22)$ modulo trivial terms. As another example, one may verify $QPLI4(7,4;4), QES(11,12;-12)$ are connected by double expansion too. Now it's natural to consider all QLI/QPLIs and all possible choices in the expansion procedure.\\

\subsection{Integral representation}
For each formula in Prop. 5 (and other integral representations), multiply both sides by appropriate harmonic factor, sum up w.r.t $n$, use Prop. 3 (or other generating functions) to simplify LHS, reduce them to known LI/PLI/QLI/QPLIs via various manipulations (sub, IBP, etc).\\

Other integral representations are those featuring single harmonic terms, for instance\\
$$H_n=\int_0^1 \frac{1-x^n}{1-x} \, dx, \widetilde {H_n^{(2)}}=\int_0^1 -\frac{\left(1-(-x)^n\right) \log (x)}{x+1} \, dx, H_{2 n}^{(3)}=\int_0^1 \frac{\left(1-x^{2 n}\right) \log ^2(x)}{2 (1-x)} \, dx$$\\
or those composite ones like (readers may complete them)\\
$$\int_0^1 \frac{x^{2 n}}{x+1} \, dx=H_n-H_{2 n}+\log (2), \int_0^1 \frac{x^n}{x+1} \, dx=(-1)^{n-1} \left(\widetilde{H_n}-\log (2)\right)$$

$$\int_0^1 x^{n-1} \log (x) \log (x+1) \, dx,\int_0^1 x^{n-1} \log (1-x) \log (x+1) \, dx,\int_0^1 x^{n-1} \log ^2(x+1) \, dx$$\\

Other generating functions are derived from Prop. 3, i.e. consider $f(\pm z^a), f(z^a)\pm f(-z^a)$ or $\Re/\Im f(iz^a)$ where $a\in\{1,2,\frac{1}{2}\}$, for instance\\
$$\sum _{n=1}^{\infty } \frac{(-1)^n x^n}{2 n+1}=\frac{\tan ^{-1}\left(\sqrt{x}\right)-\sqrt{x}}{\sqrt{x}}, \sum _{n=0}^{\infty } \frac{x^{2 n+1}}{(2 n+1)^2}=\frac{1}{2} (\text{Li}_2(x)-\text{Li}_2(-x))$$

$$\sum _{n=1}^{\infty } H_{2 n} x^{2 n}=-\frac{1}{2} \left(\frac{\log (1-x)}{1-x}+\frac{\log (x+1)}{x+1}\right),  \sum _{n=1}^{\infty } (-1)^n H_{2 n}^{(2)} x^{2 n}=\frac{\text{Li}_2\left(-x^2\right)}{4 \left(x^2+1\right)}-\frac{x \text{Ti}_2(x)}{x^2+1}$$\\

The general process is rather tricky and complicated so we only show one example. Consider $QES(21,31;22)$. Modulo trivial $QES(21;32)$ we only need to calculate $\sum _{n=0}^{\infty } \frac{H_{2 n+1} \widetilde{H_n}}{(2 n+1)^2}$. By manipulating Prop. 3(15) one have (all polylogs principal)\\
\begin{equation}
\begin{aligned}
F(x)=\sum _{n=0}^{\infty } \frac{(-1)^n H_{2 n+1} x^{2 n}}{(2 n+1)^2}=-\frac{i }{2 x}\left(-\text{Li}_3(-i x)+\text{Li}_3(i x)+\text{Li}_3(i x+1)-\text{Li}_3(1-i x)\ \ \ \ \ \ \ \ \ \right. \\ \left.-\text{Li}_2(i x+1) \log (1+i x)+\text{Li}_2(1-i x) \log (1-i x)-\frac{1}{2} \log (-i x) \log ^2(1+i x)+\frac{1}{2} \log (i x) \log ^2(1-i x)\right)
\end{aligned}
\notag
\end{equation}\\
Hence, recall $\widetilde{H_n}=\int_0^1 \frac{2 x \left(1-\left(-x^2\right)^n\right)}{x^2+1} \, dx$, it is direct that $\sum _{n=0}^{\infty } \frac{H_{2 n+1} \widetilde{H_n}}{(2 n+1)^2}=\int_0^1 \frac{2 x (-F(x)+F(i))}{x^2+1} \, dx$. Expand the integrand, apply partial integration to $\frac{\text{Li}_3(\cdots)}{x^2+1}$ terms then Prop. 7(2) for $\frac{\text{Li}_2(\cdots) \log(\cdots)}{x^2+1}$, $\frac{\text{Li}_2(\cdots) \tan^{-1}(\cdots)}{x^2+1}$ terms, simplify and plug in QLI/QPLI values, the original sum is evaluated.\\

\subsection{Contour integration}
Identical to subsection 3-4, see \cite{ref11} for examples and general process.\\

\subsection{Partial fractions}
Generalize subsection 3-3 to consideration of $\frac{(\pm 1)^{j-1} (\pm 1)^{k-1}}{f(j)^a g(k)^b}$ where $f, g\in \{2t, 2t+1\}$. For those relations involving summand $J_n^{(k)}=\sum _{j=1}^n \frac{(-1)^{j-1}}{(2 j+1)^k}$, try subsection 7-4, i.e. use generating functions/integral representations of $J_n^{(k)}$ to solve it by QLI/QPLI.\\

\subsection{Symmetric relations}
Generalize 3-2$\to$7-7 like 3-3$\to$7-6. Do the same as 7-6 for $J_n^{(k)}$ involved sums.\\

\subsection{Solution}
By combining methods listed above, all but 10 weight 4 alternating QESs are solved. As one may notice, all integrals and sums in sections above are given closed-forms by combining MZV-free methods, but it simply terminates here, suggesting us getting help from more powerful tools. Indeed, based on closed-forms of weight 4 level 4 MZVs \cite{ref24} and QES-MZV reduction, the remaining 10 sums are evaluated immediately. They will be dyed \textcolor{blue}{blue} in tables of section 9.\\

\section{Applications}

\subsection {Nonhomogeneous integrals}
\subsubsection{NLI}
Define nonhomogeneous log integrals (abbr. NLIs) as\\
$$NLI(a(0),a(1),a(2);k)=\int_0^1 x^k \prod _{m=0}^2 l(m;x)^{a(m)} dx$$\\
Whenever the order $k\in \mathbb{Z}-\{-1\}$. \\

\noindent \textbf{Proposition 11.} All convergent NLIs with weight $W=a(0)+a(1)+a(2)\leq 5$ are generated by $\cup_{w=0}^W A_w$ over $\mathbb{Q}$.\\

\noindent Proof. Integrate by parts repeatedly. Suppose $k$ is nonnegative, lift up $x^k$ to $\frac{x^{k+1}}{k+1}$ (or $\frac{x^{k+1}-1}{k+1}$ sometimes for convergence issues) and differentiate the log terms, after splitting and partial fraction decomposition on $\frac{x^{k+1}}{1\pm x}$, it boils down to several weight $W$ LIs and weight $W-1$ NLIs. Hence an induction on NLIs' weight completes the proof. For negative $k$ the technique is similar. $\hfill\square$\\

In fact, by partial fraction decomposition and method above, all convergent integrals of form $\int_0^1 \frac{ \prod _{m=0}^2 l(m;x)^{a(m)}}{(1-x)^p x^q (x+1)^r} dx$ are evaluable.\\

\subsubsection{NPLI, NQLI, NQPLI}
Similarly we define nonhomogeneous polylog integrals, nonhomogeneous quadratic log integrals and nonhomogeneous quadratic polylog integrals (abbr. NPLIs, NQLIs, NQPLIs) by replacing PLI/QLI/QPLI denominators by $x^k, k\not=-1$. They share similar structure with NLIs, that is, completely generated by $\cup_{w=0}^W A_w (W\leq5)$ or $\cup_{w=0}^W B_w(W\leq 4)$ for fixed weight $W$. The proof, using partial integration, is identical to NLI case. Generalizations with denominator $(1-x)^p x^q (x+1)^r (x^2+1)^s$ is also available for NQLI/NQPLIs. Also by manipulating these integrals some logsine variations can be deduced, for instance\\
$$\int_0^{\frac{\pi }{2}} x \cos (x) \log (2 \sin (x)) \log (2 \cos (x)) dx$$ $$=2 C-4 C \log (2)+4 \Im(\text{Li}_3(1+i))+\frac{\pi ^2}{8}-\frac{\pi ^3}{8}+\pi -6+\frac{3}{4} \pi  \log ^2(2)-2 \log ^2(2)-\frac{3}{2} \pi  \log (2)+6 \log (2)$$\\
One may define NES, NQES similarly by changing denominators of the summands. Their nonhomogeneous structures will be fully established by partial fractions and reindexing. Notwithstanding, they won't be discussed later.\\

\subsection{Nested sums and multiple log integrals}
\subsubsection{NS manipulation}
Some PLIs not evaluable directly in subsection 4-7 can be used to evaluate some nested harmonic sums (i.e. NSs). For instance, by double expansion and Cauchy summation\\
$$\int_0^1 \frac{\text{Li}_2(-x) \log ^2(1-x)}{x (x+1)} \, dx=\sum _{k=1}^{\infty } \sum _{j=1}^{\infty } \frac{(-1)^{j+k-1}\widetilde{H_j}  H_{j+k}}{k^2 (j+k)}=\sum _{m=1}^{\infty } \frac{(-1)^{m-1} H_m }{m}\sum _{k=1}^{m-1} \frac{\widetilde{H_{m-k}}}{k^2}$$\\
LHS is known. Since $\sum _{k=1}^{m-1} \frac{(-1)^{m-k}}{k^2 (m-k)}$ is elementary and that: \\
$$\sum _{n=1}^m (-1)^{n-1} \left(\frac{\widetilde{H_n}}{n^2}+\frac{\widetilde{H_n^{(2)}}}{n}\right)=\widetilde{H_n} \widetilde{H_n^{(2)}}+H_n^{(3)}$$\\
 we know \\
 $$\sum _{k=1}^{m-1} \frac{\widetilde{H_{m-k}}}{k^2}=\sum _{n=1}^m \frac{\widetilde{H_n}}{n^2}-\widetilde{H_m} \widetilde{H_m^{(2)}}-\widetilde{H_m^{(3)}}+H_m^{(3)}$$\\
 Using it to simplify the first identity and plugging in ES values, we obtain the closed-form of the following:\\
\begin{equation}
\begin{aligned}
\sum _{m=1}^{\infty } \frac{(-1)^{m-1} H_m}{m} \sum _{n=1}^m \frac{\widetilde{H_n}}{n^2}=-10 \text{Li}_5\left(\frac{1}{2}\right)+\frac{\pi ^2 \zeta (3)}{96}+\frac{277 \zeta (5)}{32}\\+\frac{9}{16} \zeta (3) \log ^2(2)+\frac{\log ^5(2)}{12}-\frac{1}{72} \pi ^2 \log ^3(2)-\frac{1}{18} \pi ^4 \log (2)\ \ \ \
\end{aligned}
\notag
\end{equation}\\

Another way of solving NSs is to exploit symmetry. We have outer and inner symmetric formula:\\
$$(O)\ \ \sum _{m=1}^{\infty } g(m) \sum _{n=1}^m f(n)+\sum _{m=1}^{\infty } f(m) \sum _{n=1}^m g(n)=\sum _{m=1}^{\infty } f(m) \sum _{m=1}^{\infty } g(m)+\sum _{m=1}^{\infty } f(m) g(m)$$

{\small $$(I)\ \ \sum _{m=1}^{\infty } f(m) \sum _{n=1}^m g(n) \sum _{k=1}^n h(k)+\sum _{m=1}^{\infty } f(m) \sum _{n=1}^m h(n) \sum _{k=1}^n g(k)=\sum _{m=1}^{\infty } f(m) \left(\sum _{n=1}^m g(n) \sum _{n=1}^m h(n)+\sum _{n=1}^m g(n) h(n)\right)$$}\\
Apply these and $n$-ple generalizations repeatedly with all $f, g, h, \cdots$ arbitrary harmonic terms, we obtain some relations between NSs. For example, all weight 5 RSs of form $$\sum _{m=1}^{\infty } f(m) \sum _{n=1}^m \frac{H_n}{n}, \sum _{m=1}^{\infty } f(m) \sum _{n=1}^m (-1)^{n-1} \frac{\widetilde{H_n}}{n}$$ with $f$ harmonic term can be reduced to ESs directly using inner $\sum _{n=1}^m \frac{H_n}{n}=\frac{1}{2} \left(\left(H_m\right){}^2+H_m^{(2)}\right)$ and its alternating analogue, while $\sum _{m=1}^{\infty } \frac{\widetilde{H_m} }{m^2}\sum _{n=1}^m (-1)^{n-1}\frac{ H_n}{n}$ is evaluated by outer symmetry, known ESs and the NS computed above. \\

\subsubsection{MLI representation}
Plenty of multiple log integrals (i.e. MLIs) are evaluable by using ES/NSs. For instance, using $\widetilde{H_n}=\int_0^1 \frac{1-(-x)^n}{x+1} \, dx$ thus $\sum _{n=1}^m \frac{\widetilde{H_n}}{n^2}=\int_0^1 \frac{H_m^{(2)}-\sum _{n=1}^m \frac{(-x)^n}{n^2}}{x+1} \, dx$, Proposition 5 and Fubini we have\\
$$\sum _{m=1}^{\infty } \frac{(-1)^{m-1} H_m}{m} \sum _{n=1}^m \frac{\widetilde{H_n}}{n^2}=\int _0^1\int _0^1\frac{\log (1-y) \left(\sum _{m=1}^{\infty } H_m^{(2)} (-y)^m- \sum _{m=1}^{\infty }\sum _{n=1}^m \frac{(-y)^m (-x)^n}{n^2}\right)}{(x+1) y}dxdy$$\\
$$=\int _0^1\int _0^1\frac{\log (1-y) (\text{Li}_2(-y)-\text{Li}_2(-x y))}{(x+1) y (y+1)}dxdy$$\\
In general, we use Proposition 5 (and its alternating analogue) then exchange order repeatedly to obtain the integral representation for certain sums.\\

\subsection{Geometric and binomial Euler sums}
\subsubsection{GES}
In section 3 we have made use of GES results i.e. geometric Euler sums of form:\\
$$GES(a(1),\cdots,a(M);p)=\sum _{n=1}^{\infty } \frac{1 }{n^p 2^n}\prod _{k=1}^M H_n^{(a(k))}$$\\
According to Xu \cite{ref19}, all GESs with weight $W=\sum _{k=1}^M a(k)+p\leq5$ and $p>0$ are generated by $A_W$ over $\mathbb{Q}$. In fact this trivially extends to the case $p=0$ by the reindex trick, for instance \\
$$\frac{1}{2} \sum _{n=1}^{\infty } \frac{\left(H_n\right){}^3 H_n^{(2)}}{2^n}=\sum _{n=1}^{\infty } \frac{\left(H_n-\frac{1}{n}\right){}^3 \left(H_n^{(2)}-\frac{1}{n^2}\right)}{2^n}$$\\
From which we may extract the value of LHS since all GESs on RHS with $p>0$ are known. Thus all 39 GESs with $W\leq 5$ and $p\geq0$ are computed (we exclude the trivial polylog case where no $a(k)$ exists).\\

\subsubsection{BES}
Another generalization is binomial Euler sum (abbr, BES):\\
$$BES(a(1),\cdots,a(M);p)=\sum_{n=1}^\infty\frac{\binom{2n}{n}}{n^p 4^n}\prod _{k=1}^M H_n^{(a(k))}$$ \\
The weight is naturally defined as $W=\sum _{k=1}^M a(k)+p$ and $p>0$ ensuring convergence. Different from GES, we take the case where no $a(k)$ exists into account.\\

\noindent \textbf{Proposition 12.} All 26 BESs with weight $W\leq 5$ are generated by $A_W$ over $\mathbb{Q}$.\\

\noindent Proof. All can be done by a single method. Denote the set of harmonic integral representations as $A$, which contains 3 classes of formulas:\\
$$\int_0^1 x^{n-1} \log ^k(x) \, dx\ (k\leq3), \int_0^1 \frac{\left(1-x^n\right) \log ^k(x)}{1-x} \, dx\ (k\leq3)$$ \\
featuring $\frac{1}{n^{k+1}}$ and $H_n^{(k+1)}$, and\\
$$\int_0^1 x^{n-1} \log ^k(x) \log ^{m-k}(1-x) \, dx\ (m\leq4, 0\leq k<m)$$\\
featuring composite harmonic terms. Part of them are recorded in Prop. 5 (see also subsection 7-4) and we tabulate the rest here.\\
$$\int_0^1 x^{n-1} \log ^2(x) \log (1-x) \, dx=-\frac{2 H_n}{n^3}-\frac{2 H_n^{(2)}}{n^2}-\frac{2 H_n^{(3)}}{n}+\frac{\pi ^2}{3 n^2}+\frac{2 \zeta (3)}{n}$$

$$\int_0^1 x^{n-1} \log (x) \log ^2(1-x) \, dx=-\frac{2 H_n H_n^{(2)}}{n}-\frac{\left(H_n\right){}^2}{n^2}+\frac{\pi ^2 H_n}{3 n}-\frac{H_n^{(2)}}{n^2}-\frac{2 H_n^{(3)}}{n}+\frac{2 \zeta (3)}{n}$$

$$\int_0^1 x^{n-1} \log ^4(1-x) \, dx=\frac{6 \left(H_n\right){}^2 H_n^{(2)}}{n}+\frac{8 H_n H_n^{(3)}}{n}+\frac{\left(H_n\right){}^4}{n}+\frac{3 \left(H_n^{(2)}\right){}^2}{n}+\frac{6 H_n^{(4)}}{n}$$

$$\int_0^1 x^{n-1} \log (x) \log ^3(1-x) \, dx=\frac{3 H_n H_n^{(2)}}{n^2}+\frac{3 \left(H_n\right){}^2 H_n^{(2)}}{n}+\frac{6 H_n H_n^{(3)}}{n}+\frac{\left(H_n\right){}^3}{n^2}$$ $$-\frac{6 \zeta (3) H_n}{n}-\frac{\pi ^2 \left(H_n\right){}^2}{2 n}+\frac{2 H_n^{(3)}}{n^2}+\frac{3 \left(H_n^{(2)}\right){}^2}{n}+\frac{6 H_n^{(4)}}{n}-\frac{\pi ^2 H_n^{(2)}}{2 n}-\frac{\pi ^4}{15 n}$$

$$\int_0^1 x^{n-1} \log ^2(x) \log ^2(1-x) \, dx=\frac{4 H_n H_n^{(2)}}{n^2}+\frac{4 H_n H_n^{(3)}}{n}+\frac{2 \left(H_n\right){}^2}{n^3}-\frac{2 \pi ^2 H_n}{3 n^2}$$

$$-\frac{4 \zeta (3) H_n}{n}+\frac{2 H_n^{(2)}}{n^3}+\frac{4 H_n^{(3)}}{n^2}+\frac{2 \left(H_n^{(2)}\right){}^2}{n}+\frac{6 H_n^{(4)}}{n}-\frac{2 \pi ^2 H_n^{(2)}}{3 n}-\frac{4 \zeta (3)}{n^2}-\frac{\pi ^4}{90 n}$$

$$\int_0^1 x^{n-1} \log ^3(x) \log (1-x) \, dx=\frac{6 H_n}{n^4}+\frac{6 H_n^{(2)}}{n^3}+\frac{6 H_n^{(3)}}{n^2}+\frac{6 H_n^{(4)}}{n}-\frac{\pi ^2}{n^3}-\frac{6 \zeta (3)}{n^2}-\frac{\pi ^4}{15 n}$$\\

These can be proved by using Beta derivatives like in Prop. 5. Now, denote \\
$$f_k(x)=\sum _{n=1}^{\infty } \frac{\binom{2 n}{n} x^n}{4^n n^k}, fh_j(x)=\sum _{n=1}^{\infty } \frac{\binom{2 n}{n} H_n x^n}{4^n n^k}$$\\
binomial generating functions $(k=0,1,2,3; j=0,1)$ and the set of these 6 functions as $B$, then we have closed-form expressions:\\
$$(1)(2)(3) \ f_0(x)=\frac{1-\sqrt{1-x}}{\sqrt{1-x}}, f_1(x)=2 \log (2)-2 \log \left(\sqrt{1-x}+1\right), fh_0(x)=\frac{2 \log \left(\frac{\sqrt{1-x}+1}{2 \sqrt{1-x}}\right)}{\sqrt{1-x}}$$

$$(4)\ fh_1(x)=2 \text{Li}_2\left(1-\sqrt{1-x}\right)+2 \text{Li}_2\left(-\sqrt{1-x}\right)+2 \text{Li}_2\left(\frac{1}{2} \left(\sqrt{1-x}+1\right)\right)$$
$$+\log ^2\left(\sqrt{1-x}+1\right)-\log (4) \log \left(1-\sqrt{1-x}\right)+\log (1-x) \log \left(1-\sqrt{1-x}\right)$$
$$+2 \log \left(\frac{1}{\sqrt{1-x}}+1\right) \log \left(1-\sqrt{1-x}\right)-2 \log \left(\frac{1}{\sqrt{1-x}}+1\right) \log \left(\sqrt{1-x}+1\right)-\frac{\pi ^2}{6}+\log ^2(2)$$

$$(5)\ f_2(x)=-2 \text{Li}_2\left(\frac{1}{2} \left(\sqrt{1-x}+1\right)\right)-\log ^2\left(\sqrt{1-x}+1\right)$$
$$+2 \log (2) \log (x)-2 \log \left(\frac{1}{2} \left(1-\sqrt{1-x}\right)\right) \log \left(\sqrt{1-x}+1\right)+\frac{\pi ^2}{3}-3 \log ^2(2)$$

{\scriptsize $$(6)\ f_3(x)=2 \text{Li}_3\left(\frac{1}{2} \left(1-\sqrt{1-x}\right)\right)-2 \text{Li}_3\left(\frac{1}{2} \left(\sqrt{1-x}+1\right)\right)-2 \text{Li}_2\left(\frac{1}{2} \left(1-\sqrt{1-x}\right)\right) \log (x)$$ $$+2 \text{Li}_2\left(\frac{1}{2} \left(1-\sqrt{1-x}\right)\right) \log \left(\sqrt{1-x}+1\right)-2 \text{Li}_2\left(\frac{1}{2} \left(\sqrt{1-x}+1\right)\right) \log (x)+2 \text{Li}_2\left(\frac{1}{2} \left(\sqrt{1-x}+1\right)\right) \log \left(\sqrt{1-x}+1\right)$$ $$-\frac{1}{3} \log ^3\left(\sqrt{1-x}+1\right)-3 \log ^2(2) \log (x)-2 \log ^2(2) \log \left(\sqrt{1-x}-1\right)+\log ^2(2) \log \left(\sqrt{1-x}+1\right)+\log (4) \log ^2(x)$$ $$-\log (4) \log ^2\left(\sqrt{1-x}-1\right)+\log \left(1-\sqrt{1-x}\right) \log ^2\left(\sqrt{1-x}+1\right)+\log (4) \log (2) \log (x)-\log (8) \log (2) \log \left(\sqrt{1-x}+1\right)$$ $$+\log (2) \log \left(1-\sqrt{1-x}\right) \log \left(\sqrt{1-x}+1\right)+\frac{1}{3} \pi ^2 \log (x)-\log (4) \log \left(1-\sqrt{1-x}\right) \log (x)-2 \log \left(1-\sqrt{1-x}\right) \log \left(\sqrt{1-x}+1\right) \log (x)$$ $$+\log (4) \log (x) \log \left(\sqrt{1-x}-1\right)+\log (4) \log \left(1-\sqrt{1-x}\right) \log \left(\sqrt{1-x}-1\right)-\log (4) \log \left(1-\sqrt{1-x}\right) \log \left(\sqrt{1-x}+1\right)$$ $$+\log (8) \log \left(1-\sqrt{1-x}\right) \log \left(\sqrt{1-x}+1\right)-\log (4) \log \left(\sqrt{1-x}-1\right) \log \left(\sqrt{1-x}+1\right)+2 \zeta (3)+\frac{4 \log ^3(2)}{3}+2 i \pi  \log ^2(2)-\frac{7}{3} \pi ^2 \log (2)$$}\\

In which we applied binomial expansion for $f_0$, the identity $fh_0(x)=\int_0^1 \frac{f_0(x)-f_0(t x)}{1-t} \, dt$, brute force and repeated use of $L(f)$ in Prop 3.\\

Fix a formula in $A$ and a (harmonic-)binomial term belongs to $\{\frac{\binom{2 n}{n}}{4^n n^k},\frac{\binom{2 n}{n} H_n}{4^n n^j}\}, (k=0,1,2,3; j=0,1)$. Muptiply both sides of the chosen formula with this term, sum w.r.t $n$, use the corresponding generating functions in $B$ to convert LHS into an integral and RHS a combination of BESs. Substitute $x\to 1-x^2$ in LHS, the integral can be decomposed into PLIs (in all cases here, with weight $\leq5$), since the integrand only contains logs and polylogs of form $\text{Li}_2(-x),\text{Li}_2(1-x),\text{Li}_2\left(\frac{1\pm x}{2}\right)$ or $\text{Li}_3\left(\frac{1\pm x}{2}\right)$ which are main focuses of section 4 (with a possible reflection on the second). Since LHS is known now, we obtain a relation of BESs on RHS. Consideration of all possible choices in $A, B$ completes the evaluation. $\hfill\square$\\

Here are two examples of the procedure above. Composite relations generated by, say the kernel $\int_0^1 x^{n-1} \log ^2(x) \log ^2(1-x) \, dx$, are cumbersome thus omitted here.\\
$$BES(4;1)=-\frac{1}{3} \int_0^1 \frac{\left(f_1(1)-f_1\left(1-x^2\right)\right) \log ^3\left(1-x^2\right)}{x} \, dx,$$
$$f_1(1)-f_1\left(1-x^2\right)=2 \log (x+1)$$

$$BES(13;1)=\int_0^1 \frac{\left(fh_1(1)-fh_1\left(1-x^2\right)\right) \log ^2\left(1-x^2\right)}{x} \, dx,$$
$$fh_1(1)-fh_1\left(1-x^2\right)=-2 \text{Li}_2\left(\frac{x+1}{2}\right)-2 \text{Li}_2(1-x)-2 \text{Li}_2(-x)-\log ^2(x+1)+\log (4) \log (1-x)$$$$-2 \log \left(\frac{1}{x}+1\right) \log (1-x)-2 \log (1-x) \log (x)+2 \log \left(\frac{1}{x}+1\right) \log (x+1)+\frac{\pi ^2}{2}-\log ^2(2)$$\\
Evidently they are evaluable via weight 5 LI/PLIs.\\

\subsubsection{IBES}
Now we deal with inverse binomial Euler sum (abbr, IBES):\\
$$IBES(a(1),\cdots,a(M);p)=\sum_{n=1}^\infty\frac{ 4^n}{n^p\binom{2n}{n}}\prod _{k=1}^M H_n^{(a(k))}$$ \\
The weight is defined as $W=\sum _{k=1}^M a(k)+p$ as usual and $p>1$ ensuring convergence. We also take the case where no $a(k)$ exists into account.\\

\noindent \textbf{Proposition 13.} All 14 IBESs with weight $W\leq 5$ are generated by $A_W$ over $\mathbb{Q}$.\\

\noindent Proof. Instead of PLI usage in Prop. 12, LSI/PLSI (subsection 5-1) will play an important role here. First of all,  recall generating functions of inverse sine power from \cite{ref8}:\\
$$(T1) (T2)\ d(x)=\sum _{n=1}^{\infty } \frac{4^n x^n}{n \binom{2 n}{n}}=\frac{2 \sqrt{x} \sin ^{-1}\left(\sqrt{x}\right)}{\sqrt{1-x}},\ f(x)=\sum _{n=1}^{\infty } \frac{4^n x^n}{n^2 \binom{2 n}{n}}=2 \sin ^{-1}\left(\sqrt{x}\right)^2$$

$$(T3)\ g(x)=\sum _{n=1}^{\infty } \frac{4^n H_n^{(2)} x^n}{n^2 \binom{2 n}{n}}-\sum _{n=1}^{\infty } \frac{4^n x^n}{n^4 \binom{2 n}{n}}=\frac{2}{3} \sin ^{-1}\left(\sqrt{x}\right)^4$$\\

$IBES(0;2)=f(1)$ is trivial. Multiply both sides of (T2) with $\frac{\log ^k(x)}{x}$ and integrate on (0,1). Let $x\to \sin^2(u)$, plug in LSI values, we get closed-forms of $IBES(0;p)$ for $3\leq p\leq 5$.\\

By using $H_n=\int_0^1 \frac{1-x^n}{1-x} \, dx$ for $IBES(1;2)$, Prop. 5(2) for $IBES(1;3)$ it is clear that\\
$$IBES(1;2)=\int_0^1 \frac{f(1)-f(x)}{1-x} \, dx, IBES(1;3)=\int_0^1 -\frac{f(x) \log (1-x)}{x} \, dx$$\\
Substitute $x\to \sin^2(u)$ in both identities, plug in LSI values, $IBES(1;2/3)$ are readily solved. Similar to $IBES(1;2)$ we have $IBES(3;2)=\int_0^1 \frac{(f(1)-f(x)) \log ^2(x)}{2 (1-x)} \, dx$, which is computed via LSI under the same substitution.\\

Moreover, by manipulating (T3) we have \\
$$IBES(2;2)=g(1)+IBES(0;4),\ IBES(2;3)=\int_0^1 \frac{2 g(x)}{3 x} \, dx+IBES(0;5)$$\\
$IBES(2;2)$ is solved directly. For the latter, let $x\to \sin^2(u)$ then plug in LSI value.\\

Furthermore, similar to solution of $IBES(1;3)$, multiply Prop. 5(3)(4)(5) by appropriate factors, sum w.r.t $n$ gives\\
{\small $$\int_0^1 \frac{d(x) \log ^2(1-x)}{x} \, dx=IBES(11;2)+IBES(2;2),\ \int_0^1 \frac{f(x) \log ^2(1-x)}{x} \, dx=IBES(11;3)+IBES(2;3)$$}
$$\int_0^1 \frac{f(x) \log (x) \log (1-x)}{x} \, dx=IBES(1;4)+IBES(2;3)-\frac{\pi^2}{6} IBES(0;3)$$

$$-\int_0^1 \frac{d(x) \log ^3(1-x)}{x} \, dx=3IBES(12;2)+IBES(111;2)+2IBES(3;2)$$\\
While another use of (T3) gives\\
$$\int_0^1 \frac{g(1)-g(x)}{1-x} \, dx=IBES(12;2)-IBES(1;4)$$\\
All 5 integrals on LHS boils down to known LSIs using $x\to \sin^2(u)$. Solving this system yields closed-forms of last 5 IBESs, i.e. $IBES(11;2/3), IBES(1;4), IBES(12;2), IBES(111;2)$, completing the proof. $\hfill\square$\\

\subsubsection{QBES, IQBES ($_pF_q$ part 1)}
We have quadratic binomial Euler sum and its inverse (abbr. QBESs, IQBESs):\\
$$QBES(a(1),\cdots,a(M);p)=\sum_{n=1}^\infty\left(\frac{\binom{2n}{n}}{4^n}\right)^2 \frac{1}{n^p}\prod _{k=1}^M H_n^{(a(k))}$$

$$IQBES(a(1),\cdots,a(M);p)=\sum_{n=1}^\infty\left(\frac{4^n}{\binom{2n}{n}}\right)^2 \frac{1}{n^p}\prod _{k=1}^M H_n^{(a(k))}$$  \\
Here $p>0$ for QBES and $p>2$ for IQBES, ensuring convergence. Defining weight $W=\sum _{k=1}^M a(k)+p$ as usual, there are 7 QBESs and 3 IQBESs for $W\leq4$, with 3 and 2 harmonic-free sums respectively.\\

Part 1: First of all we compute 5 harmonic-free sums.\\

\noindent \textbf{Proposition 14.} The following formulas hold:\\
$$(1)\ IQBES(0;3)=4 \, _4F_3\left(1,1,1,1;\frac{3}{2},\frac{3}{2},2;1\right)=8 \pi  C-14 \zeta (3)$$

$$(2)\ IQBES(0;4)=4 \, _5F_4\left(1,1,1,1,1;\frac{3}{2},\frac{3}{2},2,2;1\right)$$ $$=64 \pi  \Im(\text{Li}_3(1+i))+64 \text{Li}_4\left(\frac{1}{2}\right)-\frac{233 \pi ^4}{90}+\frac{8 \log ^4(2)}{3}-\frac{20}{3} \pi ^2 \log ^2(2)$$

$$(3)\ \pi QBES(0;1)=\frac{1}{4} \pi  \, _4F_3\left(1,1,\frac{3}{2},\frac{3}{2};2,2,2;1\right)=4 \pi  \log (2)-8 C$$

$$(4)\ \pi QBES(0;2)=\frac{1}{4} \pi  \, _5F_4\left(1,1,1,\frac{3}{2},\frac{3}{2};2,2,2,2;1\right)=64 \Im(\text{Li}_3(1+i))-2 \pi ^3-12 \pi  \log ^2(2)$$

$$(5)\ \pi QBES(0;3)=\frac{1}{4} \pi  \, _6F_5\left(1,1,1,1,\frac{3}{2},\frac{3}{2};2,2,2,2,2;1\right)$$ $$=-512 \Im(\text{Li}_4(1+i))+4 \pi  \zeta (3)+16 \pi  \log ^3(2)+8 \pi ^3 \log (2)+\frac{\psi ^{(3)}\left(\frac{1}{4}\right)}{4}-\frac{\psi ^{(3)}\left(\frac{3}{4}\right)}{4}$$\\

\noindent Proof. (1): By Clausen formula or manipulating (T1) above we have$\, _3F_2\left(1,1,1;\frac{3}{2},2;x\right)=\frac{\sin ^{-1}\left(\sqrt{x}\right)^2}{x}$. Multiply both sides by $\frac{1}{2 \sqrt{1-x}}$, integrate on (0,1). Now expand $_3F_2$ for LHS (see \cite{ref13} for general Euler integral of $_pF_q$) and let $x\to \sin^2(t)$ for RHS.\\

(2): Use the same technique in (1), integrate by parts (lifting up $\frac{1}{2 x \sqrt{1-x}}$), one have\\
 $$\, _5F_4\left(1,1,1,1,1;\frac{3}{2},\frac{3}{2},2,2;1\right)=\int_0^1 \frac{x \, _4F_3\left(1,1,1,1;\frac{3}{2},2,2;x\right)}{2 x \sqrt{1-x}} \, dx$$
 
$$=\int_0^1 -\frac{\left(\log \left(1-\sqrt{1-x}\right)-\log \left(\sqrt{1-x}+1\right)\right) \sin ^{-1}\left(\sqrt{x}\right)^2}{2 x} \, dx$$\\
Now substitute $x\to \sin^2(t)$, $t\to 2\tan^{-1}(u)$, apply partial fractions and plug in QLI values.\\

(3): By repeated integration we obtain  $\, _3F_2\left(1,1,\frac{3}{2};2,2;x\right)=-\frac{4 \log \left(\frac{1}{2} \left(\sqrt{1-x}+1\right)\right)}{x}$. Multiply both sides by $\frac{2 \sqrt{x}}{\pi  \sqrt{1-x}}$, integrate on (0,1). Now expand $_3F_2$ for LHS and let $x\to \sin^2(t)$ for RHS.\\

(4)(5): Using $L^*(f)=\frac{\int_0^x f(t) \, dt}{x}$ on (3) it's clear that\\
$$\, _4F_3\left(1,1,1,\frac{3}{2};2,2,2;x\right)=-\frac{4 \text{Li}_2\left(\frac{1}{2} \left(\sqrt{1-x}+1\right)\right)}{x}+\frac{2 \pi ^2}{3 x}$$
$$-\frac{4 \log \left(1-\sqrt{1-x}\right) \log \left(\frac{1}{2} \left(\sqrt{1-x}+1\right)\right)}{x}-\frac{2 \log \left(\frac{1}{8} \left(\sqrt{1-x}+1\right)\right) \log \left(\frac{1}{2} \left(\sqrt{1-x}+1\right)\right)}{x}$$

{\small $$\, _5F_4\left(1,1,1,1,\frac{3}{2};2,2,2,2;x\right)=\frac{4 \text{Li}_3\left(\frac{1}{2} \left(1-\sqrt{1-x}\right)\right)}{x}-\frac{4 \text{Li}_3\left(\frac{1}{2} \left(\sqrt{1-x}+1\right)\right)}{x}$$ $$+\frac{4 \text{Li}_2\left(\frac{1}{2} \left(1-\sqrt{1-x}\right)\right) \log \left(\sqrt{1-x}+1\right)}{x}+\frac{4 \text{Li}_2\left(\frac{1}{2} \left(\sqrt{1-x}+1\right)\right) \log \left(\sqrt{1-x}+1\right)}{x}$$ $$-\frac{4 \text{Li}_2\left(\frac{1}{2} \left(1-\sqrt{1-x}\right)\right) \log (-x)}{x}-\frac{4 \text{Li}_2\left(\frac{1}{2} \left(\sqrt{1-x}+1\right)\right) \log (-x)}{x}+\frac{4 \zeta (3)}{x}+\frac{2 i \pi ^3}{3 x}+\frac{2 \log ^3(2)}{3 x}$$ $$-\frac{2 \log ^3\left(\sqrt{1-x}+1\right)}{3 x}+\frac{\log (4) \log ^2(2)}{x}-\frac{2 \log ^2(2) \log \left(1-\sqrt{1-x}\right)}{x}-\frac{4 \log ^2(2) \log \left(\sqrt{1-x}-1\right)}{x}$$ $$-\frac{2 \log ^2(2) \log \left(\sqrt{1-x}+1\right)}{x}+\frac{2 \log \left(1-\sqrt{1-x}\right) \log ^2\left(\sqrt{1-x}+1\right)}{x}-\frac{\log (16) \log ^2\left(\sqrt{1-x}-1\right)}{x}$$ $$-\frac{\log ^2(4) \log \left(1-\sqrt{1-x}\right)}{x}-\frac{\log ^2(4) \log \left(\sqrt{1-x}+1\right)}{x}+\frac{\log (16) \log (2) \log \left(1-\sqrt{1-x}\right)}{x}-\frac{2 \pi ^2 \log (2)}{3 x}$$ $$+\frac{2 \pi ^2 \log (x)}{3 x}+\frac{\log (16) \log (-x) \log \left(\sqrt{1-x}-1\right)}{x}+\frac{\log (16) \log \left(1-\sqrt{1-x}\right) \log \left(\sqrt{1-x}-1\right)}{x}$$ $$+\frac{\log (16) \log (-x) \log \left(\sqrt{1-x}+1\right)}{x}+\frac{\log (16) \log \left(1-\sqrt{1-x}\right) \log \left(\sqrt{1-x}+1\right)}{x}$$ $$-\frac{4 \log \left(1-\sqrt{1-x}\right) \log \left(\sqrt{1-x}+1\right) \log (-x)}{x}-\frac{\log (16) \log \left(\sqrt{1-x}-1\right) \log \left(\sqrt{1-x}+1\right)}{x}$$}\\

We only establish (5), (4) is relevantly trivial. Multiply both sides of the $_5F_4$ formula  by $\frac{2 \sqrt{x}}{\pi  \sqrt{1-x}}$, integrate on (0,1). Expand $_5F_4$ for LHS, let $x\to \sin^2(t)$, $t\to 2\tan^{-1}(u)$ for RHS. After simplifications and modulo QLIs, there only left several nontrivial ones of form \\
$$\int_0^1 \frac{\text{Li}_2\left(p(x)\right) \log (q(x))}{x^2+1} \, dx, \int_0^1 \frac{\text{Li}_3(p(x))}{x^2+1} \, dx\to \int_0^1 \frac{\tan^{-1}(x)\text{Li}_2(p(x))}{r(x)} \, dx$$\\
Where $p\in\{\frac{1}{x^2+1}, \frac{x^2}{x^2+1}\}, q\in\{x,x^2+1\}, r\in \{x, \frac{x^2+1}{x}\}$, and the `$\to$' denotes IBP for 2 $\text{Li}_3$ terms and partial fractions. Now, use Prop. 7(2), 7(4) and 7(2) again, take care of complex arguments, modulo trivial QLI terms we reduce $\text{Li}_2(p(x))$ to ordinary $\text{Li}_2\left(-x^2\right)$:\\
$$\text{Li}_2\left(\frac{x^2}{x^2+1}\right)\to \text{Li}_2\left(\frac{1}{x^2+1}\right)\to \text{Li}_2\left(x^2+1\right)\to \text{Li}_2\left(-x^2\right)$$\\
Finally plug in QPLI to complete the proof. $\hfill\square$\\\\

Part 2: Now there left only one IQBES. We will establish that\\
$$ IQBES(1;3)=32 C^2-32 \pi  C \log (2)+128 \pi  \Im(\text{Li}_3(1+i))+64 \text{Li}_4\left(\frac{1}{2}\right)-\frac{413 \pi ^4}{90}+\frac{8 \log ^4(2)}{3}-\frac{32}{3} \pi ^2 \log ^2(2)$$\\
Indeed, by using Prop. 5(2), $\sum _{n=1}^{\infty }\left(\frac{4^n}{\binom{2 n}{n}}\right)^2  \frac{x^n}{n^2}=4 x \, _3F_2\left(1,1,1;\frac{3}{2},\frac{3}{2};x\right)$ then a shifted representation $\sum _{n=0}^{\infty } \left(\frac{4^n}{\binom{2 n}{n}}\right)^2\frac{x^n}{(2 n+1)^2}=\, _3F_2\left(1,1,1;\frac{3}{2},\frac{3}{2};x\right)$ and Prop. 5(2) again, one have\\
$$IQBES(1;3)=-4 \int_0^1 \log (1-x) \, _3F_2\left(1,1,1;\frac{3}{2},\frac{3}{2};x\right) \, dx,$$
 {\small $$ \int_0^1 \log (1-x) \, _3F_2\left(1,1,1;\frac{3}{2},\frac{3}{2};x\right) \, dx=\sum _{n=0}^{\infty }\int_{0}^{1}\left(\frac{4^n}{\binom{2n}{n}}\right)^2\frac{x^n \log(1-x)}{(2n+1)^2}\,dx=-\sum _{n=0}^{\infty }\left(\frac{4^n}{(2n+1)\binom{2n}{n}}\right)^2\frac{H_{n+1}}{n+1}$$}\\
Due to another shifted sum $\sum _{n=0}^{\infty } \frac{(4 x)^n}{(2 n+1) \binom{2 n}{n}}=\frac{\sin ^{-1}\left(\sqrt{x}\right)}{\sqrt{x (1-x)}}$, Prop. 5(2) again and $xz\to u$, we have\\
$$\sum_{n\geq 0}\frac{4^n H_{n+1}}{(n+1)(2n+1)\binom{2n}{n}}z^n=-\int_{0}^{z}\frac{\arcsin\sqrt{u}}{\sqrt{u(1-u)}}\log(1-\frac{u}z)\frac{du}z$$\\
Therefore, by using $\frac{1}{(2n+1)\binom{2n}{n}}= B(n+1,n+1)$, the generating function above and Fubini:\\
$$-\sum _{n=0}^{\infty }\left(\frac{4^n}{(2n+1)\binom{2n}{n}}\right)^2\frac{H_{n+1}}{n+1}=-\int_{0}^{1} \sum_{n\geq 0}\frac{4^n H_{n+1}}{(n+1)(2n+1)\binom{2n}{n}}(4x(1-x))^n\,dx$$

$$=\int_0^1 dx \int_{0}^{4x(1-x)}\frac{\arcsin\sqrt{u}}{\sqrt{u(1-u)}}\log(1-\frac{u}{4x(1-x)})\frac{du }{4x(1-x)}$$

$$=\int_0^1\frac{\arcsin\sqrt{u}}{\sqrt{u(1-u)}} du \left(\int_{\frac{1}{2} \left(-\sqrt{1-u}+1\right)}^{\frac{1}{2} \left(\sqrt{1-u}+1\right)}\log(1-\frac{u}{4x(1-x)})\frac{dx}{4x(1-x)}\right)$$\\
According to subsection 2-1, by brute force we have 8 times of the inner $x$ integral equals to\\
$$r(u)=-8 \text{Li}_2\left(\frac{1}{2} \left(1-\sqrt{1-u}\right)\right)+2 \text{Li}_2\left(\frac{2 \sqrt{1-u}}{\sqrt{1-u}-1}\right)+4 \text{Li}_2\left(\frac{1-\sqrt{1-u}}{\sqrt{1-u}+1}\right)-2 \text{Li}_2\left(\frac{2 \sqrt{1-u}}{\sqrt{1-u}+1}\right)$$ $$+\log ^2\left(1-\sqrt{1-u}\right)+\log ^2\left(\sqrt{1-u}+1\right)-2 \log \left(1-\sqrt{1-u}\right) \log (1-u)+\log (16) \log (u)$$ $$-\log (4) \left(-\log \left(\left(\sqrt{1-u}+1\right) u\right)+\log \left(1-\sqrt{1-u}\right)+\log (8)\right)-\log (16) \log \left(\sqrt{1-u}+1\right)$$ $$+2 \log (1-u) \log \left(\sqrt{1-u}+1\right)-6 \log \left(1-\sqrt{1-u}\right) \log \left(\sqrt{1-u}+1\right)-2 \log ^2(2)+\log ^2(4)$$\\
To calculate the resulting $\int_0^1\frac{ r(u)\arcsin\sqrt{u}}{8\sqrt{u(1-u)}} du$, apply $u\to \sin^2(t), t\to 2\tan^{-1}(v)$. After simplifications and modulo QLIs, there only left several nontrivial ones of form
$\int_0^1 \frac{\text{Li}_2\left(p(x)\right) \tan^{-1}(x)}{x^2+1} \, dx$ where $p\in \{x^2,1-x^2,1-\frac{1}{x^2},\frac{x^2}{x^2+1}\}$. By using Prop. 7 (1)(2)(4) flexibly and take care of branches, modulo trivial QLI terms we reduce $\text{Li}_2$ to ordinary arguments:\\
$$\text{Li}_2\left(1-x^2\right)\to \text{Li}_2\left(x^2\right)\to \text{Li}_2\left(\pm x\right),\ \text{Li}_2\left(1-\frac{1}{x^2}\right)\to\text{Li}_2\left(\frac{1}{x^2}\right)\to \text{Li}_2\left(x^2\right)\to \text{Li}_2\left(\pm x\right)$$

$$\text{Li}_2\left(\frac{x^2}{x^2+1}\right)\to \text{Li}_2\left(\frac{1}{x^2+1}\right)\to \text{Li}_2\left(x^2+1\right)\to \text{Li}_2\left(-x^2\right)$$\\
Finally plug in QPLI to complete evaluation of $IQBES(1;3)$.\\\\

Part 3: Now we compute the rest 4 QBESs. Firstly, in \cite{ref23} it is proved that\\
 $$\pi QBES(1;1)=32 C \log (2)-64 \Im(\text{Li}_3(1+i))+\frac{11 \pi ^3}{6}+4 \pi  \log ^2(2)$$\\
  So we only need to focus on evaluation of $\pi QBES(11;1),\pi QBES(2;1),\pi QBES(1;2)$.\\

\noindent \textbf{Proposition 15.} The following formulas hold:\\
$$(1)\ \ \pi QBES(1;2)=1024 \Im(\text{Li}_4(1+i))-256 \log (2) \Im(\text{Li}_3(1+i))$$ $$+\pi  \zeta (3)+\frac{16}{3} \pi  \log ^3(2)-8 \pi ^3 \log (2)+\frac{2}{3} \pi ^3 \log (2)-\frac{25 \psi ^{(3)}\left(\frac{1}{4}\right)}{48}+\frac{25 \psi ^{(3)}\left(\frac{3}{4}\right)}{48}$$

$$(2)\ \ \pi QBES(2;1)=\frac{20 \pi ^2 C}{3}-128 C \log ^2(2)-1024 \Im(\text{Li}_4(1+i))$$ $$+512 \log (2) \Im(\text{Li}_3(1+i))-14 \pi  \zeta (3)-\frac{64}{3}  \pi  \log ^3(2)+\frac{25 \psi ^{(3)}\left(\frac{1}{4}\right)}{48}-\frac{25 \psi ^{(3)}\left(\frac{3}{4}\right)}{48}$$

$$(3)\ \ \pi QBES(11;1)=-\frac{4 \pi ^2 C}{3}+12 \pi  \zeta (3)-\frac{\psi ^{(3)}\left(\frac{1}{4}\right)}{48}+\frac{\psi ^{(3)}\left(\frac{3}{4}\right)}{48}$$\\

\noindent Proof. In spirit of \cite{ref23}, if $f(x)=\sum_{n=0}^\infty c_n P_n(2x-1)$, then by Bonnet recursion formula \cite{ref13} and solving diffrence equation, one may prove\\
$$\frac{f(x)}x=\sum_{n=0}^\infty(-1)^n (2n+1)\left(\int_0^1 \frac{f(x)}x dx+2\sum_{m=1}^{n}\frac{1}{m}\sum_{k=m}^\infty (-1)^k c_k\right)P_n(2x-1)$$\\
Therefore, recall FL expansions of $\text{Li}_2(x), \log(x)\log(1-x), \log^2(1-x)$ in literature above, after simplification using subsection 8-2 we have\\
 $$\frac{\text{Li}_2(x)}x=\sum_{n=0}^\infty A_n P_n(2x-1),\frac{\log(x)\log(1-x)}x=\sum_{n=0}^\infty B_n P_n(2x-1),\frac{\log^2(1-x)}x=\sum_{n=0}^\infty C_n P_n(2x-1)$$
 
$$A_n=(-1)^n (2n+1)\left(\zeta (3)-\frac{\pi^2}{3} H_n +2\widetilde{H_n^{(3)}}+4H_n\widetilde{H_n^{(2)}}-4 \sum _{m=1}^n\frac{(-1)^{m-1} }{m^2}H_m\right)$$

$$B_n=(-1)^n (2n+1)\left(\zeta (3)+\frac{\pi^2}{3} H_n -2\widetilde{H_n^{(3)}}-4H_n\widetilde{H_n^{(2)}}-2H_n^{(3)}+4 \sum _{m=1}^n\frac{(-1)^{m-1} }{m^2}H_m\right)$$

$$C_n=(-1)^n (2 n+1) \left( 2 \zeta (3)+4\widetilde{H_n^{(3)}}-8 \sum _{m=1}^n\frac{(-1)^{m-1} }{m^2}H_m\right)$$\\
Now, by using Prop. 5(2)(3)(4) and that $\pi  \sum _{n=1}^{\infty } \left(\frac{\binom{2 n}{n}}{4^n}\right)^2 x^n=2 K(x)-\pi$, the following three equations are established:\\
{\small $$\pi QBES(1;2)=-\int_0^1 \frac{\pi  \log (1-x) }{x} \sum _{n=1}^{\infty } \left(\frac{\binom{2 n}{n}}{4^n}\right)^2\frac{ x^n}{n}\, dx=\pi QBES(0;1) \text{Li}_2(1)-\int_0^1 \frac{(2 K(x)-\pi ) \text{Li}_2(x)}{x} \, dx$$}

$$\int_0^1 \frac{(2 K(x)-\pi ) \log (1-x) \log (x)}{x} \, dx=\pi QBES(1;2)+\pi QBES(2;1)-\frac{\pi^3}{6} QBES(0;1)$$

$$\int_0^1 \frac{(2 K(x)-\pi ) \log ^2(1-x)}{x} \, dx=\pi QBES(11;1)+\pi QBES(2;1)$$\\
From which we observe that it suffices to compute\\
 $$\int_0^1 \frac{K(x) \text{Li}_2(x)}{x} \, dx,\int_0^1 \frac{K(x) \log (1-x) \log (x)}{x} \, dx,\int_0^1 \frac{K(x) \log ^2(1-x)}{x} \, dx$$\\
  (modulo known QBESs and trivial integrals) to complete evaluation of 3 QBESs. Using the classical $K(x)=\sum _{n=0}^{\infty } \frac{2 P_n(2 x-1)}{2 n+1}$, 3 expansions above and orthogonal identity $\int_0^1 P_p(2 x-1) P_q(2 x-1) \, dx=\frac{\delta _{p,q}}{2 p+1}$, they boils down to\\
   $$\sum _{n=0}^{\infty } \frac{2 A_n}{(2 n+1)^2},\ \sum _{n=0}^{\infty } \frac{2 B_n}{(2 n+1)^2},\ \sum _{n=0}^{\infty } \frac{2 C_n}{(2 n+1)^2}$$\\
Modulo weight 4 QESs and up to constant factors, the only notrivial sum in all three series above is a quadratic NS (i.e. QNS):\\
$$\sum_{n=1}^\infty \frac{(-1)^{n-1}}{2n+1}\sum _{m=1}^n\frac{(-1)^{m-1} }{m^2}H_m\to \sum_{n=1}^\infty \frac{(-1)^{n-1} }{n^2}H_n\sum _{m=1}^n \frac{(-1)^{m-1}}{2m+1}\to  \sum_{n=1}^\infty \frac{(-1)^{n-1} }{n^2}H_n\sum _{m=1}^n \frac{(-1)^{m-1}}{2m-1}$$\\
Where we have used outer symmetry (subsection 8-2) and reindexing for two `$\to$' repsectively, transforming the QNS into a extended QES modulo trivial sums, i.e.\\
$$\sum _{n=1}^{\infty } \frac{(-1)^{n-1}}{2 n+1}=\frac{4-\pi }{4},\ \sum _{n=1}^{\infty } \frac{(-1)^{n-1} H_n}{n^2}=\frac{5 \zeta (3)}{8},\ \sum _{n=1}^{\infty } \frac{H_n}{n^2 (2 n+1)}=-2 \zeta (3)-4 \log ^2(2)$$\\
 Since $\sum _{m=1}^n \frac{(-1)^{m-1}}{2m-1}=\int_0^1 \frac{1-\left(-x^2\right)^n}{x^2+1} \, dx$, plug in Prop. 4(15) and simplify by Prop. 7, the extended QES is equal to\\
{\small $$\int_0^1 \frac{1}{8 \left(x^2+1\right)} \left(8 \text{Li}_3\left(x^2\right)-8 \text{Li}_3\left(1-x^2\right)+8 \text{Li}_2\left(1-x^2\right) \log \left(1-x^2\right)+4 \log \left(x^2\right) \log ^2\left(1-x^2\right)+13 \zeta (3)\right)\, dx$$}\\
Lift up $\tan ^{-1}(x)$ for $\int_0^1 \frac{\text{Li}_3\left(1-x^2\right)}{x^2+1} \, dx$ to reduce $\text{Li}_3$ to $\text{Li}_2$, use Prop. 7(2) to have $\text{Li}_2(1-x^2)\to \text{Li}_2(x^2)$, IBP again, plug in QPLI values, 3 QBESs are solved. $\hfill\square$\\\\

Thus we have finished evaluation of 10 desired sums. Note that as solving $\sum _{n=0}^{\infty } \frac{2 A_n}{(2 n+1)^2}, \cdots$, 3 weight 4 QES values we invoke i.e. $QES(11,22;-12), QES(13;-12), QES(23;-12)$ are not among 10 MZV-dependent ones, that is, these QBESs are established without MZV usage.\\

\subsection{Related problems}
We present several minor problems with basically increasing complexity.\\

\subsubsection{Symmetric PLI}

Denote $D=\{(x,y): 0<x<1, 0<y<\frac{1-x}{x+1}\}$ a region symmetric about $y=x$. Let $D_1=D\cap\{(x,y): y>x\}$, we have\\
 {\small $$-\int_0^1 \frac{\text{Li}_2\left(\frac{1-x}{x+1}\right) \log (1-x)}{x} \, dx=\int_{D} \frac{\log (1-x) \log (1-y)}{x y} dxdy$$

  $$=2\int_{D_1} \frac{\log (1-x) \log (1-y)}{x y} dxdy=2\int_0^{\sqrt{2}-1} \frac{\left(\text{Li}_2(x)-\text{Li}_2\left(\frac{1-x}{x+1}\right)\right) \log (1-x)}{x} \, dx$$} Hence by using $\int \frac{\text{Li}_2(x) \log (1-x)}{x} \, dx=-\frac{\text{Li}_2(x){}^2}{2}$ and value of $PLI(1,0,0;25;1)$:\\
{\small $$\int_0^{\sqrt{2}-1} \frac{\text{Li}_2\left(\frac{1-x}{x+1}\right) \log (1-x)}{x} \, dx=-\frac{1}{2} \text{Li}_2\left(\sqrt{2}-1\right){}^2-\frac{\text{Li}_4\left(\frac{1}{2}\right)}{2}-\frac{7}{16} \zeta (3) \log (2)+\frac{7 \pi ^4}{2880}-\frac{1}{48} \log ^4(2)+\frac{1}{48} \pi ^2 \log ^2(2)$$}\\
A similar consideration for integrand $\frac{\text{Li}_n(x) \text{Li}_n(y)}{x y}$ yields the general formula:\\
{\small $$\left(\int_0^{\sqrt{2}-1}-\int_{\sqrt{2}-1}^1\right) \frac{\text{Li}_n(x) \text{Li}_{n+1}\left(\frac{1-x}{x+1}\right)}{x}dx=2\left(\int_0^{\sqrt{2}-1}-\int_{\sqrt{2}-1}^1\right) \frac{\text{Li}_n(x) \text{Li}_{n+1}\left(\frac{1-x}{x+1}\right)}{x^2-1}dx=\text{Li}_{n+1}\left(\sqrt{2}-1\right){}^2$$}\\

\subsubsection{MT of PLI}
We present a master formula of PLI:\\
$$\int_0^1 F(x) \left(\text{Li}_2\left(\frac{x-1}{x+1}\right)+\text{Li}_2\left(\frac{1-x}{x+1}\right)\right) \, dx=\int_0^1 \frac{\left(1-\sqrt{1-x^2}\right) F\left(\frac{1-\sqrt{1-x^2}}{x}\right) \text{Li}_2\left(\frac{1-x}{x+1}\right)}{2 x^2 \sqrt{1-x^2}} \, dx$$\\

For the proof, rewrite LHS as $\frac{1}{2} \int_0^1 F(x) \text{Li}_2\left(\frac{x^2-2 x+1}{x^2+2 x+1}\right) \, dx$ then let $\frac{1}{t}\to \frac{x^2+1}{2x}$. For instance, take $F=\frac{\text{Li}_2\left(\frac{1-x}{x+1}\right)}{x+1}$ and plug in $PLI(0,0,0;25,25/26;2)$ gives:\\
$$\int_0^1 \frac{\left(1-\sqrt{1-x^2}\right) \text{Li}_2\left(\frac{1-x}{x+1}\right) \text{Li}_2\left(\frac{x+\sqrt{1-x^2}-1}{x-\sqrt{1-x^2}+1}\right)}{x \sqrt{1-x^2} \left(x-\sqrt{1-x^2}+1\right)} \, dx$$ $$=8 \text{Li}_5\left(\frac{1}{2}\right)+8 \text{Li}_4\left(\frac{1}{2}\right) \log (2)-\frac{13 \pi ^2 \zeta (3)}{16}-\frac{\zeta (5)}{4}+\frac{7}{2} \zeta (3) \log ^2(2)+\frac{4 \log ^5(2)}{15}-\frac{2}{9} \pi ^2 \log ^3(2)+\frac{1}{36} \pi ^4 \log (2)$$\\

\subsubsection{Inverse ES}

Some inverse sums are computed below, where by inverse we mean the summand contains $\zeta (k)-H_n^{(k)}$ or $\eta (k)-\widetilde{H_n^{(k)}}$, that is, the tail of generalized harmonic sums.\\

Consider $\sum _{n=1}^{\infty } (-1)^n \left(\log (2)-\widetilde{H_n}\right)^5$. By Abel's summation (lift up $(-1)^n$ by partial summation and take difference of the other part) we have it equals to:\\
$$\sum _{n=1}^{\infty } \frac{1}{2} \left((-1)^{n-1}-1\right) \left(\left(\frac{(-1)^{n-1}}{n}+\log (2)-\widetilde{H_n}\right){}^5-\left(\log (2)-\widetilde{H_n}\right){}^5\right)$$\\
Expand the RHS, all terms are weight 5 ESs except another inverse sum $\sum _{n=1}^{\infty } \frac{\left(\log (2)-\widetilde{H_n}\right)^4}{n}$. Apply Abel summation again and split the summand, the new sum can be reduced to ESs either. Therefore we have:

{\small$$\sum _{n=1}^{\infty } (-1)^n \left(\log (2)-\widetilde{H_n}\right){}^5=10 \text{Li}_5\left(\frac{1}{2}\right)+10 \text{Li}_4\left(\frac{1}{2}\right) \log (2)-\frac{157 \zeta (5)}{16}+\frac{35}{8} \zeta (3) \log ^2(2)-\frac{2 \log ^5(2)}{3}-\frac{5}{18} \pi ^2 \log ^3(2)$$}

Using the same method we are able to deduce that:\\
\begin{equation}
\begin{aligned}
\sum _{n=0}^{\infty } \left(\log (2)-\widetilde{H_n}\right){}^5=-20 \text{Li}_4\left(\frac{1}{2}\right)-\frac{45}{4} \zeta (3) \log (2)+\frac{259 \pi ^4}{1440}+\frac{5 \log ^4(2)}{3}+\frac{5}{12} \pi ^2 \log ^2(2)
\end{aligned}
\notag
\end{equation}\\
Here the reason why a 'quintic' sum yields a 'quartic' value is the existence of $n$ after summation by parts once.  Combining 2 results and make use of $\int_0^1 \frac{x^n}{x+1} \, dx=(-1)^n \left(\log (2)-\widetilde{H_n}\right)$, we conclude:\\

\begin{equation}
\begin{aligned}
\sum _{n=0}^{\infty } \left(\log (2)-\widetilde{H_{2 n}}\right){}^5=\int_{(0,1)^5}\frac{dV}{(u+1) (v+1) (x+1) (y+1) (z+1) (1-(u v x y z)^2)}\\=-10 \text{Li}_4\left(\frac{1}{2}\right)+5 \text{Li}_5\left(\frac{1}{2}\right)+5 \text{Li}_4\left(\frac{1}{2}\right) \log (2)-\frac{157 \zeta (5)}{32}+\frac{35}{16} \zeta (3) \log ^2(2)\ \ \ \ \\-\frac{45}{8} \zeta (3) \log (2)+\frac{259 \pi ^4}{2880}+\frac{\log ^5(2)}{6}+\frac{5 \log ^4(2)}{6}-\frac{5}{36} \pi ^2 \log ^3(2)+\frac{5}{24} \pi ^2 \log ^2(2)\ \ \ 
\end{aligned}
\notag
\end{equation}\\

\subsubsection{High weight ES}
Consider the integral $\int_0^1 \frac{\log (x) \log \left(1-x^2\right) \sin ^{-1}(x)^4}{\sqrt{1-x^2}} \, dx$. Using $\sin ^{-1}(x)^4=\frac{3}{2} \sum _{n=1}^{\infty } \frac{H_{n-1}^{(2)} (2 x)^{2 n}}{n^2 \binom{2 n}{n}}$ \cite{ref8}, $t\to x^2$, Beta derivatives and that \\
$$\psi ^{(0)}(n+1)=H_n-\gamma, \psi ^{(0)}\left(n+\frac{1}{2}\right)=2 \left(H_{2 n}-H_n-\log (2)\right)$$ \\
The integral is converted to harmonic sum. On the other hand, by $x\to \sin(t)$ and expanding $2 \int_0^{\frac{\pi }{2}} t^4 (\log (2 \sin (t))-\log (2)) (\log (2 \cos (t))-\log (2)) \, dt$, the integral is also reduced to several LSIs with weight at most 6.  Evaluate them using weight 6 LI \cite{ref1} and contour deformation (subsection 5-1) gives:\\
{\small $$\sum _{n=1}^{\infty } \frac{H_{n-1}^{(2)} \left(2 \left(H_n-H_{2 n}+\log (2)\right) \left(H_n+2 \log (2)\right)+H_n^{(2)}-\frac{\pi ^2}{6}\right)}{n^2}$$ $$=-4 ES(1;-5)-\frac{4}{3} \pi ^2 \text{Li}_4\left(\frac{1}{2}\right)\\-\frac{15 \zeta (3)^2}{8}-\pi ^2 \zeta (3) \log (2)-\frac{1}{4} \zeta (5) \log (2)+\frac{281 \pi ^6}{15120}-\frac{1}{18} \pi ^2 \log ^4(2)+\frac{4}{45} \pi ^4 \log ^2(2)$$}\\

\subsubsection{High weight PLI}
Using multiple expansion the following 4 integrals are reducible to Zeta values:\\
{\small $$\int_0^1 \frac{\text{Li}_3(x) \log (x) \log (1-x)}{x} \, dx, \int_0^1 \frac{\text{Li}_3(x) \log (x) \log (1-x)}{1-x} \, dx, \int_0^1 \frac{\text{Li}_3(x) \log ^2(x)}{1-x} \, dx, \int_0^1 \frac{\text{Li}_3(x) \log ^2(1-x)}{x} \, dx$$}\\
By Prop. 7(7)(8) $\text{Li}_3\left(1-\frac{1}{x}\right)$ is expressible in $\text{Li}_3(x), \text{Li}_3(1-x)$ and log terms, so 4 corresponding integrals $\int_0^1 \frac{\text{Li}_3\left(1-\frac{1}{x}\right) \log (x) \log (1-x)}{x} \, dx, \cdots$ also Zeta expressible. Let $-t\to1-\frac{1}{x}$, we arrive at 4 integrals on $(0,\infty)$ with integrands composed of $\text{Li}_3(-x)$ and logs. Combining 2 of them gives\\
 $$\int_0^{\infty } \frac{\text{Li}_3(-x) \log (x+1) \log \left(\frac{1}{x}+1\right)}{x} \, dx=-2 \zeta (3)^2-\frac{31 \pi ^6}{5670}$$\\
  Split the integral, map back to (0,1), using Prop. 7(7) and weight 6 LI values \cite{ref1} we readily arrive at:\\
$$\int_0^1 \frac{\text{Li}_3(-x) \log (x+1) \log \left(\frac{1}{x}+1\right)}{x} \, dx$$
 $$=-ES(1;-5)+\frac{1}{3} \pi ^2 \text{Li}_4\left(\frac{1}{2}\right)+\frac{7}{24} \pi ^2 \zeta (3) \log (2)-\frac{137 \pi ^6}{90720}+\frac{1}{72} \pi ^2 \log ^4(2)-\frac{1}{72} \pi ^4 \log ^2(2)-\zeta (3)^2$$\\

\subsubsection{Parseval LSI}
We have an expanded logsine integral:\\

{\footnotesize $$\int_0^{2 \pi } \left(\frac{1}{12} \pi ^2 f(x)-\frac{1}{3} f(x)^3+\frac{1}{2} (x-\pi ) g(x)+h(x)\right)^2+\left(\frac{(\pi -x)}{2}  f(x)^2-f(x)g(x)+\frac{\pi ^2 (\pi -x)}{24} \right)^2 dx= 6 \pi  \zeta (3)^2+\frac{979 \pi ^7}{11340}$$}\\

Where $f(x)=\log \left(2 \sin \left(\frac{x}{2}\right)\right), g(x)=\text{Cl}_2(x), h(x)=\text{Cl}_3(x)$. And a compact one:\\

{\footnotesize $$\int_0^{2 \pi } \left| -\frac{1}{2} \text{Li}_2\left(e^{i x}\right){}^2-\log \left(1-e^{i x}\right) \text{Li}_3\left(e^{i x}\right)+\text{Li}_4\left(e^{i x}\right)\right|^2 \, dx=15 \pi  ES(2; 6)+\frac{4 \pi ^3 \zeta (3)^2}{3}-70 \pi  \zeta (5) \zeta (3)+\frac{677 \pi ^9}{113400}$$}\\

The first one is a direct consequence of Parseval identity applied to $f(x)=\sum _{n=1}^{\infty } \frac{1}{n} H_n^2 e^{i n x}$ on $(0, 2\pi)$, Prop. 3(13), closed-form of $ES(1,1,1,1;2)$ and elementary expansions, while the latter makes use of Prop. 4(5) and weight 8 non-alternating $ES(3,3;2)$. See \cite{ref15, ref20} for evaluation of 2 ESs respectively.\\

\subsubsection{Convolutional LI} 
Consider $\int _0^1\int _0^1\frac{\log ^n((1-x) (1-y))}{(x+1) (y+1)}dxdy$. By binomial theorem and Prop. 1(5) we may express this integral in terms of polylogs. On the other hand, let $u\to(1-x) (1-y), x\to x$, integrate w.r.t $x$, we find it equals to $\int_0^1 \frac{(2 \log (2-u)-\log (u)) \log ^n(u)}{4-u} \, du$. Using substitution $u\to 1-x$ and $\int_0^1 \frac{\log ^{n+1}(u)}{4-u} \, du=(-1)^{n-1} (n+1)! \text{Li}_{n+2}\left(\frac{1}{4}\right)$ we have\\
$$\int_0^1 \frac{\log (x+1) \log ^n(1-x)}{x+3} \, dx=\frac{(-1)^n (n+1)!}{2}  \left(\frac{1}{n+1}\sum _{k=0}^n \text{Li}_{k+1}\left(\frac{1}{2}\right) \text{Li}_{-k+n+1}\left(\frac{1}{2}\right)-\text{Li}_{n+2}\left(\frac{1}{4}\right)\right)$$\\

Similarly, consider evaluating $\int _0^1\int _0^1\frac{\log (x+1) \log (y+1) (\log(x+1)+\log(y+1))^n}{x y}dxdy$ by Prop. 1(8) and substitution yields: \\
$$\int_1^2 f(x) \log ^m(x) \, dx+\int_2^4 g(x) \log ^m(x) \, dx=m!\sum _{k=0}^m (k+1)(-k+m+1)J(k) J(m-k)$$\\
Where $J, f, g$ are \\
$$J(k)=-\sum _{j=0}^k \frac{\log ^j(2) \text{Li}_{-j+k+2}\left(\frac{1}{2}\right)}{j!}-\frac{\log ^{k+2}(2)}{(k+2) k!}+\zeta (k+2)$$ 

\begin{equation}
\begin{aligned}
f(x)=\frac{1}{6(x-1)}\left(12 \text{Li}_3\left(\frac{1}{x}\right)+12 \text{Li}_3(x)-6 \text{Li}_2(1-x) \log (x)+6 \text{Li}_2\left(\frac{1}{x}\right) \log (x) \right. \\ \left. -12 \text{Li}_2(x) \log (x)-6 \log (x-1) \log ^2(x)-6 i \pi  \log ^2(x)+\pi ^2 \log (x)-24 \zeta (3)\right)\ \ \ \
\end{aligned}
\notag
\end{equation}

\begin{equation}
\begin{aligned}
g(x)=\frac{1}{12(x-1)} \left(24 \text{Li}_3\left(\frac{x}{2}\right)+24 \text{Li}_3\left(\frac{2}{x}\right)-24 \text{Li}_2\left(\frac{x}{2}\right) \log \left(\frac{x}{2}\right)-12 \text{Li}_2\left(1-\frac{x}{2}\right) \log (x)\right. \\ \left. +12 \text{Li}_2\left(\frac{2}{x}\right) \log \left(\frac{x}{4}\right)-24 \text{Li}_3(2)+12 \text{Li}_2(2) \log (4)+6 \log ^2(2) \log \left(\frac{x}{2}\right)+12 i \pi  \log ^2(2)\ \ \right. \\ \left. -12 \log \left(1-\frac{x}{2}\right) \log ^2\left(\frac{x}{2}\right)+12 \log (2) \log \left(\frac{x}{2}\right) \log \left(\frac{x-2}{x}\right)-21 \zeta (3)+2 \log ^3(2)\right)\ \ \ \
\end{aligned}
\notag
\end{equation}\\
One may furtherly simplify Mathematica-generated $F,G$. Note that both integrals share the structure of 'convolution'.\\

\subsubsection{Valean's hexagram}

We generalize Valean's results \cite{ref17}. Denote

{\footnotesize $$X=(-1)^k \left(1-2^{k-2 n}\right) \zeta (-k+2 n+1)+\sum _{j=0}^{-k+2 n-1} \frac{(-1)^{j-1} \zeta (j+2) (i \pi )^{-j-k+2 n-1}}{(-j-k+2 n-1)!},\ Y=\Re\left(\frac{X}{(2 i)^{2 n-k}}\right)-\frac{\log (2) \left(\frac{\pi }{2}\right)^{2 n-k}}{(2 n-k)!}$$}\\
{\small $$A=2^{1-4 n} (2n-2)! \left(2^{2 n+1} \sum _{k=1}^{n-1} \left(1-2^{1-2 k}\right) \zeta (2 k) \zeta (-2 k+2n+1)-\left(2n \left(2^{2n}-1\right) -2^{2 n+1}+1\right) \zeta (2 n+1)\right)$$

$$C_1=(2 n-2)! \zeta (2 n+1),\ C_2=\left(2^{-2 n}-1\right) (2 n-2)! \zeta (2 n+1), \ U=2^{-2 (2 n-2)-3} \left(1-2^{2 n}\right) (2 n-2)! \zeta (2 n+1)$$}
\begin{equation}
\begin{aligned}
P=(2 n-2)! \left(-2^{-2 n} \left(1-2^{1-2 n}\right) \log (2) \zeta (2 n)-\sum _{k=0}^{2 n-1} \beta (k+1) \beta (2 n-k)\right. \\ \left.-\sum _{k=1}^{2 n-2} \left(1-2^{-k}\right) 2^{-2 n-1} \left(1-2^{k-2 n+1}\right) \zeta (k+1) \zeta (2 n-k)+2 n \zeta (2 n+1)\right)
\end{aligned}
\notag
\end{equation}

\begin{equation}
\begin{aligned}
Q=-2^{-4 n} (2n-2)! \left(-\left(2^{2n}-2\right) \sum _{k=1}^{2 (n-1)} \zeta (k+1) \zeta (-k+2n)\ \ \ \ \ \ \ \ \ \ \ \ \right. \\ \left.+2^{2 n+1} \sum _{k=1}^{n-1} \left(1-2^{1-2 k}\right) \zeta (2 k) \zeta (-2 k+2n+1)+\left(-2 n+2^{2 n+1}-1\right) \zeta (2 n+1)\right)
\end{aligned}
\notag
\end{equation}\\
Then we have six general formulas ($B_n$ denotes Bernoulli numbers \cite{ref13}):\\
$$(1)\ V_1=\int_0^1 \frac{\text{Li}_2\left(\frac{2 x}{x^2+1}\right) \log ^{2 n-2}(x)}{x} \, dx=(-1)^{n-1} (2 n-2)! \sum _{k=0}^{2 n-1}\frac{ B_k (2 \pi )^k Y}{k!}$$

$$(2)\ V_2=\int_0^1 \frac{\text{Li}_2\left(-\frac{2 x}{x^2+1}\right) \log ^{2 n-2}(x)}{x} \, dx=(-1)^n (2 n-2)! \sum _{k=0}^{2 n-1}\frac{ B_k (2 \pi )^k (1-2^{1-k}) Y }{k!}$$

$$(3)\ W_1=\int_0^1 \frac{x \log (1-x) \log ^{2 n-1}(x)}{x^2+1} \, dx=-(2 n-1) \left(\frac{A}{8}-\frac{C_1}{2}+\frac{U}{8}+\frac{V_1}{4}\right)$$

$$(4)\ W_2=\int_0^1 \frac{x \log (x+1) \log ^{2 n-1}(x)}{x^2+1} \, dx=-(2 n-1) \left(\frac{A}{8}-\frac{C_2}{2}+\frac{U}{8}+\frac{V_2}{4}\right)$$

$$(5)\ T_1=\int_0^1 \frac{\log (1-x) \log \left(x^2+1\right) \log ^{2 n-2}(x)}{x} \, dx=-\log (2) (2 n-2)! \zeta (2 n)-\frac{2 W_1}{2 n-1}+P$$

$$(6)\ T_2=\int_0^1 \frac{\log (1+x) \log \left(x^2+1\right) \log ^{2 n-2}(x)}{x} \, dx=Q-T_1$$\\

We establish $V_1$ below. Following \cite{ref17}, by applying double integration, Fourier expansion of Poisson kernel and logsine functions, we have the formula:\\
$$\int_0^1 \frac{\text{Li}_2\left(\frac{2 x}{x^2+1}\right) \log ^{2 n-2}(x)}{x} \, dx=2 (-1)^{2 n-1} (2 n-1)! \int_0^{\frac{\pi }{2}} \log (\cos (t)) \left(\sum _{k=1}^{\infty } \frac{\sin (k t)}{k^{2 n-1}}\right) \, dt$$\\
Now, \cite{ref13} offers the Fourier expansion\\
 $$\sum _{k=1}^{\infty } \frac{\sin (k t)}{k^{2 n-1}}=\frac{(2 \pi )^{2 n-1} (-1)^n B_{2 n-1}\left(\frac{t}{2 \pi }\right)}{2 (2 n-1)!}$$\\
 where $B_n(x)=\sum _{k=0}^n B_k \binom{n}{k} x^{n-k}$ are Bernoulli polynomials. Plugging in these results we find that:\\
$$\cdots=\frac{(2 \pi )^{2 n-1} (-1)^{n-1} }{2 n-1}\sum _{k=0}^{2 n-1} \frac{B_k \binom{2 n-1}{k} \int_0^{\frac{\pi }{2}} t^{-k+2 n-1} \log (\cos (t)) \, dt}{(2 \pi )^{-k+2 n-1}}$$\\
Finally, by simple induction we have \\
$$\int \frac{\log (z+1) \log ^n(z)}{z} \, dz=\sum _{k=0}^n \frac{(-1)^{k-1} n! \text{Li}_{k+2}(-z) \log ^{n-k}(z)}{(n-k)!}$$ \\
Integrate $\frac{\log (z+1) \log ^n(z)}{z}$ along the upper unit circle (counter clockwise), parametrize $z=e^{2 i x}$, it is direct that $\int_0^{\frac{\pi }{2}} x^n \log (2 \cos (x)) \, dx=\Re\left(\frac{Z}{(2 i)^{n+1}}\right)$, where \\
$$Z=(-1)^{n-1} \left(1-2^{-n-1}\right) n! \zeta (n+2)+\sum _{j=0}^n \frac{(-1)^{j-1} n! (i \pi )^{n-j} \zeta (j+2)}{(n-j)!}$$\\
Making use of this result completes the evaluation of $V_1$. $V_2$ is established similarly. $W_1, W_2, T_1, T_2$ are relevantly easy (using Prop. 1, ES/PLI/QLI general formulas and hints in \cite{ref17}) so we omit the details. These six formulas hold to arbitrary order.\\

\subsubsection{Hypergeometric rhapsody ($_pF_q$ part 2)}
First of all, take a look at radical polylog ntegral $\int_0^1 \frac{\text{Li}_3\left(\sqrt{x}\right)}{\sqrt{x (1-x)}} \, dx$ i.e. $\int_0^1 \frac{\text{Li}_3(x)}{\sqrt{1-x^2}} \, dx, \int_0^{\frac{\pi }{2}} \text{Li}_3(\sin (t)) \, dt$. Expanding $\text{Li}_3$ and taking care of parity, we obtain a hypergeometric form. Alternatively, following subsection 6-10-4, by double integration one have $\int_0^1 \frac{\text{Li}_3(x)}{\sqrt{1-x^2}} \, dx=\int_0^1 \frac{\log ^2(y) \left(-\pi  \sqrt{1-y^2}+2 \sin ^{-1}(y)+\pi \right)}{2 \left(2 y \sqrt{1-y^2}\right)} \, dy$. Now let $y\to \sin(t), t\to 2\tan^{-1}(u)$ then plug in QLI values to express it by extended basis $B_4$. Equating 2 forms yields\\
{\small $$(1)\ \, _5F_4\left(\frac{1}{2},\frac{1}{2},\frac{1}{2},1,1;\frac{3}{2},\frac{3}{2},\frac{3}{2},\frac{3}{2};1\right)=8 \Im(\text{Li}_4(1+i))-\frac{1}{12} \pi  \log ^3(2)-\frac{3}{16} \pi ^3 \log (2)-\frac{1}{384} \left(\psi ^{(3)}\left(\frac{1}{4}\right)-\psi ^{(3)}\left(\frac{3}{4}\right)\right)$$}\\

Now we focuse on QLSIs. Consider $\int_0^1 \frac{\log ^3(x)}{\sqrt{2-x^2}} \, dx=\int_0^{\frac{\pi }{4}} \left(\log (2 \sin (t))-\frac{\log (2)}{2}\right)^3 \, dt$. The RHS reduces to QLSIs solved in subsection 5-1, while LHS is transformed into hypergeometric sum by expanding $\frac{1}{\sqrt{2-x^2}}$ into power series. After simplification one have\\
{\footnotesize $$(2)\ \, _5F_4\left(\frac{1}{2},\frac{1}{2},\frac{1}{2},\frac{1}{2},\frac{1}{2};\frac{3}{2},\frac{3}{2},\frac{3}{2},\frac{3}{2};\frac{1}{2}\right)=\frac{1}{2} \sqrt{2} \Im(\text{Li}_4(1+i))+\frac{1}{16} \pi  \sqrt{2} \zeta (3)+\frac{1}{768} \pi ^3 \sqrt{2} \log (2)-\frac{\sqrt{2} \left(\psi ^{(3)}\left(\frac{1}{4}\right)-\psi ^{(3)}\left(\frac{3}{4}\right)\right)}{12288}$$}\\
Similar consideration for $\int_0^1 \frac{\log ^2(x) \sin ^{-1}\left(\frac{x}{\sqrt{2}}\right)}{\sqrt{2-x^2}} \, dx=\int_0^{\frac{\pi }{4}} t \left(\log (2 \sin (t))-\frac{\log (2)}{2}\right)^2 \, dt$ yields (we may obtain the expansion for $\frac{\sin ^{-1}(x)}{\sqrt{1-x^2}}$ by differentiating $\sin ^{-1}(x)^2=\frac{1}{2} \sum _{n=1}^{\infty } \frac{(2 x)^{2 n}}{n^2 \binom{2 n}{n}}$):\\
$$(3)\ \, _5F_4\left(1,1,1,1,1;\frac{3}{2},2,2,2;\frac{1}{2}\right)=2 \pi  \Im(\text{Li}_3(1+i))+\frac{5 \text{Li}_4\left(\frac{1}{2}\right)}{2}-\frac{11 \pi ^4}{144}+\frac{5 \log ^4(2)}{48}-\frac{1}{6} \pi ^2 \log ^2(2)$$\\

For higher weights not all QLI/QLSIs can be computed analytically. Nevertheless, the method of expansion above can be naturally generalized here, yielding connection formulas between $_pF_q$ and Borwein LSIs (for notations of $\text{Ls}$ see \cite{ref10}, which is a bit different from our QLSIs). Either of these constants can be added to $B_5$ representing irreducible constants:\\
$$(4)\ \, _6F_5\left(\frac{1}{2},\frac{1}{2},\frac{1}{2},\frac{1}{2},\frac{1}{2},\frac{1}{2};\frac{3}{2},\frac{3}{2},\frac{3}{2},\frac{3}{2},\frac{3}{2};\frac{1}{2}\right)$$ $$=\frac{1}{96} \sqrt{2} C \log ^3(2)-\frac{1}{16} \sqrt{2} \log ^2(2) \Im(\text{Li}_3(1+i))+\frac{1}{4} \sqrt{2} \log (2) \Im(\text{Li}_4(1+i))-\frac{1}{48} \sqrt{2} \text{Ls}_5\left(\frac{\pi}{2}\right)$$ $$+\frac{1}{32} \pi  \sqrt{2} \zeta (3) \log (2)+\frac{1}{512} \pi  \sqrt{2} \log ^4(2)+\frac{\pi ^3 \sqrt{2} \log ^2(2)}{3072}-\frac{\sqrt{2} \log (2) \left(\psi ^{(3)}\left(\frac{1}{4}\right)-\psi ^{(3)}\left(\frac{3}{4}\right)\right)}{24576}$$

$$(5)\ \, _6F_5\left(1,1,1,1,1,1;\frac{3}{2},2,2,2,2;\frac{1}{2}\right)=-\frac{1}{2} \pi  C \log ^2(2)+2 \pi  \log (2) \Im(\text{Li}_3(1+i))$$ $$+\frac{4}{3}\text{Ls}_5^{(1)}\left(\frac{\pi}{2}\right)+\frac{5}{2} \text{Li}_4\left(\frac{1}{2}\right) \log (2)+\frac{35}{32} \zeta (3) \log ^2(2)+\frac{5 \log ^5(2)}{48}-\frac{5}{24} \pi ^2 \log ^3(2)-\frac{11}{144} \pi ^4 \log (2)$$\\

Moreover, multiply both sides of $\sinh ^{-1}\left(\sqrt{\frac{x}{8}}\right)^2=\frac{1}{2} \sum _{n=1}^{\infty } \frac{(-1)^{n-1} x^n}{n^2 2^n \binom{2 n}{n}}$ by $\frac{\log ^2(x)}{x}$ then integrate on $(0,1)$. Let $t\to \sinh ^{-1}\left(\sqrt{\frac{x}{8}}\right)$ for LHS, use Prop. 1(1) for RHS, we deduce $\sum _{k=1}^{\infty } \frac{(-1)^{k-1}}{k^5 2^k \binom{2 k}{k}}=2 \int_0^{\frac{\log (2)}{2}} x^2 \coth (x) \log ^2\left(8 \sinh ^2(x)\right) \, dx$. Clearly it boils down to $\int_0^{\frac{\log (2)}{2}} x^2 \coth (x) \log ^k(\sinh (x)) \, dx$ for $k=0, 1, 2$, all of which are directly solved under substitution $x\to \frac{1}{2} \log (y+1)$ and several $LI(0,a,b;1/2)$ values. After simplification we arrive at\\
 $$(6)\ \, _6F_5\left(1,1,1,1,1,1;\frac{3}{2},2,2,2,2;-\frac{1}{8}\right)$$ $$=40 \text{Li}_5\left(\frac{1}{2}\right)+24 \text{Li}_4\left(\frac{1}{2}\right) \log (2)-38 \zeta (5)+4 \zeta (3) \log ^2(2)+\frac{19 \log ^5(2)}{30}-\frac{4}{9} \pi ^2 \log ^3(2)+\frac{7}{45} \pi ^4 \log (2)$$\\

Furtherly, consider calculating $\sum _{n=1}^{\infty } \frac{\binom{2 n}{n}}{n^6 4^n}$ via Beta derivatives like in section of BES,  and evaluating $\int_0^1 \frac{\log ^5(x)}{\sqrt{1-x^2}} \, dx$, $\int_0^1 \frac{\log ^4(x) \sin ^{-1}(x)}{\sqrt{1-x^2}} \, dx$ (essentially IBES $\sum _{n=1}^{\infty } \frac{4^n}{n^6 \binom{2 n}{n}}$) and $\sum _{n=1}^{\infty } \frac{(-1)^{n-1}}{n^6 2^n \binom{2 n}{n}}$ in 2 ways using methods above, with the help of \cite{ref1} one may deduce: \\
{\small $$(7)\ \, _8F_7\left(1,1,1,1,1,1,1,\frac{3}{2};2,2,2,2,2,2,2;1\right)$$ $$=-4 \zeta (3)^2-\frac{16}{3} \zeta (3) \log ^3(2)+\frac{4}{3} \pi ^2 \zeta (3) \log (2)-24 \zeta (5) \log (2)+\frac{79 \pi ^6}{7560}-\frac{8}{45}  \log ^6(2)+\frac{2}{9} \pi ^2 \log ^4(2)+\frac{1}{10} \pi ^4 \log ^2(2)$$}
$$(8)\ \, _7F_6\left(\frac{1}{2},\frac{1}{2},\frac{1}{2},\frac{1}{2},\frac{1}{2},\frac{1}{2},\frac{1}{2};\frac{3}{2},\frac{3}{2},\frac{3}{2},\frac{3}{2},\frac{3}{2},\frac{3}{2};1\right)$$ $$=\frac{\pi ^3 \zeta (3)}{192}+\frac{3 \pi  \zeta (5)}{32}+\frac{1}{16} \pi  \zeta (3) \log ^2(2)+\frac{1}{240} \pi  \log ^5(2)+\frac{1}{288} \pi ^3 \log ^3(2)+\frac{19 \pi ^5 \log (2)}{11520}$$

$$(9)\ \, _7F_6\left(1,1,1,1,1,1,1;\frac{3}{2},2,2,2,2,2;1\right)$$ {\small$$=2 ES(1;-5)+16 \text{Li}_6\left(\frac{1}{2}\right)-\zeta (3)^2+\pi ^2 \zeta (3) \log (2)+\frac{31}{8} \zeta (5) \log (2)-\frac{19 \pi ^6}{1120}+\frac{\log ^6(2)}{45}+\frac{1}{9} \pi ^2 \log ^4(2)-\frac{19}{360} \pi ^4 \log ^2(2)$$}
$$(10)\ \, _7F_6\left(1,1,1,1,1,1,1;\frac{3}{2},2,2,2,2,2;-\frac{1}{8}\right)=8 ES(1;-5)+64 \text{Li}_6\left(\frac{1}{2}\right)+24 \text{Li}_5\left(\frac{1}{2}\right) \log (2)$$ $$+4 \zeta (3)^2-\frac{4}{3} \zeta (3) \log ^3(2)+38 \zeta (5) \log (2)-\frac{43 \pi ^6}{540}-\frac{19}{180}  \log ^6(2)+\frac{1}{9} \pi ^2 \log ^4(2)-\frac{7}{90} \pi ^4 \log ^2(2)$$\\
All 4 general hypergeometric sums are generated by $A_6$, in which methods for non-alternating (7)(8) can be generalized to $_{n+1}F_n\left(\{1\}_n,\frac{3}{2};\{2\}_n;1\right),\  _{n+1}F_n\left(\{\frac{1}{2}\}_{n+1};\{\frac{3}{2}\}_n;1\right)$. As a corollary, all radical polylog integrals of form $\int_0^1 \frac{\text{Li}_{n-1}(x)}{\sqrt{x (1-x)}} \, dx$ i.e. $\int_0^1 \frac{\text{Li}_{n-1}\left(x^2\right)}{\sqrt{1-x^2}} \, dx, \int_0^{\frac{\pi }{2}} \text{Li}_{n-1}\left(\sin ^2(t)\right) \, dt$ are trivial since up to a constant factor it's equivalent to $_{n+1}F_n\left(\{1\}_n,\frac{3}{2};\{2\}_n;1\right)$ by expanding numerators. Moreover, the shifted $\int_0^1 \frac{\text{Li}_{n-1}(x)}{x \sqrt{1-x^2}} \, dx$ can be reduced to generalization of (9) i.e. $_{n+1}F_n\left(\{1\}_{n+1};\frac{3}{2},\{2\}_{n-1};1\right)$ modulo trivial $_pF_q$ values, which means these integrals are also solved if we plug in high weight LI \cite{ref1} then LSI values. For further examples, see (4)(5) of next subsection.\\

\subsubsection{FL fantasia ($_pF_q$ part 3)}
Apart from solving QBESs, we present some other use of FL theory (all FL expansions below are recorded in/deducible from \cite{ref23}). Firstly, in literature above 2 connection formulas between a QBES i.e. $\, _5F_4\left(\frac{1}{2},\frac{1}{2},\frac{1}{2},\frac{1}{2},\frac{1}{2};1,\frac{3}{2},\frac{3}{2},\frac{3}{2};1\right)$, $\sum _{n=1}^{\infty } \frac{(-1)^n \left(H_{2 n}-\frac{H_n}{2}\right){}^3}{2 n+1}$ and $\int_0^{\frac{1}{2}} \frac{\text{Li}_3(z)}{\sqrt{z (1-z)}} \, dz$ are proved, in which the middle sum is made up of $QES(\{11\}_k, \{31\}_{3-k};-12)$ $(0\leq k \leq 3)$ hence trivial. Therefore\\
{\small $$(1)\ \pi  \, _5F_4\left(\frac{1}{2},\frac{1}{2},\frac{1}{2},\frac{1}{2},\frac{1}{2};1,\frac{3}{2},\frac{3}{2},\frac{3}{2};1\right)=64 \Im(\text{Li}_4(1+i))-\frac{2}{3}  \pi  \log ^3(2)-\pi ^3 \log (2)-\frac{1}{32} \left(\psi ^{(3)}\left(\frac{1}{4}\right)-\psi ^{(3)}\left(\frac{3}{4}\right)\right)$$}\\
is confirmed true. Making use of another relation and value of $\int_0^1 \frac{\text{Li}_{3}(x)}{\sqrt{x (1-x)}} \, dx$ above yields a pair of by-products:\\
$$(2)\ \int_0^{\frac{1}{2}} \frac{\text{Li}_3(1-z)}{\sqrt{z (1-z)}} \, dz=\frac{\pi ^2 C}{6}-C \log ^2(2)+56 \Im(\text{Li}_4(1+i))-4 \log (2) \Im(\text{Li}_3(1+i))$$ $$+\pi  \zeta (3)+\frac{1}{3} \pi  \log ^3(2)-\frac{25}{24} \pi ^3 \log (2)-\frac{7 \psi ^{(3)}\left(\frac{1}{4}\right)}{256}+\frac{7 \psi ^{(3)}\left(\frac{3}{4}\right)}{256}$$

$$(2')\ \int_0^{\frac{1}{2}} \frac{\text{Li}_3(z)}{\sqrt{z (1-z)}} \, dz=-\frac{\pi ^2 C}{6}+C \log ^2(2)-56 \Im(\text{Li}_4(1+i))+4 \log (2) \Im(\text{Li}_3(1+i))$$ $$+\pi  \zeta (3)+\pi  \log ^3(2)+\frac{17}{24} \pi ^3 \log (2)+\frac{7 \psi ^{(3)}\left(\frac{1}{4}\right)}{256}-\frac{7 \psi ^{(3)}\left(\frac{3}{4}\right)}{256}$$\\

Secondly, by orthogonal identity $\int_0^1 P_p(2 x-1) P_q(2 x-1) \, dx=\frac{\delta _{p,q}}{2 p+1}$ and following FL series\\
$$K(x)=\sum _{n=0}^{\infty } \frac{2 P_n(2 x-1)}{2 n+1}, E(x)=\sum _{n=0}^{\infty } -\frac{4 P_n(2 x-1)}{(2 n-1) (2 n+1) (2 n+3)},\ \text{Li}_3(x)=\sum _{n=0}^{\infty } A_n P_n(2n-1)$$

$$A_0=\zeta (3)-\frac{\pi ^2}{6}+1,\ A_1=\frac{\pi ^2}{4}-\frac{15}{8}$$

$$A_n (n\geq2)=\frac{2 (-1)^n (2 n+1) \left(\sum _{k=n+3}^{\infty } \frac{(-1)^k}{k^2}\right)}{n (n+1)}+\frac{\left(n^4+n^3+3 n^2+8 n+4\right) (2 n+1)}{n^3 (n+1)^3 (n+2)^2}$$\\
We may compute polylog elliptic pair $\int_0^1 \text{Li}_3(x) K(x) \, dx,\int_0^1 \text{Li}_3(x) E(x) \, dx$ in closed-forms. The literature already computes the former so we only deal with the latter here. Modulo trivial rational sums, the representation $\sum _{k=n+3}^{\infty } \frac{(-1)^k}{k^2}=\int_0^1 \frac{(-x)^{n+2} \log (x)}{x+1} \, dx$ and constant factors, the problem boils down to evaluation of \\
$$\sum _{n=2}^{\infty } \frac{A_n}{(2 n-1) (2 n+1)^2 (2 n+3)}\ \text{then}\ \int_0^1 \frac{x^2 \log (x) }{x+1}\sum _{n=2}^{\infty } \frac{x^n}{n (n+1) (2 n-1) (2 n+1) (2 n+3)} \, dx$$\\
Let $x\to x^2$ and do partial fractions, the latter integral is decomposed into several NQLIs of form $\int_0^1 \frac{x^k \log (x) \log (1\pm x)}{x^2+1} \, dx$, all of which are trivially evaluated due to subsection 8-1-2 and QLI ($W\leq3$) results. Therefore we readily deduce:\\
{\small $$(3)\ \int_0^1 \text{Li}_3(x) E(x) \, dx=\frac{64}{3} C \log (2)-\frac{128}{3} \Im(\text{Li}_3(1+i))+\frac{4 \zeta (3)}{3}+\frac{4 \pi ^3}{3}-\frac{16 \pi ^2}{27}+\frac{2336}{81}+\frac{8}{3} \pi  \log ^2(2)-\frac{1120 \log (2)}{27}$$}\\

Moreover we have nonhomegeneous hypergeometric results. Indeed, by $\text{Li}_3(x)$ expansion above and that\\
$$\frac{1}{\sqrt[4]{x (1-x)}}=B\left(\frac{3}{4},\frac{3}{4}\right) \sum _{n=0}^{\infty } \frac{\binom{2 n}{n} P_{2 n}(2 x-1)}{4^n}$$ \\
The quartic-radical polylog integral $\int_0^1 \frac{\text{Li}_3(x)}{\sqrt[4]{x (1-x)}} \, dx$ is evaluable. Indeed, modulo trivial rational sums, the identity $-\sum _{k=2 n+3}^{\infty } \frac{4 (-1)^k}{k^2}=H_{n+1}^{(2)}-H_{n+\frac{1}{2}}^{(2)}$ and constant factors, the problem boils down to evaluation of \\
$$\sum _{n=1}^{\infty } \frac{A_{2n} \binom{2 n}{n}}{2^n (4 n+1)}\ \text{then} \ \sum _{n=1}^{\infty } \frac{\left(\frac{1}{2 n}-\frac{1}{2 n+1}\right) \binom{2 n}{n} \left(H_{n+1}^{(2)}-H_{n+\frac{1}{2}}^{(2)}\right)}{4^n}$$\\
Substituting $H_{n+1}^{(2)}-H_{n+\frac{1}{2}}^{(2)}=\int_0^1 \frac{\left(x-\sqrt{x}\right) x^n \log (x)}{1-x} \, dx$ and sum w.r.t $n$ in latter sum yields an elementary intrgral, which reduces to QLI/NQLIs ($W\leq3$) after $x\to \sin^2(t), t\to 2\tan^{-1}(u)$ and partial fractions. By subsection 8-1-2 again we obtain\\
$$(4)\ \int_0^1 \frac{\text{Li}_3(x)}{\sqrt[4]{x (1-x)}} \, dx=\frac{\Gamma \left(\frac{3}{4}\right)^2 }{\sqrt{\pi }}\ _5F_4\left(1,1,1,1,\frac{7}{4};2,2,2,\frac{5}{2};1\right)$$ 
\begin{equation}
\begin{aligned}
=\frac{2 \Gamma \left(\frac{3}{4}\right)^2}{\sqrt{\pi }} \left(2 \pi  C-4 C \log (2)+8 \Im(\text{Li}_3(1+i))-\frac{21 \zeta (3)}{8}-\frac{5 \pi ^2}{12}-\frac{3 \pi ^3}{16}-2 \pi \right. \\ \left.+8+\frac{\log ^3(2)}{6}-\frac{3}{4} \pi  \log ^2(2)+\log ^2(2)-\frac{5}{24} \pi ^2 \log (2)-\pi  \log (2)+4 \log (2)\right)
\end{aligned}
\notag
\end{equation}\\
In which the $_5F_4$ form comes from Taylor expansion of $\text{Li}_3(x)$ and Beta integrals. Now consider $\int_0^1 \frac{\text{Li}_3\left(\sqrt{x}\right)}{\sqrt[4]{x (1-x)}} \, dx$: Expand $\text{Li}_3\left(\sqrt{x}\right)$ as above, taking care of parity yields a Gamma combination of known $_5F_4\left(1,1,1,1,\frac{7}{4};2,2,2,\frac{5}{2};1\right)$ and $\, _5F_4\left(\frac{1}{2},\frac{1}{2},\frac{1}{2},\frac{5}{4},1;2,\frac{3}{2},\frac{3}{2},\frac{3}{2};1\right)$. According to parameter pairs $(1;2)$, $\left(\frac{1}{2}; \frac{3}{2}\right)$ in latter $_5F_4$, via partial fractions it's furtherly decomposed into\\
$$\, _3F_2\left(\frac{1}{2},\frac{5}{4},1;2,\frac{3}{2};1\right),\, _3F_2\left(\frac{1}{2},\frac{1}{2},\frac{5}{4};\frac{3}{2},\frac{3}{2};1\right),\, _4F_3\left(\frac{1}{2},\frac{1}{2},\frac{1}{2},\frac{5}{4};\frac{3}{2},\frac{3}{2},\frac{3}{2};1\right)$$\\
Using limiting argument for the first and Beta derivative for latter two, all three sums hence $\, _5F_4\left(\frac{1}{2},\frac{1}{2},\frac{1}{2},\frac{5}{4},1;2,\frac{3}{2},\frac{3}{2},\frac{3}{2};1\right)$ is evaluated. Thus plugging back the $_5F_4$ in (4) gives\\
\begin{equation}
\begin{aligned}
(5)\ \int_0^1 \frac{\text{Li}_3\left(\sqrt{x}\right)}{\sqrt[4]{x (1-x)}} \, dx=\pi  \sqrt{2}+\frac{\Gamma \left(\frac{3}{4}\right)^2}{32 \sqrt{\pi }} \left(48 \pi  C-32 C \log (2)+64 \Im(\text{Li}_3(1+i))-21 \zeta (3)-\frac{\pi ^3}{2}\ \ \ \right. \\ \left.-\frac{34 \pi ^2}{3} -48 \pi+64+\frac{4 \log ^3(2)}{3}-10 \pi  \log ^2(2)+8 \log ^2(2)-\frac{17}{3} \pi ^2 \log (2)-24 \pi  \log (2)+32 \log (2)\right)
\end{aligned}
\notag
\end{equation}\\
Here, by limiting argument we mean to apply following asymptotics:\\
$$_3F_2\left(\frac{1}{2},\frac{5}{4},1;2,\frac{3}{2};1\right)=\underset{z\to1^-}\lim z \, _3F_2\left(\frac{1}{2},1,\frac{5}{4};2,\frac{3}{2};z\right)=\underset{z\to1^-}\lim\left(2 z \, _2F_1\left(\frac{1}{2},\frac{5}{4};\frac{3}{2};z\right)-\frac{4}{\sqrt[4]{1-z}}+4\right)$$

$$_2F_1\left(\frac{1}{2},\frac{5}{4};\frac{3}{2};z\right)=\frac{2 }{\sqrt[4]{1-z}}-\frac{\pi ^{3/2}}{\sqrt{2} \Gamma \left(\frac{1}{4}\right) \Gamma \left(\frac{5}{4}\right)}+O((1-z)^{\frac34})\  (z\to1^-)$$\\
And by Beta derivative we mean to choose appropriate $r, p, a$ in:\\
$$_{r+2}F_{r+1}\left(\{\frac{n+1}{2}\}_{r+1}, \frac{1-p}{2}; \{\frac{n+3}{2}\}_{r+1};a\right)=\frac{(-1)^r (n+1)^{r+1}}{r!} \int_0^1 t^n \left(1-a t^2\right)^{\frac{p-1}{2}} \log ^r(t) \, dt$$\\
Similarly, with the help of\\
$$\frac{1}{(x (1-x))^{3/4}}=B\left(\frac{1}{4},\frac{1}{4}\right) \sum _{n=0}^{\infty } \frac{\binom{2 n}{n}(4n+1) P_{2 n}(2 x-1)}{4^n}$$\\
One may deduce a pair of dual result:\\
{\small$$(4')\ \int_0^1 \frac{\text{Li}_3(x)}{(x (1-x))^{3/4}} \, dx=\frac{\Gamma \left(\frac{1}{4}\right)^2 }{2 \sqrt{\pi }}\, _5F_4\left(1,1,1,1,\frac{5}{4};\frac{3}{2},2,2,2;1\right)$$
$$=\frac{\Gamma \left(\frac{1}{4}\right)^2 }{\sqrt{\pi }}\left(2 \pi  C+4 C \log (2)-8 \Im(\text{Li}_3(1+i))-\frac{21 \zeta (3)}{8}+\frac{3 \pi ^3}{16}+\frac{\log ^3(2)}{6}+\frac{3}{4} \pi  \log ^2(2)-\frac{5}{24} \pi ^2 \log (2)\right)$$

$$(5')\ \int_0^1 \frac{\text{Li}_3\left(\sqrt{x}\right)}{(x (1-x))^{3/4}} \, dx$$ $$=\frac{\Gamma \left(\frac{1}{4}\right)^2}{\sqrt{\pi }} \left(\frac{3 \pi  C}{4}+\frac{1}{2} C \log (2)-\Im(\text{Li}_3(1+i))-\frac{21 \zeta (3)}{64}+\frac{\pi ^3}{128}+\frac{\log ^3(2)}{48}+\frac{5}{32} \pi  \log ^2(2)-\frac{17}{192} \pi ^2 \log (2)\right)$$}\\
Using exactly the same method in (4)(5). Note that while expanding $\text{Li}_3\left(\sqrt{x}\right)$ in $(5')$ we have a Gamma conbination of known $_5F_4\left(1,1,1,1,\frac{5}{4};\frac{3}{2},2,2,2;1\right)$ and $\, _4F_3\left(\frac{1}{2},\frac{1}{2},\frac{1}{2},\frac{3}{4};\frac{3}{2},\frac{3}{2},\frac{3}{2};1\right)$, in which the latter, again, is trivially evaluated by Beta derivative.\\\\

Finally we have another application of (4). Let $x\to 1-x$ in Prop. 7(8), multiple both sides by $\frac{1}{\sqrt[4]{x(1-x)}}$, integrate on (0,1). Since $\int_0^1 \frac{\log ^k(x)\log^m(1-x)}{\sqrt[4]{x (1-x)}} \, dx,\int_0^1 \frac{\text{Li}_3(x)}{\sqrt[4]{x (1-x)}} \, dx,\int_0^1 \frac{\text{Li}_3(1-x)}{\sqrt[4]{x (1-x)}} \, dx$ are trivial due to Beta derivatives, formula (4) and reflection, we extract the value of $\int_0^1 \frac{\text{Li}_3\left(1-\frac{1}{x}\right)}{\sqrt[4]{x (1-x)}} \, dx$. Now a substitution $x\to \cos^2(t)$ yields\\
{\footnotesize $$(6)\ \int_0^{\frac{\pi }{2}} \sqrt{\sin (2 t)} \text{Li}_3\left(-\tan ^2(t)\right) \, dt=\sqrt{\frac{2}{\pi }} \Gamma \left(\frac{3}{4}\right)^2 \left(-4 \pi  C+8 C \log (2)-16 \Im(\text{Li}_3(1+i))-\frac{31 \zeta (3)}{4}+\pi ^3+\pi  \log ^2(2)\right)$$}\\
Note the corresponding derivation of (4') is:\\
{\small $$(6')\ \int_0^{\frac{\pi }{2}} \frac{\text{Li}_3\left(-\tan ^2(t)\right)}{\sqrt{\sin (2 t)}} \, dt=\frac{\Gamma \left(\frac{1}{4}\right)^2}{\sqrt{2 \pi }} \left(-2 \pi  C-4 C \log (2)+8 \Im(\text{Li}_3(1+i))-\frac{31 \zeta (3)}{8}-\frac{\pi ^3}{2}-\frac{1}{2} \pi  \log ^2(2)\right)$$}\\

Other nonhomegeneous radical polylog integrals (i.e. hypergeometric series), including the double parametric class $\int_0^1 x^{m-\frac{1}{2}} (1-x)^{n-\frac{1}{2}} \text{Li}_3(x) \, dx$ where $m,n\in \mathbb{N}$, are also solvable using FL expansion and results of NPLI/NQLIs. Since they are relevantly trivial, we simply pause here.\\

\subsection{Two conjectures solved}
We prove two related conjectures in \cite{ref1, ref3}.\\

\subsubsection{Borwein's conjecture on binomial sum ($_pF_q$ part 4)}
In \cite{ref3} Borwein conjectured that $\sum _{n=1}^{\infty } \frac{1}{n^4 2^n \binom{3 n}{n}}$ is expressible via single and multiple polylogarithms $\text{Li}_4, \text{Li}_{3,1}$. Based on $\text{Li}_4(z)+\text{Li}_{3,1}(z)=\sum _{n=1}^{\infty } \frac{H_n z^n}{n^3}$, Prop. 4(3) and subsection 5-3, his statement is equivalent to:\\
$$(*) \sum _{n=1}^{\infty } \frac{1}{n^4 2^n \binom{3 n}{n}}=-2 \pi  \Im(\text{Li}_3(1+i))-\frac{21 \text{Li}_4\left(\frac{1}{2}\right)}{2}-\frac{57}{8} \zeta (3) \log (2)+\frac{83 \pi ^4}{480}-\frac{23}{48}  \log ^4(2)+\frac{7}{12} \pi ^2 \log ^2(2)$$\\
Now we prove this identity. We have\\
{\footnotesize $$\sum _{n=1}^{\infty } \frac{1}{n^4 2^n \binom{3 n}{n}}=\sum _{n=1}^{\infty } \frac{\int_0^1 \frac{2 \left(x^2 (1-x)\right)^n}{x} \, dx}{n^3 2^n}=\int_0^1 \frac{2 \text{Li}_3\left(\frac{ x^2 (1-x)}{2}\right)}{x} \, dx=\int_0^1 \text{Li}_2\left(\frac{ x^2 (1-x)}{2}\right) \left(\frac{2 \log (x)}{1-x}-\frac{4 \log (x)}{x}\right) \, dx$$}\\
Where we've used Beta integral, Fubini, definition of polylog and integration by parts, lifting up $\log(x)$. Now the original problem is reduced to the evaluation of 2 integrals:\\
$$J_1\overset{\triangle}{=}\int_0^1 \frac{\text{Li}_2\left(\frac{1}{2} x^2 (1-x)\right) \log (x)}{x} \, dx, J_2\overset{\triangle}{=}\int_0^1 \frac{\text{Li}_2\left(\frac{1}{2} x^2 (1-x)\right) \log (x)}{1-x} \, dx$$
Indeed, we will establish that:\\
{\small $$J_1=-C^2+\frac{1}{4} \pi  C \log (2)-\frac{1}{2} \pi  \Im(\text{Li}_3(1+i))-\frac{15 \text{Li}_4\left(\frac{1}{2}\right)}{8}-\frac{83}{64} \zeta (3) \log (2)+\frac{31 \pi ^4}{768}-\frac{13}{192}  \log ^4(2)+\frac{1}{12} \pi ^2 \log ^2(2)$$

$$J_2=-2 C^2+\frac{1}{2} \pi  C \log (2)-2 \pi  \Im(\text{Li}_3(1+i))-9 \text{Li}_4\left(\frac{1}{2}\right)-\frac{197}{32} \zeta (3) \log (2)+\frac{107 \pi ^4}{640}-\frac{3}{8} \log ^4(2)+\frac{11}{24} \pi ^2 \log ^2(2)$$}\\
For $J_1$, lift up $\frac{1}{2} \log^2(x)$, simplify then factorize the log:\\
$$J_1=\int_0^1 \frac{1}{2} \left(\frac{2}{x}-\frac{1}{1-x}\right) \log ^2(x) \left(\log \left(\frac{ (1-x)^2+1}{2}\right)+\log (x+1)\right) \, dx$$\\
Via possible reflection $x\to 1-x$, this expression is decomposed into several LIs and QLIs, that is:\\
$$LI(0,2,1;1)-\frac{1}{2}LI(0,2,1;0)+\frac{\log(2)}{2}LI(2,0,0;1)-\frac{1}{2}QLI(114;2)+\int_0^1 \frac{\log ^2(1-x) \log \left(\frac{1}{2} \left(x^2+1\right)\right)}{1-x} \, dx$$
Lifting up $-\frac{1}{3}{\log^3(1-x)}$, the last one equals to $\frac{2}{3}QLI(111;5)$ so $J_1$ is solved. For $J_2$, firstly we lift up $\text{Li}_2(1-x)(=\int \frac{\log(x)}{1-x}dx)$ to arrive at:\\
$$J_2=\int_0^1 \left(\frac{2}{x}-\frac{1}{1-x}\right) \text{Li}_2(1-x) \left(\log \left(\frac{ (1-x)^2+1}{2}\right)+\log (x+1)\right) \, dx$$\\
Replacing $\text{Li}_2(1-x)$ by $\frac{\pi ^2}{6}-\text{Li}_2(x)-\log (1-x) \log (x)$, two parts of $J_2$ i.e. $\int_0^1 \frac{\text{Li}_2(1-x) \log (x+1)}{x} \, dx$, $\int_0^1 \frac{\text{Li}_2(1-x) \log (x+1)}{1-x} \, dx$ are evaluated by using integration by parts and LI, PLI results. For the other two parts involving $\log \left(\frac{ (1-x)^2+1}{2}\right)$, use $x\to 1-x$ again. Modulo known LI/PLI/QLIs, the only two remains are:\\
$$\int_0^1 \frac{\text{Li}_2(x) \log \left(x^2+1\right)}{x} \, dx, \int_0^1 \frac{\text{Li}_2(x) \log \left(\frac{1}{2} \left(x^2+1\right)\right)}{1-x} \, dx$$\\
The former is simply $QPLI4(1;4;2)$ and the latter boils down to $QPLI4(1;1;5)$ after lifting up $-\log(1-x)$ (modulo QLI values), thus $J_2$ is solved and Borwein's conjecture is confirmed true.\\

Note that the equivalent hypergeometric form is:\\
$$\, _6F_5\left(1,1,1,1,1,\frac{3}{2};\frac{4}{3},\frac{5}{3},2,2,2;\frac{2}{27}\right)$$ $$=-12 \pi  \Im(\text{Li}_3(1+i))-63 \text{Li}_4\left(\frac{1}{2}\right)-\frac{171}{4} \zeta (3) \log (2)+\frac{83 \pi ^4}{80}-\frac{23 }{8} \log ^4(2)+\frac{7}{2} \pi ^2 \log ^2(2)$$\\

\subsubsection{Au's conjecture on general LI}
In \cite{ref1} Au conjectured that $I_n=\int_0^1 \frac{\log^n (1-x) \log^{n-1} (1+x)}{1+x} dx$ lies in the algebra generated by $\log(2), \pi^2, $ $\zeta(3), \zeta(5), \zeta(7), \cdots$ over $\mathbb{Q}$ and verified the case $n\leq10$. Now we prove the general statement by generalizing method of  Prop. 2(7). Denote \\
$$f(k,j)=\int_0^{\frac{1}{2}} \frac{\log ^j(1-y) \log ^k(y)}{1-y} \, dy$$\\
Consider the integral $jf(j,k-1)$. By partial integration, separation, Beta derivatives and reflection $y\to 1-y$, we have\\
 $$jf(k,j-1)= -(-\log (2))^{j+k}+ k \int_0^{\frac{1}{2}} \frac{\log ^j(1-y) \log ^{k-1}(y)}{y} \, dy$$ $$ =-(-\log (2))^{j+k}+ k \left(\int_0^{1}-\int_{\frac{1}{2}}^1 \right) \frac{\log ^j(1-y) \log ^{k-1}(y)}{y} \, dy=U(j,k)-kf(j,k-1)$$\\
Where\\
$$U(k,j)= jf(k,j-1)+kf(j,k-1)=-(-\log(2))^{j+k}+ k \underset{a\to 0}{\text{lim}}\underset{b\to 1}{\text{lim}}\frac{\partial ^{k+j-1}B(a,b)}{\partial a^{k-1}\, \partial b^{j}}$$\\
Noticing $\frac{\binom{n-1}{j-1} \binom{n}{k}}{\binom{n}{j} \binom{n-1}{k-1}}=\frac{j}{k}$, we have a crucial identity\\
$$\binom{n}{k} \binom{n-1}{m-k} f(k,m-k)+\binom{n}{m+1-k} \binom{n-1}{k-1} f(m+1-k,k-1)=\frac{\binom{n}{k}\binom{n-1}{m-k} }{-k+m+1}U(k,m+1-k)$$\\
Now, using substitution $y\to\frac{1-x}{2}$ and Binomial theorem twice on $I_n$:\\
$$I_n=\int_0^{\frac{1}{2}} \frac{\log ^n(2 y) \log ^{n-1}(2 (1-y))}{1-y} \, dy=\sum _{k=0}^n \sum _{j=0}^{n-1} \binom{n}{k} \binom{n-1}{j} f(k,j) \log ^{2n-j-k-1}(2)$$\\
Extract $k=0$\\
$$I_n=\sum _{k=1}^n \sum _{j=0}^{n-1} \binom{n}{k} \binom{n-1}{j} f(k,j) \log ^{2n-j-k-1}(2)+\int_0^{\frac{1}{2}} \frac{\log ^n(2) \log ^{n-1}(2 (1-y))}{1-y} \, dy$$\\
Take Cauchy product and simplify: \\
$$I_n=\sum _{m=1}^{2n-1} \sum _{k+j=m}\binom{n}{k} \binom{n-1}{j} f(k,j) \log ^{2n-m-1}(2)+\frac{\log ^{2 n}(2)}{n}$$ $$=\left( \sum _{m=1}^n \sum _{k=1}^m+ \sum _{m=n+1}^{2 n-1} \sum _{k=m-n+1}^n\right) \binom{n}{k} \binom{n-1}{m-k} f(k,m-k) \log ^{2n-m-1}(2) +\frac{\log ^{2 n}(2)}{n}$$\\
 Let $k\to m+1-k$ then take averages: \\
 $$I_n=\frac{\log ^{2 n}(2)}{n}+\frac{1}{2}\left( \sum _{m=1}^n \sum _{k=1}^m+ \sum _{m=n+1}^{2 n-1} \sum _{k=m-n+1}^n\right) $$ $$\left(\binom{n}{k} \binom{n-1}{m-k} f(k,m-k)+\binom{n}{m+1-k} \binom{n-1}{k-1} f(m+1-k,k-1)\right) \log ^{2n-m-1}(2) $$\\
Plug in the crucial identity\\
$$(*)I_n=\frac{\log ^{2 n}(2)}{n}+\frac{1}{2}\left( \sum _{m=1}^n \sum _{k=1}^m+ \sum _{m=n+1}^{2 n-1} \sum _{k=m-n+1}^n\right) \frac{\binom{n}{k}\binom{n-1}{m-k} \log ^{-m+2 n-1}(2) }{-k+m+1}  U(k,m+1-k)$$\\
According to \cite{ref1}, Lemma $2.3$, all Beta derivatives in $U(k,m+1-k)$ above lies in the algebra generated by $\pi^2, $ $\zeta(3), \zeta(5), \zeta(7), \cdots$ over $\mathbb{Q}$. Adding up $\log(2)$ factors in $U$ concludes the proof.\\

\clearpage

\section{Table of subjects}

\subsection{LI}

Table of 85 LIs with weight $\leq5$. For instance:\\
\begin{center}
\begin{tabular}{|c|c|c|c|}
\hline
Weight 3 LI & $\pi^2 \log(2)$ & $\log^3(2)$ & $\zeta(3)$ \\ \hline
(1,0,1;2)& $-\frac{1}{12}$ & $\frac{1}{3}$ & $\frac{1}{8}$ \\ \hline
\end{tabular}
\end{center}
Implies that $LI(1,0,1;2)=\int_0^1 \frac{\log (1-x) \log (x+1)}{x+1} \, dx=\frac{\zeta (3)}{8}+\frac{\log ^3(2)}{3}-\frac{1}{12} \pi ^2 \log (2)$.\\

$$\text{Weight 1 LI:}\ LI(0,0,0;1)=\log(2)$$

\begin{center}
\begin{tabular}{|c|c|c|}
\hline
Weight 2 LI & $\pi^2 $ & $\log^2(2)$ \\ \hline
(0,1,0;0) & $-\frac{1}{6}$ & $0$ \\ \hline
(1,0,0;1) & $-\frac{1}{6}$ & $0$ \\ \hline
(0,0,1;1) & $\frac{1}{12}$ & $0$ \\ \hline
(1,0,0;2) & $-\frac{1}{12}$ & $\frac{1}{2}$ \\ \hline
(0,1,0;2) & $-\frac{1}{12}$ & $0$ \\ \hline
(0,0,1;2) & $0$ & $\frac{1}{2}$ \\ \hline
\end{tabular}
\end{center}

\begin{center}
\begin{tabular}{|c|c|c|c|}
\hline
Weight 3 LI & $\pi^2 \log(2)$ & $\log^3(2)$ & $\zeta(3)$ \\ \hline
(0,2,0;0)& $0$ & $0$ & $2$ \\ \hline
(1,1,0;0)& $0$ & $0$ & $1$ \\ \hline
(0,1,1;0)& $-\frac{1}{4}$ & $0$ & $1$ \\ \hline
(2,0,0;1)& $0$ & $0$ & $2$ \\ \hline
(0,0,2;1)& $0$ & $0$ & $\frac{1}{4}$ \\ \hline
(1,1,0;1)& $0$ & $0$ & $1$ \\ \hline
(0,1,1;1)& $0$ & $0$ & $-\frac{3}{4}$ \\ \hline
(1,0,1;1)& $0$ & $0$ & $-\frac{5}{8}$ \\ \hline
(2,0,0;2)& $-\frac{1}{6}$ & $\frac{1}{3}$ & $\frac{7}{4}$ \\ \hline
(0,2,0;2)& $0$ & $0$ & $\frac{3}{2}$ \\ \hline
(0,0,2;2)& $0$ & $\frac{1}{3}$ & $0$ \\ \hline
(1,1,0;2)& $-\frac{1}{4}$ & $0$ & $\frac{13}{8}$ \\ \hline
(0,1,1;2)& $0$ & $0$ & $-\frac{1}{8}$ \\ \hline
(1,0,1;2)& $-\frac{1}{12}$ & $\frac{1}{3}$ & $\frac{1}{8}$ \\ \hline
\end{tabular}
\end{center}

\begin{center}
\begin{tabular}{|c|c|c|c|c|c|}
\hline
Weight 4 LI & $\pi^4$ & $\pi^2 \log^2(2)$ & $\log^4(2)$ & $\zeta(3) \log(2)$ & $\text{Li}_4(\frac{1}{2})$ \\ \hline
(0,3,0;0) & $-\frac{1}{15}$ & $0$ & $0$ & $0$ & $0$ \\ \hline
(2,1,0;0) & $-\frac{1}{45}$ & $0$ & $0$ & $0$ & $0$ \\ \hline
(1,2,0;0) & $-\frac{1}{180}$ & $0$ & $0$ & $0$ & $0$ \\ \hline
(0,2,1;0) & $-\frac{19}{720}$ & $0$ & $0$ & $\frac{7}{2}$ & $0$ \\ \hline
(0,1,2;0) & $-\frac{7}{144}$ & $-\frac{5}{12}$ & $\frac{1}{6}$ & $\frac{21}{4}$ & $4$ \\ \hline
(1,1,1;0) & $\frac{17}{1440}$ & $-\frac{1}{24}$ & $-\frac{1}{12}$ & $\frac{7}{8}$ & $-2$ \\ \hline
(3,0,0;1) & $-\frac{1}{15}$ & $0$ & $0$ & $0$ & $0$ \\
\hline
\end{tabular}
\end{center}

{\small
\begin{center}
\begin{tabular}{|c|c|c|c|c|c|c|c|c|}
\hline
Weight 5 LI & $\pi^4 \log(2)$ & $\pi^2 \log^3(2)$ & $\log^5(2)$ & $\pi^2 \zeta(3)$ & $\zeta(3) \log^2(2)$ & $\zeta(5)$ & $\text{Li}_4(\frac{1}{2}) \log(2)$ & $\text{Li}_5(\frac{1}{2})$ \\ \hline
(0,4,0;0) & $0$ & $0$ & $0$ & $0$ & $0$ & $24$ & $0$ & $0$ \\ \hline
(3,1,0;0) & $0$ & $0$ & $0$ & $0$ & $0$ & $6$ & $0$ & $0$ \\ \hline
(1,3,0;0) & $0$ & $0$ & $0$ & $-1$ & $0$ & $12$ & $0$ & $0$ \\ \hline
(0,3,1;0) & $-\frac{1}{8}$ & $0$ & $0$ & $-\frac{3}{8}$ & $0$ & $12$ & $0$ & $0$ \\ \hline
(0,1,3;0) & $-\frac{17}{80}$ & $-\frac{1}{3}$ & $-\frac{1}{20}$ & $-\frac{7}{16}$ & $\frac{63}{8}$ & $\frac{117}{4}$ & $-6$ & $-24$ \\ \hline
(2,2,0;0) & $0$ & $0$ & $0$ & $-\frac{2}{3}$ & $0$ & $8$ & $0$ & $0$ \\ \hline
(0,2,2;0) & $-\frac{49}{360}$ & $\frac{1}{9}$ & $-\frac{1}{5}$ & $-\frac{13}{24}$ & $\frac{7}{2}$ & $\frac{47}{2}$ & $-8$ & $-16$ \\
\hline
\end{tabular}
\end{center}}

\clearpage

\subsection{ES} 

Table of 89 ESs with weight $\leq5$.\\

$$\text{Weight 2 ESs:} \ ES(1;-1)=\frac{\pi ^2}{12}-\frac{\log ^2(2)}{2}, ES(-1;-1)=\frac{\pi ^2}{12}+\frac{\log ^2(2)}{2}$$\\

\begin{center}
\begin{tabular}{|c|c|c|c|}
\hline
Weight 3 ES & $\pi^2 \log(2)$ & $\log^3(2)$ & $\zeta(3)$ \\ \hline
$(1;2)$ & $0$ & $0$ & $2$ \\ \hline
$(-1;2)$ & $\frac{1}{4}$ & $0$ & $-\frac{1}{4}$ \\ \hline
$(1;-2)$ & $0$ & $0$ & $\frac{5}{8}$ \\ \hline
$(-1;-2)$ & $\frac{1}{4}$ & $0$ & $-\frac{5}{8}$ \\ \hline
$(2;-1)$ & $-\frac{1}{12}$ & $0$ & $1$ \\ \hline
$(-2;-1)$ & $-\frac{1}{6}$ & $0$ & $\frac{13}{8}$ \\ \hline
$(1,1;-1)$ & $-\frac{1}{12}$ & $\frac{1}{3}$ & $\frac{3}{4}$ \\ \hline
$(1,-1;-1)$ & $\frac{1}{6}$ & $-\frac{1}{3}$ & $\frac{1}{4}$ \\ \hline
$(-1,-1;-1)$ & $\frac{1}{4}$ & $\frac{1}{3}$ & $-\frac{1}{2}$ \\ \hline
\end{tabular}
\end{center}

\begin{center}
\begin{tabular}{|c|c|c|c|c|c|}
\hline
Weight 4 ES & $\pi^4$ & $\pi^2 \log^2(2)$ & $\log^4(2)$ & $\zeta(3) \log(2)$ & $\text{Li}_4(\frac{1}{2})$ \\ \hline
$(1;3)$ & $\frac{1}{72}$ & $0$ & $0$ & $0$ & $0$ \\ \hline
$(-1;3)$ & $-\frac{1}{288}$ & $0$ & $0$ & $\frac{7}{4}$ & $0$ \\ \hline
$(1;-3)$ & $\frac{11}{360}$ & $\frac{1}{12}$ & $-\frac{1}{12}$ & $-\frac{7}{4}$ & $-2$ \\ \hline
$(-1;-3)$ & $\frac{1}{60}$ & $\frac{1}{12}$ & $-\frac{1}{12}$ & $0$ & $-2$ \\ \hline
$(2;2)$ & $\frac{7}{360}$ & $0$ & $0$ & $0$ & $0$ \\ \hline
$(-2;2)$ & $\frac{17}{288}$ & $\frac{1}{6}$ & $-\frac{1}{6}$ & $-\frac{7}{2}$ & $-4$ \\ \hline
$(2;-2)$ & $-\frac{17}{480}$ & $-\frac{1}{6}$ & $\frac{1}{6}$ & $\frac{7}{2}$ & $4$ \\
\hline
\end{tabular}
\end{center}

{\scriptsize\begin{center}
\begin{tabular}{|c|c|c|c|c|c|c|c|c|}
\hline
Weight 5 ES & $\pi^4 \log(2)$ & $\pi^2 \log^3(2)$ & $\log^5(2)$ & $\pi^2 \zeta(3)$ & $\zeta(3) \log^2(2)$ & $\zeta(5)$ & $\text{Li}_4(\frac{1}{2}) \log(2)$ & $\text{Li}_5(\frac{1}{2})$ \\ \hline
$(1;4)$ & $0$ & $0$ & $0$ & $-\frac{1}{6}$ & $0$ & $3$ & $0$ & $0$ \\ \hline
$(-1;4)$ & $\frac{1}{48}$ & $0$ & $0$ & $\frac{1}{16}$ & $0$ & $-\frac{17}{16}$ & $0$ & $0$ \\ \hline
$(1;-4)$ & $0$ & $0$ & $0$ & $-\frac{1}{12}$ & $0$ & $\frac{59}{32}$ & $0$ & $0$ \\ \hline
$(-1;-4)$ & $\frac{1}{48}$ & $0$ & $0$ & $\frac{1}{8}$ & $0$ & $-\frac{59}{32}$ & $0$ & $0$ \\ \hline
$(2;3)$ & $0$ & $0$ & $0$ & $\frac{1}{2}$ & $0$ & $-\frac{9}{2}$ & $0$ & $0$ \\ \hline
$(-2;3)$ & $0$ & $0$ & $0$ & $-\frac{1}{24}$ & $0$ & $\frac{51}{32}$ & $0$ & $0$ \\ \hline
$(2;-3)$ & $0$ & $0$ & $0$ & $\frac{5}{48}$ & $0$ & $-\frac{11}{32}$ & $0$ & $0$ \\
\hline
\end{tabular}
\end{center}}

\clearpage

\subsection{PLI}

Table of 263 PLIs with weight $\leq5$. For sake of simplicity, we abbreviate $PLI(a,b,c;\cdots)$ as $PLI(abc;\cdots)$.\\

\begin{center}
\begin{tabular}{|c|c|c|c|}
\hline
Weight 3 PLI & $\pi^2 \log(2)$ & $\log^3(2)$ & $\zeta(3)$ \\ \hline
$(000;21;1)$ & $0$ & $0$ & $1$ \\ \hline
$(000;21;2)$ & $\frac{1}{6}$ & $0$ & $-\frac{5}{8}$ \\ \hline
$(000;22;1)$ & $0$ & $0$ & $-\frac{3}{4}$ \\ \hline
$(000;22;2)$ & $-\frac{1}{12}$ & $0$ & $\frac{1}{4}$ \\ \hline
$(000;23;0)$ & $-\frac{1}{12}$ & $\frac{1}{6}$ & $\frac{7}{8}$ \\ \hline
$(000;23;2)$ & $\frac{1}{12}$ & $-\frac{1}{6}$ & $-\frac{1}{4}$ \\ \hline
$(000;24;2)$ & $\frac{1}{12}$ & $-\frac{1}{6}$ & $\frac{1}{8}$ \\ \hline
$(000;25;0)$ & $-\frac{1}{6}$ & $0$ & $\frac{13}{8}$ \\ \hline
$(000;25;2)$ & $\frac{1}{6}$ & $0$ & $-\frac{5}{8}$ \\ \hline
$(000;26;0)$ & $\frac{1}{12}$ & $0$ & $-1$ \\ \hline
$(000;26;2)$ & $-\frac{1}{12}$ & $0$ & $\frac{1}{4}$ \\ \hline
\end{tabular}
\end{center}
$$$$ $$$$
\begin{center}
\begin{tabular}{|c|c|c|c|c|c|}
\hline
Weight 4 PLI & $\pi^4$ & $\pi^2 \log^2(2)$ & $\log^4(2)$ & $\zeta(3) \log(2)$ & $\text{Li}_4(\frac{1}{2})$ \\ \hline
$(010;21;0)$ & $-\frac{1}{120}$ & $0$ & $0$ & $0$ & $0$ \\ \hline
$(010;21;1)$ & $-\frac{1}{90}$ & $0$ & $0$ & $0$ & $0$ \\ \hline
$(010;21;2)$ & $-\frac{1}{480}$ & $0$ & $0$ & $0$ & $0$ \\ \hline
$(001;21;1)$ & $-\frac{1}{60}$ & $-\frac{1}{12}$ & $\frac{1}{12}$ & $\frac{7}{4}$ & $2$ \\ \hline
$(001;21;2)$ & $-\frac{1}{480}$ & $\frac{1}{12}$ & $0$ & $0$ & $0$ \\ \hline
$(100;21;1)$ & $-\frac{1}{72}$ & $0$ & $0$ & $0$ & $0$ \\ \hline
$(100;21;2)$ & $-\frac{29}{1440}$ & $-\frac{1}{24}$ & $\frac{1}{8}$ & $0$ & $3$ \\ \hline
$(010;22;0)$ & $\frac{71}{1440}$ & $\frac{1}{6}$ & $-\frac{1}{6}$ & $-\frac{7}{2}$ & $-4$ \\ \hline
$(010;22;1)$ & $\frac{7}{720}$ & $0$ & $0$ & $0$ & $0$ \\ \hline
$(010;22;2)$ & $\frac{13}{288}$ & $\frac{1}{6}$ & $-\frac{1}{6}$ & $-\frac{7}{2}$ & $-4$ \\ \hline
$(001;22;1)$ & $-\frac{1}{288}$ & $0$ & $0$ & $0$ & $0$ \\ \hline
$(001;22;2)$ & $\frac{1}{30}$ & $\frac{1}{12}$ & $-\frac{1}{8}$ & $-\frac{21}{8}$ & $-3$ \\ \hline
$(100;22;1)$ & $\frac{11}{360}$ & $\frac{1}{12}$ & $-\frac{1}{12}$ & $-\frac{7}{4}$ & $-2$ \\ \hline
$(100;22;2)$ & $\frac{1}{24}$ & $\frac{1}{8}$ & $-\frac{1}{6}$ & $-\frac{21}{8}$ & $-4$ \\ \hline
$(010;23;0)$ & $-\frac{1}{720}$ & $0$ & $-\frac{1}{24}$ & $\frac{1}{8}$ & $-1$ \\ \hline
$(010;23;2)$ & $\frac{1}{360}$ & $\frac{1}{24}$ & $0$ & $-1$ & $0$ \\ \hline
$(001;23;0)$ & $-\frac{1}{288}$ & $-\frac{1}{24}$ & $\frac{1}{24}$ & $\frac{7}{8}$ & $0$ \\ \hline
$(001;23;1)$ & $\frac{11}{240}$ & $\frac{1}{8}$ & $-\frac{1}{6}$ & $-\frac{13}{4}$ & $-4$ \\ \hline
$(001;23;2)$ & $\frac{1}{30}$ & $\frac{1}{6}$ & $-\frac{1}{6}$ & $-\frac{23}{8}$ & $-3$ \\ \hline
\end{tabular}
\end{center}

\begin{center}
\begin{tabular}{|c|c|c|c|c|c|}
\hline
Weight 4 PLI & $\pi^4$ & $\pi^2 \log^2(2)$ & $\log^4(2)$ & $\zeta(3) \log(2)$ & $\text{Li}_4(\frac{1}{2})$ \\ \hline
$(100;23;0)$ & $0$ & $0$ & $0$ & $0$ & $-1$ \\ \hline
$(100;23;1)$ & $-\frac{19}{1440}$ & $\frac{1}{24}$ & $\frac{1}{24}$ & $\frac{1}{4}$ & $1$ \\ \hline
$(100;23;2)$ & $-\frac{1}{1440}$ & $\frac{1}{24}$ & $-\frac{1}{24}$ & $-\frac{1}{4}$ & $0$ \\ \hline
$(010;24;0)$ & $-\frac{11}{288}$ & $-\frac{1}{24}$ & $\frac{1}{8}$ & $1$ & $3$ \\ \hline
$(010;24;2)$ & $\frac{1}{180}$ & $0$ & $-\frac{1}{12}$ & $-\frac{1}{8}$ & $-2$ \\ \hline
$(001;24;1)$ & $-\frac{1}{36}$ & $-\frac{5}{24}$ & $\frac{1}{6}$ & $\frac{23}{8}$ & $4$ \\ \hline
$(001;24;2)$ & $-\frac{1}{90}$ & $0$ & $0$ & $1$ & $1$ \\
\hline
\end{tabular}
\end{center}

{\small \begin{center}
\begin{tabular}{|c|c|c|c|c|c|c|c|c|}
\hline
Weight 5 PLI & $\pi^4 \log(2)$ & $\pi^2 \log^3(2)$ & $\log^5(2)$ & $\pi^2 \zeta(3)$ & $\zeta(3) \log^2(2)$ & $\zeta(5)$ & $\text{Li}_4(\frac{1}{2}) \log(2)$ & $\text{Li}_5(\frac{1}{2})$ \\ \hline
$(020;21;0)$ & $0$ & $0$ & $0$ & $1$ & $0$ & $-11$ & $0$ & $0$ \\ \hline
$(020;21;1)$ & $0$ & $0$ & $0$ & $0$ & $0$ & $2$ & $0$ & $0$ \\ \hline
$(020;21;2)$ & $0$ & $0$ & $0$ & $\frac{3}{4}$ & $0$ & $-\frac{67}{8}$ & $0$ & $0$ \\ \hline
$(002;21;1)$ & $0$ & $-\frac{2}{9}$ & $\frac{4}{15}$ & $-\frac{1}{4}$ & $\frac{7}{2}$ & $-\frac{39}{8}$ & $8$ & $8$ \\ \hline
$(002;21;2)$ & $0$ & $\frac{1}{9}$ & $-\frac{1}{15}$ & $\frac{7}{48}$ & $-\frac{7}{8}$ & $\frac{1}{4}$ & $-2$ & $-2$ \\ \hline
$(200;21;1)$ & $0$ & $0$ & $0$ & $\frac{1}{3}$ & $0$ & $-1$ & $0$ & $0$ \\ \hline
$(200;21;2)$ & $\frac{19}{720}$ & $-\frac{1}{36}$ & $\frac{1}{20}$ & $\frac{7}{24}$ & $0$ & $-\frac{23}{32}$ & $0$ & $-6$ \\
\hline
\end{tabular}
\end{center}

\begin{center}
\begin{tabular}{|c|c|c|c|c|c|c|c|c|}
\hline
Weight 5 PLI & $\pi^4 \log(2)$ & $\pi^2 \log^3(2)$ & $\log^5(2)$ & $\pi^2 \zeta(3)$ & $\zeta(3) \log^2(2)$ & $\zeta(5)$ & $\text{Li}_4(\frac{1}{2}) \log(2)$ & $\text{Li}_5(\frac{1}{2})$ \\ \hline
$(200;23;1)$ & $-\frac{1}{144}$ & $\frac{5}{36}$ & $-\frac{11}{60}$ & $\frac{11}{24}$ & $-\frac{11}{4}$ & $\frac{81}{32}$ & $-6$ & $-8$ \\ \hline
$(200;23;2)$ & $-\frac{1}{720}$ & $\frac{7}{36}$ & $-\frac{13}{60}$ & $\frac{23}{48}$ & $-\frac{23}{8}$ & $\frac{23}{32}$ & $-6$ & $-6$ \\ \hline
$(110;23;0)$ & $\frac{17}{1440}$ & $\frac{1}{24}$ & $-\frac{1}{12}$ & $\frac{1}{96}$ & $-\frac{15}{16}$ & $\frac{25}{64}$ & $-2$ & $0$ \\ \hline
$(110;23;1)$ & $-\frac{3}{160}$ & $\frac{1}{12}$ & $-\frac{7}{60}$ & $\frac{3}{8}$ & $-\frac{11}{8}$ & $\frac{143}{64}$ & $-4$ & $-6$ \\ \hline
$(110;23;2)$ & $-\frac{1}{96}$ & $\frac{5}{36}$ & $-\frac{1}{8}$ & $\frac{5}{12}$ & $-\frac{17}{8}$ & $\frac{37}{64}$ & $-4$ & $-5$ \\ \hline
$(011;23;0)$ & $-\frac{17}{360}$ & $\frac{7}{72}$ & $-\frac{17}{120}$ & $\frac{1}{12}$ & $-\frac{1}{2}$ & $\frac{245}{32}$ & $-5$ & $-8$ \\ \hline
$(011;23;1)$ & $\frac{1}{24}$ & $-\frac{1}{18}$ & $\frac{1}{10}$ & $\frac{1}{8}$ & $\frac{3}{8}$ & $-\frac{79}{8}$ & $4$ & $8$ \\
\hline
\end{tabular}
\end{center}

\begin{center}
\begin{tabular}{|c|c|c|c|c|c|c|c|c|}
\hline
Weight 5 PLI & $\pi^4 \log(2)$ & $\pi^2 \log^3(2)$ & $\log^5(2)$ & $\pi^2 \zeta(3)$ & $\zeta(3) \log^2(2)$ & $\zeta(5)$ & $\text{Li}_4(\frac{1}{2}) \log(2)$ & $\text{Li}_5(\frac{1}{2})$ \\ \hline
$(011;25;0)$ & $\frac{1}{180}$ & $-\frac{1}{18}$ & $\frac{1}{10}$ & $-\frac{5}{96}$ & $\frac{21}{16}$ & $-\frac{93}{16}$ & $4$ & $8$ \\ \hline
$(011;25;1)$ & $0$ & $\frac{1}{18}$ & $-\frac{1}{15}$ & $\frac{1}{48}$ & $-\frac{7}{8}$ & $\frac{31}{32}$ & $-2$ & $-2$ \\ \hline
$(011;25;2)$ & $-\frac{1}{360}$ & $-\frac{1}{12}$ & $\frac{1}{12}$ & $-\frac{3}{32}$ & $\frac{7}{16}$ & $\frac{31}{64}$ & $2$ & $0$ \\ \hline
$(101;25;0)$ & $\frac{29}{1440}$ & $-\frac{1}{18}$ & $\frac{1}{30}$ & $\frac{11}{96}$ & $\frac{3}{4}$ & $-\frac{109}{16}$ & $2$ & $6$ \\ \hline
$(101;25;1)$ & $0$ & $\frac{1}{18}$ & $-\frac{1}{15}$ & $-\frac{17}{96}$ & $-\frac{7}{8}$ & $\frac{31}{8}$ & $-2$ & $-2$ \\ \hline
$(101;25;2)$ & $\frac{1}{480}$ & $0$ & $\frac{1}{15}$ & $-\frac{5}{96}$ & $\frac{1}{4}$ & $-\frac{45}{32}$ & $2$ & $2$ \\ \hline
$(020;26;0)$ & $\frac{49}{720}$ & $-\frac{1}{6}$ & $\frac{7}{30}$ & $\frac{1}{8}$ & $\frac{7}{4}$ & $-\frac{527}{32}$ & $8$ & $12$ \\
\hline
\end{tabular}
\end{center}

\begin{center}
\begin{tabular}{|c|c|c|c|c|c|c|c|c|}
\hline
Weight 5 PLI & $\pi^4 \log(2)$ & $\pi^2 \log^3(2)$ & $\log^5(2)$ & $\pi^2 \zeta(3)$ & $\zeta(3) \log^2(2)$ & $\zeta(5)$ & $\text{Li}_4(\frac{1}{2}) \log(2)$ & $\text{Li}_5(\frac{1}{2})$ \\ \hline
$(001;32;2)$ & $0$ & $\frac{1}{18}$ & $-\frac{1}{15}$ & $\frac{1}{96}$ & $-\frac{5}{4}$ & $\frac{125}{64}$ & $-2$ & $-2$ \\ \hline
$(010;33;0)$ & $\frac{1}{720}$ & $\frac{1}{36}$ & $-\frac{1}{40}$ & $\frac{1}{12}$ & $-\frac{1}{2}$ & $-\frac{1}{32}$ & $-1$ & $-2$ \\ \hline
$(010;33;2)$ & $-\frac{1}{720}$ & $\frac{1}{18}$ & $-\frac{3}{40}$ & $\frac{1}{16}$ & $-\frac{13}{16}$ & $\frac{33}{64}$ & $-2$ & $-1$ \\ \hline
$(100;33;0)$ & $0$ & $0$ & $0$ & $0$ & $0$ & $0$ & $0$ & $-1$ \\ \hline
$(100;33;1)$ & $\frac{19}{1440}$ & $\frac{5}{72}$ & $-\frac{13}{120}$ & $\frac{3}{32}$ & $-\frac{23}{16}$ & $\frac{23}{64}$ & $-3$ & $-2$ \\ \hline
$(100;33;2)$ & $-\frac{1}{288}$ & $\frac{5}{72}$ & $-\frac{11}{120}$ & $\frac{7}{48}$ & $-\frac{7}{8}$ & $\frac{81}{64}$ & $-3$ & $-3$ \\ \hline
$(001;33;0)$ & $0$ & $\frac{1}{36}$ & $-\frac{1}{30}$ & $\frac{7}{96}$ & $-\frac{7}{16}$ & $\frac{27}{32}$ & $-1$ & $-2$ \\
\hline
\end{tabular}
\end{center}

\begin{center}
\begin{tabular}{|c|c|c|c|c|c|c|c|c|}
\hline
Weight 5 PLI & $\pi^4 \log(2)$ & $\pi^2 \log^3(2)$ & $\log^5(2)$ & $\pi^2 \zeta(3)$ & $\zeta(3) \log^2(2)$ & $\zeta(5)$ & $\text{Li}_4(\frac{1}{2}) \log(2)$ & $\text{Li}_5(\frac{1}{2})$ \\ \hline
$(000;43;2)$ & $0$ & $\frac{1}{36}$ & $-\frac{1}{30}$ & $\frac{7}{96}$ & $-\frac{7}{16}$ & $\frac{27}{32}$ & $-1$ & $-2$ \\ \hline
$(000;44;2)$ & $0$ & $0$ & $0$ & $0$ & $0$ & $1$ & $0$ & $-1$ \\ \hline
$(000;45;0)$ & $-\frac{1}{90}$ & $0$ & $0$ & $-\frac{1}{8}$ & $0$ & $\frac{91}{32}$ & $0$ & $0$ \\ \hline
$(000;45;2)$ & $\frac{1}{90}$ & $0$ & $0$ & $\frac{1}{8}$ & $0$ & $-\frac{59}{32}$ & $0$ & $0$ \\ \hline
$(000;46;0)$ & $\frac{7}{720}$ & $0$ & $0$ & $\frac{1}{16}$ & $0$ & $-2$ & $0$ & $0$ \\ \hline
$(000;46;2)$ & $-\frac{7}{720}$ & $0$ & $0$ & $-\frac{1}{16}$ & $0$ & $\frac{17}{16}$ & $0$ & $0$ \\ \hline
$(000;21,21;1)$ & $0$ & $0$ & $0$ & $\frac{1}{3}$ & $0$ & $-3$ & $0$ & $0$ \\
\hline
\end{tabular}
\end{center}

\begin{center}
\begin{tabular}{|c|c|c|c|c|c|c|c|c|}
\hline
Weight 5 PLI & $\pi^4 \log(2)$ & $\pi^2 \log^3(2)$ & $\log^5(2)$ & $\pi^2 \zeta(3)$ & $\zeta(3) \log^2(2)$ & $\zeta(5)$ & $\text{Li}_4(\frac{1}{2}) \log(2)$ & $\text{Li}_5(\frac{1}{2})$ \\ \hline
$(000;22,25;1)$ & $0$ & $\frac{1}{18}$ & $-\frac{1}{15}$ & $\frac{3}{32}$ & $-\frac{7}{8}$ & $\frac{31}{64}$ & $-2$ & $-2$ \\ \hline
$(000;22,25;2)$ & $-\frac{11}{720}$ & $\frac{5}{72}$ & $-\frac{11}{120}$ & $\frac{23}{96}$ & $-\frac{7}{8}$ & $\frac{93}{64}$ & $-3$ & $-4$ \\ \hline
$(000;22,26;0)$ & $-\frac{19}{240}$ & $\frac{7}{72}$ & $-\frac{7}{40}$ & $\frac{1}{32}$ & $0$ & $\frac{899}{64}$ & $-7$ & $-14$ \\ \hline
$(000;22,26;1)$ & $0$ & $\frac{1}{18}$ & $-\frac{1}{15}$ & $-\frac{1}{96}$ & $-\frac{7}{8}$ & $\frac{155}{64}$ & $-2$ & $-2$ \\ \hline
$(000;22,26;2)$ & $\frac{19}{240}$ & $-\frac{1}{24}$ & $\frac{13}{120}$ & $-\frac{1}{24}$ & $-\frac{7}{8}$ & $-\frac{93}{8}$ & $5$ & $12$ \\ \hline
$(000;23,24;0)$ & $-\frac{1}{144}$ & $-\frac{1}{18}$ & $\frac{1}{20}$ & $-\frac{7}{96}$ & $\frac{21}{16}$ & $-\frac{81}{64}$ & $3$ & $3$ \\ \hline
$(000;23,24;2)$ & $\frac{1}{144}$ & $\frac{1}{18}$ & $-\frac{1}{20}$ & $-\frac{25}{96}$ & $-\frac{21}{16}$ & $\frac{369}{64}$ & $-3$ & $-3$ \\
\hline
\end{tabular}
\end{center}}

\clearpage

\subsection{QLI}

Table of 193 QLIs with weight $\leq4$.\\

\begin{center}
\begin{tabular}{|c|c|c|}
\hline
Weight 2 QLI & $\pi^2$ & $\log^2(2)$ \\ \hline
$(4;2)$ & $\frac{1}{24}$ & $0$ \\ \hline
$(4;3)$ & $-\frac{1}{48}$ & $\frac{3}{4}$ \\ \hline 
$(5;4)$ & $\frac{1}{32}$ & $0$ \\ \hline 
$(1;5)$ & $-\frac{5}{96}$ & $\frac{1}{8}$ \\ \hline 
$(2;5)$ & $-\frac{1}{48}$ & $0$ \\ \hline 
$(3;5)$ & $\frac{1}{96}$ & $\frac{1}{8}$ \\ \hline 
$(4;5)$ & $0$ & $\frac{1}{4}$ \\ \hline 
\end{tabular}
\end{center}
$$$$ $$$$
\begin{center}
\begin{tabular}{|c|c|c|}
\hline
Weight 2 QLI & $\pi \log(2)$ & $C$ \\ \hline
$(5;2)$ & $0$ & $1$ \\ \hline 
$(5;3)$ & $\frac{1}{8}$ & $0$ \\ \hline 
$(1;4)$ & $\frac{1}{8}$ & $-1$ \\ \hline
$(2;4)$ & $0$ & $-1$ \\ \hline 
$(3;4)$ & $\frac{1}{8}$ & $0$ \\ \hline 
$(4;4)$ & $\frac{1}{2}$ & $-1$ \\ \hline 
$(5;5)$ & $-\frac{1}{8}$ & $\frac{1}{2}$ \\ \hline 
\end{tabular}
\end{center}

\begin{center}
\begin{tabular}{|c|c|c|c|c|}
\hline
Weight 3 QLI & $\pi^2 \log(2)$ & $\log^3 (2)$  & $\zeta(3)$ & $\pi C$ \\ \hline
$(24;1)$ & $-\frac{3}{16}$ & $0$  & $2$ & $-\frac{1}{2}$ \\ \hline
$(14;2)$ & $0$ & $0$  & $\frac{23}{32}$ & $-\frac{1}{2}$ \\ \hline
$(24;2)$ & $0$ & $0$  & $-\frac{3}{16}$ & $0$ \\ \hline
$(34;2)$ & $0$ & $0$  & $-\frac{33}{32}$ & $\frac{1}{2}$ \\ \hline
$(44;2)$ & $0$ & $0$  & $\frac{1}{8}$ & $0$ \\ \hline
$(55;2)$ & $0$ & $0$  & $-\frac{7}{8}$ & $\frac{1}{2}$ \\ \hline
$(14;3)$ & $-\frac{5}{32}$ & $\frac{11}{24}$  & $\frac{7}{4}$ & $-\frac{1}{2}$ \\ \hline
$(24;3)$ & $-\frac{1}{16}$ & $0$  & $\frac{3}{2}$ & $-\frac{1}{2}$ \\ \hline
$(34;3)$ & $-\frac{1}{96}$ & $\frac{11}{24}$  & $0$ & $0$ \\ \hline
$(44;3)$ & $-\frac{1}{24}$ & $\frac{2}{3}$  & $\frac{5}{2}$ & $-1$ \\ \hline
$(55;3)$ & $\frac{1}{32}$ & $0$  & $-\frac{21}{32}$ & $\frac{1}{4}$ \\ \hline
$(15;4)$ & $\frac{1}{64}$ & $0$  & $-\frac{7}{64}$ & $-\frac{1}{8}$ \\ \hline
$(25;4)$ & $0$ & $0$  & $\frac{7}{16}$ & $-\frac{1}{4}$ \\ \hline
$(35;4)$ & $\frac{1}{64}$ & $0$  & $\frac{21}{64}$ & $-\frac{1}{8}$ \\ \hline
$(45;4)$ & $\frac{1}{16}$ & $0$  & $\frac{21}{64}$ & $-\frac{1}{4}$ \\ \hline
$(11;5)$ & $-\frac{5}{96}$ & $\frac{1}{24}$  & $\frac{35}{32}$ & $0$ \\ \hline
$(12;5)$ & $-\frac{3}{32}$ & $0$  & $\frac{41}{64}$ & $0$ \\ \hline
$(13;5)$ & $-\frac{1}{48}$ & $\frac{1}{24}$  & $-\frac{7}{64}$ & $0$ \\ \hline
$(14;5)$ & $-\frac{5}{96}$ & $\frac{1}{12}$  & $\frac{11}{16}$ & $-\frac{1}{4}$ \\ \hline
$(22;5)$ & $0$ & $0$  & $\frac{3}{16}$ & $0$ \\ \hline
$(23;5)$ & $\frac{1}{32}$ & $0$  & $-\frac{15}{64}$ & $0$ \\ \hline
$(24;5)$ & $0$ & $0$  & $-\frac{1}{32}$ & $0$ \\ \hline
$(33;5)$ & $\frac{1}{96}$ & $\frac{1}{24}$  & $0$ & $0$ \\ \hline
$(34;5)$ & $\frac{1}{96}$ & $\frac{1}{12}$  & $-\frac{5}{8}$ & $\frac{1}{4}$ \\ \hline
$(44;5)$ & $0$ & $\frac{1}{6}$  & $0$ & $0$ \\ \hline
$(55;5)$ & $-\frac{1}{32}$ & $0$  & $-\frac{21}{64}$ & $\frac{1}{4}$ \\ \hline
\end{tabular}
\end{center}

\begin{center}
\begin{tabular}{|c|c|c|c|c|}
\hline
Weight 3 QLI & $\pi^3$ & $\pi \log^2 (2)$  & $C \log(2)$ & $P_3$ \\ \hline
$(25;1)$ & $\frac{3}{64}$ & $\frac{1}{8}$  & $\frac{1}{2}$ & $-2$ \\ \hline
$(15;2)$ & $0$ & $\frac{1}{16}$  & $\frac{1}{2}$ & $-1$ \\ \hline
$(25;2)$ & $-\frac{1}{32}$ & $0$  & $0$ & $0$ \\ \hline
$(35;2)$ & $\frac{3}{32}$ & $\frac{3}{16}$  & $\frac{3}{2}$ & $-3$ \\ \hline
$(45;2)$ & $\frac{1}{16}$ & $\frac{1}{8}$  & $1$ & $-2$ \\ \hline
$(15;3)$ & $\frac{19}{384}$ & $\frac{7}{32}$  & $\frac{1}{2}$ & $-2$ \\ \hline
$(25;3)$ & $-\frac{1}{64}$ & $0$  & $\frac{1}{2}$ & $0$ \\ \hline
$(35;3)$ & $\frac{7}{128}$ & $\frac{7}{32}$  & $1$ & $-2$ \\ \hline
$(45;3)$ & $\frac{5}{96}$ & $\frac{1}{4}$  & $1$ & $-2$ \\ \hline
$(11;4)$ & $\frac{7}{64}$ & $\frac{3}{16}$  & $0$ & $-2$ \\ \hline
$(12;4)$ & $\frac{3}{64}$ & $\frac{1}{16}$  & $0$ & $-1$ \\ \hline
$(13;4)$ & $-\frac{1}{32}$ & $0$  & $-1$ & $1$ \\ \hline
$(14;4)$ & $\frac{9}{64}$ & $\frac{1}{2}$  & $-\frac{1}{2}$ & $-4$ \\ \hline
$(22;4)$ & $\frac{1}{16}$ & $0$  & $0$ & $0$ \\ \hline
$(23;4)$ & $-\frac{5}{64}$ & $-\frac{3}{16}$  & $-2$ & $3$ \\ \hline
$(24;4)$ & $\frac{3}{32}$ & $\frac{1}{8}$  & $-1$ & $-2$ \\ \hline
$(33;4)$ & $-\frac{7}{64}$ & $-\frac{3}{16}$  & -2$$ & $4$ \\ \hline
$(34;4)$ & $-\frac{7}{64}$ & $0$  & $-\frac{5}{2}$ & $4$ \\ \hline
$(44;4)$ & $\frac{7}{48}$ & $\frac{5}{4}$  & $-2$ & $-4$ \\ \hline
$(55;4)$ & $\frac{1}{192}$ & $0$  & $0$ & $0$ \\ \hline
$(15;5)$ & $-\frac{19}{384}$ & $-\frac{1}{8}$  & $\frac{1}{4}$ & $1$ \\ \hline
$(25;5)$ & $-\frac{5}{64}$ & $-\frac{1}{8}$  & $0$ & $2$ \\ \hline
$(35;5)$ & $\frac{11}{384}$ & $0$  & $\frac{3}{4}$ & $-1$ \\ \hline
$(45;5)$ & $-\frac{7}{192}$ & $-\frac{1}{4}$  & $\frac{1}{2}$ & $1$ \\ \hline
\end{tabular}
\end{center}

\begin{center}
\begin{tabular}{|c|c|c|c|c|c|c|c|c|}
\hline
Weight 4 QLI & $\pi^4$ & $\pi^2 \log^2(2)$  & $\log^4(2)$ & $\zeta(3) \log(2)$ & $\text{Li}_4(\frac{1}{2})$ & $\pi P_3$ & $\pi C \log(2)$ & $C^2$ \\ \hline
$(124;1)$ & $-\frac{59}{3840}$ & $-\frac{5}{192}$  & $-\frac{5}{96}$ & $\frac{35}{64}$ & $-\frac{5}{4}$ & $\frac{1}{2}$ & $-\frac{1}{4}$ & $1$ \\ \hline
$(224;1)$ & $-\frac{199}{5760}$ & $0$  & $0$ & $\frac{35}{16}$ & $0$ & $0$ & $0$ & $2$ \\ \hline
$(234;1)$ & $\frac{211}{11520}$ & $-\frac{11}{64}$  & $\frac{1}{32}$ & $\frac{315}{64}$ & $\frac{3}{4}$ & $-\frac{3}{2}$ & $-\frac{1}{4}$ & $1$ \\ \hline
$(244;1)$ & $-\frac{61}{576}$ & $-\frac{17}{48}$  & $\frac{1}{24}$ & $\frac{77}{16}$ & $1$ & $2$ & $-2$ & $4$ \\ \hline
$(255;1)$ & $\frac{53}{5760}$ & $\frac{1}{24}$  & $-\frac{1}{24}$ & $-\frac{7}{16}$ & $-1$ & $0$ & $-\frac{1}{4}$ & $0$ \\ \hline
$(114;2)$ & $-\frac{521}{11520}$ & $-\frac{11}{96}$  & $\frac{5}{96}$ & $\frac{35}{32}$ & $\frac{5}{4}$ & $1$ & $-\frac{1}{2}$ & $2$ \\ \hline
$(124;2)$ & $-\frac{119}{5760}$ & $-\frac{5}{96}$  & $\frac{5}{96}$ & $\frac{35}{32}$ & $\frac{5}{4}$ & $0$ & $0$ & $1$ \\
\hline
\end{tabular}
\end{center}

\begin{center}
\begin{tabular}{|c|c|c|c|c|c|c|c|c|}
\hline
Weight 4 QLI & $\pi^4$ & $\pi^2 \log^2(2)$  & $\log^4(2)$ & $\zeta(3) \log(2)$ & $\text{Li}_4(\frac{1}{2})$ & $\pi P_3$ & $\pi C \log(2)$ & $C^2$ \\ \hline
$(115;4)$ & $\frac{143}{11520}$ & $\frac{11}{384}$  & $-\frac{1}{48}$ & $-\frac{35}{64}$ & $-\frac{1}{2}$ & $0$ & $-\frac{1}{8}$ & $\frac{1}{2}$ \\ \hline
$(125;4)$ & $\frac{121}{15360}$ & $\frac{1}{32}$  & $-\frac{1}{32}$ & $-\frac{7}{16}$ & $-\frac{3}{4}$ & $0$ & $-\frac{1}{8}$ & $\frac{1}{4}$ \\ \hline
$(135;4)$ & $-\frac{29}{15360}$ & $\frac{3}{128}$  & $-\frac{1}{32}$ & $-\frac{35}{64}$ & $-\frac{3}{4}$ & $\frac{1}{4}$ & $-\frac{1}{4}$ & $\frac{1}{4}$ \\ \hline
$(145;4)$ & $\frac{97}{5760}$ & $\frac{1}{12}$  & $-\frac{7}{192}$ & $-\frac{105}{128}$ & $-\frac{7}{8}$ & $-\frac{1}{4}$ & $-\frac{1}{4}$ & $\frac{1}{2}$ \\ \hline
$(225;4)$ & $\frac{151}{11520}$ & $\frac{1}{24}$  & $-\frac{1}{24}$ & $-\frac{7}{8}$ & $-1$ & $0$ & $0$ & $0$ \\ \hline
$(235;4)$ & $-\frac{119}{15360}$ & $0$  & $-\frac{1}{32}$ & $-\frac{7}{16}$ & $-\frac{3}{4}$ & $\frac{1}{2}$ & $-\frac{3}{8}$ & $\frac{1}{4}$ \\ \hline
$(245;4)$ & $\frac{437}{11520}$ & $\frac{11}{96}$  & $-\frac{1}{12}$ & $-\frac{7}{8}$ & $-2$ & $-\frac{1}{2}$ & $-\frac{1}{4}$ & $0$ \\
\hline
\end{tabular}
\end{center}

\begin{center}
\begin{tabular}{|c|c|c|c|c|c|c|c|c|}
\hline
Weight 4 QLI & $\pi^3 \log(2)$ & $\pi \log^3(2)$  & $\pi^2 C$ & $C \log^2(2)$ & $\pi \zeta(3)$ & $\log(2) P_3$ & $ P_4$ & $\Psi_3$ \\ \hline
$(125;1)$ & $\frac{3}{128}$ & $\frac{1}{32}$  & $\frac{11}{96}$ & $\frac{1}{8}$ & $\frac{35}{256}$ & $-\frac{1}{2}$ & $0$ & $-\frac{1}{1536}$ \\ \hline
$(225;1)$ & $\frac{1}{32}$ & $0$  & $-\frac{1}{24}$ & $0$ & $\frac{35}{64}$ & $0$ & $0$ & $-\frac{1}{768}$ \\ \hline
$(235;1)$ & $-\frac{1}{128}$ & $\frac{1}{32}$  & $-\frac{1}{96}$ & $\frac{1}{8}$ & $\frac{91}{256}$ & $-\frac{1}{2}$ & $0$ & $-\frac{1}{1536}$ \\ \hline
$(245;1)$ & $\frac{1}{64}$ & $\frac{1}{16}$  & $-\frac{5}{48}$ & $\frac{1}{4}$ & $\frac{105}{128}$ & $-1$ & $0$ & $-\frac{1}{768}$ \\ \hline
$(115;2)$ & $0$ & $\frac{1}{24}$  & $-\frac{5}{48}$ & $\frac{1}{4}$ & $\frac{35}{128}$ & $-1$ & $2$ & $0$ \\ \hline
$(125;2)$ & $0$ & $0$  & $-\frac{3}{16}$ & $0$ & $\frac{35}{128}$ & $0$ & $0$ & $\frac{1}{1536}$ \\ \hline
$(135;2)$ & $0$ & $\frac{1}{6}$  & $-\frac{5}{48}$ & $1$ & $\frac{7}{128}$ & $-4$ & $8$ & $-\frac{1}{256}$ \\
\hline
\end{tabular}
\end{center}

\begin{center}
\begin{tabular}{|c|c|c|c|c|c|c|c|c|}
\hline
Weight 4 QLI & $\pi^3 \log(2)$ & $\pi \log^3(2)$  & $\pi^2 C$ & $C \log^2(2)$ & $\pi \zeta(3)$ & $\log(2) P_3$ & $ P_4$ & $\Psi_3$ \\ \hline
$(111;4)$ & $\frac{21}{128}$ & $\frac{3}{32}$  & $0$ & $0$ & $0$ & $0$ & $-6$ & $0$ \\ \hline
$(112;4)$ & $\frac{3}{64}$ & $\frac{1}{48}$  & $\frac{1}{3}$ & $0$ & $0$ & $0$ & $-2$ & $-\frac{1}{768}$ \\ \hline
$(113;4)$ & $\frac{3}{128}$ & $-\frac{5}{96}$  & $\frac{1}{6}$ & $-1$ & $0$ & $2$ & $-4$ & $\frac{1}{768}$ \\ \hline
$(114;4)$ & $\frac{23}{64}$ & $\frac{7}{24}$  & $\frac{5}{48}$ & $-\frac{1}{4}$ & $\frac{35}{64}$ & $0$ & $-16$ & $\frac{1}{192}$ \\ \hline
$(122;4)$ & $\frac{1}{32}$ & $0$  & $\frac{1}{3}$ & $0$ & $0$ & $0$ & $0$ & $-\frac{1}{384}$ \\ \hline
$(123;4)$ & $-\frac{1}{64}$ & $-\frac{5}{48}$  & $\frac{1}{6}$ & $-1$ & $0$ & $2$ & $-2$ & $0$ \\ \hline
$(124;4)$ & $\frac{9}{64}$ & $\frac{1}{16}$  & $\frac{1}{4}$ & $-\frac{1}{4}$ & $\frac{7}{16}$ & $0$ & $-6$ & $0$ \\
\hline
\end{tabular}
\end{center}

\clearpage

\subsection{QPLI}

Table of 172 QPLIs with weight $\leq4$.\\

\begin{center}
\begin{tabular}{|c|c|c|c|c|}
\hline
Weight 3 QPLI & $\pi^3$ & $\pi \log^2 (2)$  & $C \log(2)$ & $P_3$ \\ \hline
$(1;4)$ & $\frac{1}{24}$ & $\frac{1}{16}$  & $\frac{1}{2}$ & $-1$ \\ \hline
$(2;4)$ & $\frac{7}{96}$ & $\frac{3}{16}$  & $\frac{3}{2}$ & $-3$ \\ \hline
$(3;4)$ & $-\frac{7}{384}$ & $-\frac{3}{32}$  & $-\frac{1}{2}$ & $1$ \\ \hline
$(4;4)$ & $\frac{35}{384}$ & $\frac{3}{32}$  & $\frac{1}{2}$ & $-2$ \\ \hline
$(5;4)$ & $\frac{1}{24}$ & $\frac{1}{16}$  & $\frac{1}{2}$ & $-1$ \\ \hline
$(6;4)$ & $\frac{7}{96}$ & $\frac{3}{16}$  & $\frac{3}{2}$ & $-3$ \\ \hline
$(7;4)$ & $\frac{5}{48}$ & $\frac{1}{4}$  & $2$ & $-4$ \\ \hline
$(8;2)$ & $\frac{1}{32}$ & $0$  & $0$ & $0$ \\ \hline
$(8;3)$ & $-\frac{3}{32}$ & $-\frac{3}{16}$  & $-\frac{1}{2}$ & $3$ \\ \hline
$(8;5)$ & $-\frac{1}{32}$ & $-\frac{1}{16}$  & $0$ & $1$ \\ \hline
\end{tabular}
\end{center}

\begin{center}
\begin{tabular}{|c|c|c|c|c|}
\hline
Weight 3 QPLI & $\pi^2 \log(2)$ & $\log^3 (2)$  & $\zeta(3)$ & $\pi C$ \\ \hline
$(1;5)$ & $\frac{1}{12}$ & $0$  & $\frac{23}{64}$ & $-\frac{1}{4}$ \\ \hline
$(2;5)$ & $-\frac{1}{24}$ & $0$  & $-\frac{33}{64}$ & $\frac{1}{4}$ \\ \hline
$(3;5)$ & $\frac{3}{64}$ & $-\frac{7}{48}$ & $-\frac{49}{64}$ & $\frac{1}{4}$ \\ \hline
$(4;5)$ & $\frac{1}{64}$ & $-\frac{7}{48}$ & $\frac{7}{8}$ & $-\frac{1}{4}$ \\ \hline
$(5;5)$ & $\frac{1}{12}$ & $0$  & $-\frac{63}{64}$ & $\frac{1}{4}$ \\ \hline
$(6;5)$ & $-\frac{1}{24}$ & $0$  & $\frac{49}{64}$ & $-\frac{1}{4}$ \\ \hline
$(7;2)$ & $0$ & $0$  & $-\frac{3}{8}$ & $0$ \\ \hline
$(7;3)$ & $-\frac{1}{12}$ & $0$  & $-\frac{33}{16}$ & $1$ \\ \hline
$(7;5)$ & $-\frac{1}{24}$ & $0$  & $\frac{1}{8}$ & $0$ \\ \hline
$(8;4)$ & $0$ & $0$  & $\frac{7}{8}$ & $-\frac{1}{4}$ \\ \hline
\end{tabular}
\end{center}

{\small \begin{center}
\begin{tabular}{|c|c|c|c|c|c|c|c|c|}
\hline
Weight 4 QPLI & $\pi^4$ & $\pi^2 \log^2(2)$  & $\log^4(2)$ & $\zeta(3) \log(2)$ & $\text{Li}_4(\frac{1}{2})$ & $\pi P_3$ & $\pi C \log(2)$ & $C^2$ \\ \hline
$(1;1;5)$ & $-\frac{511}{46080}$ & $-\frac{7}{384}$ & $\frac{5}{128}$ & $0$ & $\frac{15}{16}$ & $0$ & $-\frac{1}{8}$ & $\frac{1}{2}$ \\ \hline
$(1;2;5)$ & $-\frac{41}{7680}$ & $0$ & $0$ & $0$ & $0$ & $0$ & $0$ & $\frac{1}{2}$ \\ \hline
$(1;3;5)$ & $-\frac{269}{46080}$ & $\frac{1}{384}$ & $\frac{7}{384}$ & $\frac{7}{8}$ & $\frac{7}{16}$ & $0$ & $-\frac{1}{8}$ & $0$ \\ \hline
$(1;4;2)$ & $-\frac{23}{2304}$ & $-\frac{5}{96}$ & $\frac{5}{96}$ & $\frac{35}{32}$ & $\frac{5}{4}$ & $0$ & $0$ & $0$ \\ \hline
$(1;4;3)$ & $-\frac{35}{4608}$ & $\frac{13}{192}$ & $\frac{11}{192}$ & $\frac{7}{32}$ & $\frac{11}{8}$ & $0$ & $\frac{1}{4}$ & $-1$ \\ \hline
$(1;4;5)$ & $-\frac{877}{23040}$ & $-\frac{1}{24}$ & $\frac{1}{48}$ & $\frac{7}{16}$ & $\frac{1}{2}$ & $1$ & $-\frac{1}{2}$ & $\frac{1}{2}$ \\ \hline
$(1;5;4)$ & $\frac{89}{46080}$ & $-\frac{1}{96}$ & $\frac{1}{96}$ & $\frac{7}{32}$ & $\frac{1}{4}$ & $0$ & $0$ & $-\frac{1}{4}$ \\
\hline
\end{tabular}
\end{center}

\begin{center}
\begin{tabular}{|c|c|c|c|c|c|c|c|c|}
\hline
Weight 4 QPLI & $\pi^4$ & $\pi^2 \log^2(2)$  & $\log^4(2)$ & $\zeta(3) \log(2)$ & $\text{Li}_4(\frac{1}{2})$ & $\pi P_3$ & $\pi C \log(2)$ & $C^2$ \\ \hline
$(5;4;3)$ & $-\frac{79}{23040}$ & $\frac{13}{192}$ & $\frac{11}{192}$ & $-\frac{13}{32}$ & $\frac{11}{8}$ & $0$ & $\frac{1}{4}$ & $-1$ \\ \hline
$(5;4;5)$ & $\frac{529}{23040}$ & $\frac{1}{32}$ & $\frac{7}{96}$ & $\frac{35}{64}$ & $\frac{7}{4}$ & $-1$ & $\frac{3}{4}$ & $-\frac{3}{2}$ \\ \hline
$(5;5;4)$ & $\frac{391}{46080}$ & $\frac{5}{192}$ & $-\frac{1}{96}$ & $-\frac{7}{32}$ & $-\frac{1}{4}$ & $-\frac{1}{4}$ & $\frac{1}{8}$ & $\frac{1}{4}$ \\ \hline
$(5;5)$ & $\frac{269}{23040}$ & $\frac{11}{192}$ & $-\frac{11}{192}$ & $-\frac{45}{64}$ & $-\frac{11}{8}$ & $0$ & $0$ & $0$ \\ \hline
$(6;1;5)$ & $\frac{101}{46080}$ & $\frac{5}{384}$ & $-\frac{3}{128}$ & $-\frac{7}{64}$ & $-\frac{9}{16}$ & $0$ & $-\frac{1}{8}$ & $\frac{1}{2}$ \\ \hline
$(6;2;5)$ & $-\frac{97}{7680}$ & $-\frac{1}{48}$ & $\frac{1}{48}$ & $\frac{7}{8}$ & $\frac{1}{2}$ & $0$ & $0$ & $\frac{1}{2}$ \\ \hline
$(6;3;5)$ & $-\frac{349}{9216}$ & $-\frac{59}{384}$ & $\frac{55}{384}$ & $\frac{217}{64}$ & $\frac{55}{16}$ & $0$ & $-\frac{1}{8}$ & $0$ \\
\hline
\end{tabular}
\end{center}

\begin{center}
\begin{tabular}{|c|c|c|c|c|c|c|c|c|}
\hline
Weight 4 QPLI & $\pi^3 \log(2)$ & $\pi \log^3(2)$  & $\pi^2 C$ & $C \log^2(2)$ & $\pi \zeta(3)$ & $\log(2) P_3$ & $ P_4$ & $\Psi_3$ \\ \hline
$(1;1;4)$ & $\frac{1}{48}$ & $\frac{1}{96}$ & $-\frac{7}{32}$ & $\frac{1}{8}$ & $0$ & $0$ & $-1$ & $\frac{1}{768}$ \\ \hline
$(1;2;4)$ & $0$ & $0$ & $-\frac{17}{48}$ & $0$ & $0$ & $0$ & $0$ & $\frac{1}{512}$ \\ \hline
$(1;3;4)$ & $\frac{1}{48}$ & $\frac{17}{96}$ & $-\frac{5}{32}$ & $\frac{9}{8}$ & $0$ & $-4$ & $7$ & $-\frac{1}{384}$ \\ \hline
$(1;4;4)$ & $\frac{1}{12}$ & $\frac{1}{24}$ & $-\frac{13}{48}$ & $\frac{1}{2}$ & $-\frac{7}{32}$ & $0$ & $-4$ & $\frac{1}{256}$ \\ \hline
$(1;5;2)$ & $0$ & $0$ & $\frac{1}{6}$ & $0$ & $\frac{35}{128}$ & $0$ & $0$ & $-\frac{1}{768}$ \\ \hline
$(1;5;3)$ & $\frac{1}{48}$ & $-\frac{1}{96}$ & $\frac{5}{96}$ & $-\frac{1}{8}$ & $\frac{7}{128}$ & $0$ & $1$ & $-\frac{1}{768}$ \\ \hline
$(1;5;5)$ & $-\frac{1}{48}$ & $\frac{1}{48}$ & $\frac{5}{96}$ & $0$ & $\frac{7}{64}$ & $-1$ & $4$ & $-\frac{1}{384}$ \\
\hline
\end{tabular}
\end{center}

\begin{center}
\begin{tabular}{|c|c|c|c|c|c|c|c|c|}
\hline
Weight 4 QPLI & $\pi^3 \log(2)$ & $\pi \log^3(2)$  & $\pi^2 C$ & $C \log^2(2)$ & $\pi \zeta(3)$ & $\log(2) P_3$ & $ P_4$ & $\Psi_3$ \\ \hline
$(5;5;2)$ & $0$ & $\frac{1}{24}$ & $\frac{1}{16}$ & $\frac{1}{4}$ & $-\frac{7}{16}$ & $-1$ & $2$ & $0$ \\ \hline
$(5;5;3)$ & $\frac{1}{48}$ & $\frac{1}{96}$ & $-\frac{5}{96}$ & $\frac{1}{8}$ & $-\frac{27}{128}$ & $0$ & $-1$ & $\frac{1}{768}$ \\ \hline
$(5;5;5)$ & $-\frac{1}{48}$ & $\frac{1}{32}$ & $\frac{1}{16}$ & $\frac{1}{8}$ & $-\frac{49}{256}$ & $-1$ & $3$ & $-\frac{1}{768}$ \\ \hline
$(5;4)$ & $0$ & $0$ & $-\frac{1}{6}$ & $0$ & $-\frac{3}{128}$ & $0$ & $0$ & $\frac{1}{768}$ \\ \hline
$(6;1;4)$ & $\frac{7}{192}$ & $-\frac{3}{32}$ & $\frac{1}{32}$ & $-\frac{3}{8}$ & $0$ & $3$ & $-9$ & $\frac{7}{1536}$ \\ \hline 
$(6;2;4)$ & $0$ & $-\frac{1}{6}$ & $-\frac{1}{48}$ & $-1$ & $0$ & $4$ & $-8$ & $\frac{7}{1536}$ \\ \hline
$(6;3;4)$ & $\frac{7}{192}$ & $-\frac{3}{32}$ & $\frac{1}{96}$ & $-\frac{3}{8}$ & $0$ & $3$ & $-9$ & $\frac{7}{1536}$ \\
\hline
\end{tabular}
\end{center}}

\clearpage

\subsection{QES}
Table of 93 QESs with weight $\leq4$, 28 non-alternating QESs with weight 5. 10 results established only by using level 4 MZVs \cite{ref24} are dyed \textcolor{blue}{blue}.\\

$$\text{Weight 2 QES}: QES(31;-11)=\frac{5}{48} \pi^2-\frac{1}{4} \log^2(2)$$
\begin{center}
\begin{tabular}{|c|c|c|}
\hline
Weight 2 QES & $\pi \log(2)$ & $C$ \\ \hline
$(11;-12)$ & $\frac{1}{2}$ & $-1$ \\ \hline
$(21;-12)$ & $-\frac{1}{4}$ & $1$ \\ \hline
$(31;-12)$ & $\frac{1}{8}$ & $0$ \\ \hline
\end{tabular}
\end{center}

\begin{center}
\begin{tabular}{|c|c|c|c|c|}
\hline
Weight 3 QES & $\pi^2 \log(2)$ & $\log^3 (2)$  & $\zeta(3)$ & $\pi C$ \\ \hline
$(11,31;-11)$ & $-\frac{5}{48}$ & $\frac{1}{6}$ & $-\frac{1}{16}$ & $\frac{1}{2}$ \\ \hline
$(21,31;-11)$ & $\frac{7}{48}$ & $-\frac{1}{6}$ & $\frac{25}{32}$ & $0$ \\ \hline
$(31,31;-11)$ & $-\frac{1}{12}$ & $\frac{1}{12}$ & $\frac{3}{16}$ & $\frac{1}{2}$ \\ \hline
$(32;-11)$ & $-\frac{1}{48}$ & $0$ & $2$ & $-\frac{1}{2}$ \\ \hline
$(31;21)$ & $0$ & $0$ & $\frac{11}{4}$ & $0$ \\ \hline
$(31;-21)$ & $0$ & $0$ & $-\frac{23}{16}$ & $1$ \\ \hline
$(11;22)$ & $-\frac{1}{4}$ & $0$ & $\frac{7}{4}$ & $0$ \\ \hline
$(21;22)$ & $\frac{1}{8}$ & $0$ & $-\frac{7}{4}$ & $\frac{1}{2}$ \\ \hline
$(31;22)$ & $0$ & $0$ & $\frac{7}{16}$ & $0$ \\ \hline
\end{tabular}
\end{center}

\begin{center}
\begin{tabular}{|c|c|c|c|c|}
\hline
Weight 3 QES & $\pi^3$ & $\pi \log^2 (2)$  & $C \log(2)$ & $P_3$ \\ \hline
$(11,11;-12)$ & $-\frac{1}{4}$ & $-\frac{3}{2}$ & $0$ & $8$ \\ \hline
$(11,21;-12)$ & $\frac{7}{32}$ & $\frac{7}{8}$ & $0$ & $-6$ \\ \hline
$(11,31;-12)$ & $-\frac{3}{64}$ & $-\frac{3}{8}$ & $-\frac{1}{2}$ & $2$ \\ \hline
$(12;-12)$ & $-\frac{5}{48}$ & $-\frac{1}{4}$ & $-2$ & $4$ \\ \hline
$(21,21;-12)$ & $-\frac{1}{8}$ & $-\frac{1}{2}$ & $0$ & $4$ \\ \hline
$(21,31;-12)$ & $\frac{7}{64}$ & $\frac{5}{16}$ & $1$ & $-3$ \\ \hline
$(22;-12)$ & $\frac{11}{48}$ & $\frac{1}{2}$ & $4$ & $-8$ \\ \hline
$(31,31;-12)$ & $-\frac{5}{96}$ & $-\frac{3}{16}$ & $-\frac{1}{2}$ & $2$ \\ \hline
$(32;-12)$ & $-\frac{11}{192}$ & $-\frac{1}{8}$ & $-\frac{1}{2}$ & $2$ \\ \hline
$(11;-22)$ & $-\frac{3}{32}$ & $-\frac{1}{8}$ & $1$ & $2$ \\ \hline
$(21;-22)$ & $-\frac{1}{32}$ & $-\frac{1}{8}$ & $-2$ & $2$ \\ \hline
$(31;-22)$ & $\frac{1}{32}$ & $\frac{1}{16}$ & $\frac{1}{2}$ & $-1$ \\ \hline
\end{tabular}
\end{center}

\begin{center}
\begin{tabular}{|c|c|c|c|c|c|c|c|c|}
\hline
Weight 4 QES & $\pi^4$ & $\pi^2 \log^2(2)$  & $\log^4(2)$ & $\zeta(3) \log(2)$ & $\text{Li}_4(\frac{1}{2})$ & $\pi P_3$ & $\pi C \log(2)$ & $C^2$ \\ \hline
$(11,11,31;-11)$ & $-\frac{67}{576}$ & $-\frac{1}{8}$ & $-\frac{1}{8}$ & $-\frac{23}{16}$ & $0$ & $4$ & $-2$ & $2$ \\ \hline
$\textcolor{blue}{(11,21,31;-11)}$ & $\frac{29}{1152}$ & $-\frac{7}{96}$ & $\frac{3}{32}$ & $\frac{21}{16}$ & $-\frac{3}{4}$ & $0$ & $\frac{1}{2}$ & $-1$ \\ \hline
$(11,31,31;-11)$ & $-\frac{191}{2304}$ & $-\frac{3}{32}$ & $-\frac{1}{16}$ & $-\frac{51}{32}$ & $0$ & $3$ & $-\frac{3}{2}$ & $2$ \\ \hline
$(11,32;-11)$ & $-\frac{1}{120}$ & $-\frac{11}{96}$ & $\frac{3}{16}$ & $\frac{63}{32}$ & $\frac{9}{2}$ & $-1$ & $1$ & $0$ \\ \hline
$(12,31;-11)$ & $\frac{5}{48}$ & $\frac{13}{48}$ & $-\frac{1}{8}$ & $-\frac{21}{16}$ & $-3$ & $-2$ & $0$ & $0$ \\ \hline
$\textcolor{blue}{(21,21,31;-11)}$ & $-\frac{23}{384}$ & $\frac{1}{12}$ & $-\frac{1}{8}$ & $-\frac{1}{4}$ & $0$ & $2$ & $0$ & $0$ \\ \hline
$(21,31,31;-11)$ & $\frac{133}{5760}$ & $-\frac{5}{48}$ & $\frac{1}{12}$ & $\frac{35}{16}$ & $\frac{1}{2}$ & $0$ & $\frac{1}{2}$ & $-1$ \\
\hline
\end{tabular}
\end{center}

\begin{center}
\begin{tabular}{|c|c|c|c|c|c|c|c|c|}
\hline
Weight 4 QES & $\pi^3 \log(2)$ & $\pi \log^3(2)$  & $\pi^2 C$ & $C \log^2(2)$ & $\pi \zeta(3)$ & $\log(2) P_3$ & $P_4$ & $\Psi_3$ \\ \hline
$(11,11,11;-12)$ & $\frac{3}{2}$ & $3$ & $\frac{1}{12}$ & $0$ & $0$ & $0$ & $-96$ & $\frac{37}{768}$ \\ \hline
$\textcolor{blue}{(11,11,21;-12)}$ & $-\frac{9}{8}$ & $-\frac{5}{3}$ & $\frac{1}{6}$ & $0$ & $0$ & $0$ & $64$ & $-\frac{1}{32}$ \\ \hline
$(11,11,31;-12)$ & $\frac{5}{16}$ & $\frac{2}{3}$ & $0$ & $\frac{1}{2}$ & $\frac{7}{32}$ & $2$ & $-28$ & $\frac{11}{768}$ \\ \hline
$(11,12;-12)$ & $\frac{5}{24}$ & $0$ & $-\frac{1}{12}$ & $1$ & $0$ & $4$ & $-24$ & $\frac{5}{384}$ \\ \hline
$\textcolor{blue}{(11,21,21;-12)}$ & $\frac{11}{16}$ & $\frac{11}{12}$ & $\frac{1}{12}$ & $0$ & $0$ & $0$ & $-40$ & $\frac{5}{256}$ \\ \hline
$\textcolor{blue}{(11,21,31;-12)}$ & $-\frac{3}{8}$ & $-\frac{23}{48}$ & $\frac{1}{8}$ & $-\frac{5}{4}$ & $0$ & $0$ & $22$ & $-\frac{1}{96}$ \\ \hline
$(11,22;-12)$ & $-\frac{11}{24}$ & $-\frac{1}{3}$ & $\frac{1}{6}$ & $-4$ & $0$ & $0$ & $32$ & $-\frac{13}{768}$ \\
\hline
\end{tabular}
\end{center}

{\footnotesize\begin{center}
\begin{tabular}{|c|c|c|c|c|c|c|c|c|}
\hline
Weight 5 QES & $\pi^4 \log(2)$ & $\pi^2 \log^3(2)$ & $\log^5(2)$ & $\pi^2 \zeta(3)$ & $\zeta(3) \log^2(2)$ & $\zeta(5)$ & $\text{Li}_4(\frac{1}{2}) \log(2)$ & $\text{Li}_5(\frac{1}{2})$ \\ \hline
$(31;41)$ & $0$ & $0$ & $0$ & $-\frac{2}{3}$ & $0$ & $\frac{37}{4}$ & $0$ & $0$ \\ \hline
$(11;42)$ & $-\frac{1}{48}$ & $0$ & $0$ & $-\frac{7}{32}$ & $0$ & $\frac{31}{8}$ & $0$ & $0$ \\ \hline
$(31;42)$ & $0$ & $0$ & $0$ & $-\frac{1}{8}$ & $0$ & $\frac{93}{64}$ & $0$ & $0$ \\ \hline
$(11,31;31)$ & $0$ & $\frac{4}{9}$ & $-\frac{8}{15}$ & $-\frac{1}{12}$ & $-7$ & $\frac{307}{16}$ & $-16$ & $-16$ \\ \hline
$(31,31;31)$ & $0$ & $\frac{4}{9}$ & $-\frac{8}{15}$ & $\frac{1}{4}$ & $-7$ & $\frac{131}{8}$ & $-16$ & $-16$ \\ \hline
$(32;31)$ & $0$ & $0$ & $0$ & $\frac{19}{12}$ & $0$ & $-\frac{133}{8}$ & $0$ & $0$ \\ \hline
$(11,11;32)$ & $-\frac{1}{16}$ & $0$ & $0$ & $-\frac{7}{48}$ & $\frac{7}{2}$ & $\frac{31}{8}$ & $0$ & $0$ \\
\hline
\end{tabular}
\end{center}}

\clearpage

\subsection{GES}
Table of 39 GESs with weight $\leq5$. We abbreviate $GES(a, b, c; d)$ as $GES(abc;d)$ and so on.\\

$$\text{Weight 1 GES}:\ \  GES(1;0)=2\log(2)$$

\begin{center}
\begin{tabular}{|c|c|c|}
\hline
Weight 2 GES & $\pi^2 $ & $\log^2(2)$ \\ \hline
$(2;0)$ & $\frac{1}{6}$ & $-1$ \\ \hline
$(11;0)$ & $\frac{1}{6}$ & $1$ \\ \hline
$(1;1)$ & $\frac{1}{12}$ & $0$ \\ \hline
\end{tabular}
\end{center}

\begin{center}
\begin{tabular}{|c|c|c|c|}
\hline
Weight 3 GES & $\pi^2 \log(2)$ & $\log^3(2)$ & $\zeta(3)$ \\ \hline
$(3;0)$ & $-\frac{1}{6}$ & $\frac{1}{3}$ & $\frac{7}{4}$ \\ \hline
$(12;0)$ & $0$ & $-\frac{1}{3}$ & $\frac{3}{2}$ \\ \hline
$(111;0)$ & $\frac{1}{3}$ & $\frac{1}{3}$ & $1$ \\ \hline
$(2;1)$ & $0$ & $0$ & $\frac{5}{8}$ \\ \hline
$(11;1)$ & $0$ & $0$ & $\frac{7}{8}$ \\ \hline
$(1;2)$ & $-\frac{1}{12}$ & $0$ & $1$ \\ \hline
\end{tabular}
\end{center}

\begin{center}
\begin{tabular}{|c|c|c|c|c|c|}
\hline
Weight 4 GES & $\pi^4$ & $\pi^2 \log^2(2)$ & $\log^4(2)$ & $\zeta(3) \log(2)$ & $\text{Li}_4(\frac{1}{2})$ \\ \hline
$(4;0)$ & $0$ & $0$ & $0$ & $0$ & $2$ \\ \hline
$(13;0)$ & $-\frac{1}{240}$ & $-\frac{1}{12}$ & $\frac{1}{6}$ & $\frac{3}{2}$ & $2$ \\ \hline
$(22;0)$ & $\frac{1}{360}$ & $-\frac{1}{6}$ & $\frac{1}{6}$ & $1$ & $2$ \\ \hline
$(112;0)$ & $\frac{7}{180}$ & $0$ & $-\frac{1}{6}$ & $0$ & $-2$ \\ \hline
$(1111;0)$ & $\frac{31}{120}$ & $\frac{7}{6}$ & $-\frac{5}{6}$ & $-15$ & $-22$ \\ \hline
$(3;1)$ & $-\frac{1}{288}$ & $-\frac{1}{24}$ & $\frac{1}{24}$ & $\frac{7}{8}$ & $1$ \\ \hline
$(12;1)$ & $-\frac{1}{720}$ & $-\frac{1}{24}$ & $\frac{1}{24}$ & $\frac{7}{8}$ & $1$ \\ \hline
$(111;1)$ & $\frac{5}{72}$ & $\frac{5}{24}$ & $-\frac{5}{24}$ & $-\frac{35}{8}$ & $-5$ \\ \hline
$(2;2)$ & $\frac{1}{1440}$ & $-\frac{1}{24}$ & $\frac{1}{24}$ & $\frac{1}{4}$ & $1$ \\ \hline
$(11;2)$ & $\frac{37}{1440}$ & $\frac{1}{24}$ & $-\frac{1}{24}$ & $-\frac{7}{4}$ & $-1$ \\ \hline
$(1;3)$ & $\frac{1}{720}$ & $0$ & $\frac{1}{24}$ & $-\frac{1}{8}$ & $1$ \\ \hline
\end{tabular}
\end{center}

{\small
\begin{center}
\begin{tabular}{|c|c|c|c|c|c|c|c|c|}
\hline
Weight 5 GES & $\pi^4 \log(2)$ & $\pi^2 \log^3(2)$ & $\log^5(2)$ & $\pi^2 \zeta(3)$ & $\zeta(3) \log^2(2)$ & $\zeta(5)$ & $\text{Li}_4(\frac{1}{2}) \log(2)$ & $\text{Li}_5(\frac{1}{2})$ \\ \hline
$(5;0)$ & $0$ & $0$ & $0$ & $0$ & $0$ & $0$ & $0$ & $2$ \\ \hline
$(14;0)$ & $-\frac{1}{360}$ & $0$ & $-\frac{1}{60}$ & $-\frac{1}{48}$ & $\frac{1}{8}$ & $\frac{7}{4}$ & $0$ & $0$ \\ \hline
$(23;0)$ & $\frac{1}{180}$ & $\frac{1}{18}$ & $-\frac{1}{30}$ & $\frac{3}{16}$ & $-\frac{9}{8}$ & $-\frac{29}{16}$ & $0$ & $2$ \\ \hline
$(113;0)$ & $-\frac{19}{360}$ & $-\frac{1}{18}$ & $\frac{1}{15}$ & $-\frac{7}{48}$ & $\frac{19}{8}$ & $\frac{147}{16}$ & $0$ & $-6$ \\ \hline
$(122;0)$ & $\frac{1}{90}$ & $0$ & $\frac{1}{60}$ & $\frac{5}{16}$ & $-\frac{11}{8}$ & $-\frac{3}{2}$ & $0$ & $0$ \\ \hline
$(1112;0)$ & $-\frac{1}{360}$ & $-\frac{1}{9}$ & $\frac{1}{60}$ & $-\frac{5}{16}$ & $\frac{15}{8}$ & $\frac{77}{8}$ & $0$ & $-4$ \\ \hline
$(11111;0)$ & $\frac{8}{9}$ & $\frac{10}{9}$ & $-\frac{7}{12}$ & $-\frac{85}{48}$ & $-\frac{215}{8}$ & $-\frac{195}{4}$ & $0$ & $72$ \\
\hline
\end{tabular}
\end{center}}

\clearpage

\subsection{BES}
Table of 26 BESs with weight $\leq5$. We abbreviate $BES(a, b, c; d)$ as $BES(abc;d)$ and so on. Also we write $BES(0;k)=\sum _{n=1}^{\infty } \frac{\binom{2 n}{n}}{4^n n^k}$.\\

$$\text{Weight 1 BES}:\ \  BES(0;1)=2\log(2)$$

$$\text{Weight 2 BES}:\ \  BES(0;2)=\frac{\pi ^2}{6}-2 \log ^2(2),\  BES(1;1)=\frac{\pi ^2}{3}$$

\begin{center}
\begin{tabular}{|c|c|c|c|}
\hline
Weight 3 BES & $\pi^2 \log(2)$ & $\log^3(2)$ & $\zeta(3)$ \\ \hline
$(0;3)$ & $-\frac{1}{3}$ & $\frac{4}{3}$ & $2$ \\ \hline
$(1;2)$ & $-\frac{2}{3}$ & $0$ & $\frac{9}{2}$ \\ \hline
$(2;1)$ & $0$ & $0$ & $\frac{3}{2}$ \\ \hline
$(11;1)$ & $0$ & $0$ & $\frac{21}{2}$ \\ \hline
\end{tabular}
\end{center}

\begin{center}
\begin{tabular}{|c|c|c|c|c|c|}
\hline
Weight 4 BES & $\pi^4$ & $\pi^2 \log^2(2)$ & $\log^4(2)$ & $\zeta(3) \log(2)$ & $\text{Li}_4(\frac{1}{2})$ \\ \hline
$(0;4)$ & $\frac{1}{40}$ & $\frac{1}{3}$ & $-\frac{2}{3}$ & $-4$ & $0$ \\ \hline
$(1;3)$ & $-\frac{13}{360}$ & $\frac{1}{3}$ & $\frac{1}{3}$ & $-2$ & $8$ \\ \hline
$(2;2)$ & $\frac{1}{30}$ & $0$ & $0$ & $-3$ & $0$ \\ \hline
$(11;2)$ & $-\frac{7}{45}$ & $-\frac{4}{3}$ & $\frac{4}{3}$ & $7$ & $32$ \\ \hline
$(3;1)$ & $\frac{37}{360}$ & $\frac{1}{3}$ & $-\frac{1}{3}$ & $-7$ & $-8$ \\ \hline
$(12;1)$ & $\frac{49}{360}$ & $\frac{1}{3}$ & $-\frac{1}{3}$ & $-7$ & $-8$ \\ \hline
$(111;1)$ & $\frac{23}{72}$ & $-\frac{5}{3}$ & $\frac{5}{3}$ & $35$ & $40$ \\ \hline
\end{tabular}
\end{center}

{\small
\begin{center}
\begin{tabular}{|c|c|c|c|c|c|c|c|c|}
\hline
Weight 5 BES & $\pi^4 \log(2)$ & $\pi^2 \log^3(2)$ & $\log^5(2)$ & $\pi^2 \zeta(3)$ & $\zeta(3) \log^2(2)$ & $\zeta(5)$ & $\text{Li}_4(\frac{1}{2}) \log(2)$ & $\text{Li}_5(\frac{1}{2})$ \\ \hline
$(0;5)$ & $-\frac{1}{20}$ & $-\frac{2}{9}$ & $\frac{4}{15}$ & $-\frac{1}{3}$ & $4$ & $6$ & $0$ & $0$ \\ \hline
$(1;4)$ & $\frac{13}{180}$ & $-\frac{2}{9}$ & $-\frac{2}{15}$ & $-\frac{5}{6}$ & $2$ & $-\frac{23}{8}$ & $0$ & $16$ \\ \hline
$(2;3)$ & $-\frac{1}{15}$ & $\frac{4}{9}$ & $-\frac{8}{15}$ & $\frac{1}{3}$ & $-4$ & $\frac{123}{8}$ & $-16$ & $-16$ \\ \hline
$(11;3)$ & $\frac{14}{45}$ & $\frac{4}{9}$ & $0$ & $-\frac{2}{3}$ & $0$ & $-\frac{465}{8}$ & $16$ & $80$ \\ \hline
$(3;2)$ & $-\frac{37}{180}$ & $-\frac{2}{9}$ & $\frac{2}{15}$ & $-\frac{7}{12}$ & $7$ & $\frac{203}{8}$ & $0$ & $-16$ \\ \hline
$(12;2)$ & $-\frac{49}{180}$ & $\frac{2}{3}$ & $-\frac{14}{15}$ & $\frac{1}{12}$ & $-7$ & $\frac{217}{4}$ & $-32$ & $-48$ \\ \hline
$(111;2)$ & $-\frac{23}{36}$ & $-\frac{14}{9}$ & $\frac{38}{15}$ & $\frac{35}{12}$ & $7$ & $-\frac{217}{2}$ & $96$ & $176$ \\ \hline
$(4;1)$ & $0$ & $-\frac{4}{9}$ & $\frac{8}{15}$ & $-\frac{7}{12}$ & $7$ & $-\frac{63}{8}$ & $16$ & $16$ \\ \hline
$(13;1)$ & $0$ & $-\frac{4}{9}$ & $\frac{8}{15}$ & $-\frac{7}{4}$ & $7$ & $\frac{31}{4}$ & $16$ & $16$ \\ \hline
$(22;1)$ & $0$ & $\frac{4}{9}$ & $-\frac{8}{15}$ & $\frac{7}{12}$ & $-7$ & $\frac{93}{8}$ & $-16$ & $-16$ \\ \hline
$(112;1)$ & $0$ & $\frac{4}{3}$ & $-\frac{8}{5}$ & $-\frac{7}{12}$ & $-21$ & $\frac{589}{8}$ & $-48$ & $-48$ \\ \hline
$(1111;1)$ & $0$ & $-\frac{28}{9}$ & $\frac{56}{15}$ & $\frac{77}{4}$ & $49$ & $\frac{1829}{8}$ & $112$ & $112$ \\ \hline
\end{tabular}
\end{center}}

\clearpage

\subsection{IBES}
Table of 14 IBESs with weight $\leq5$. We abbreviate $IBES(a, b, c; d)$ as $IBES(abc;d)$ and so on. Also we write $IBES(0;k)=\sum _{n=1}^{\infty } \frac{4^n }{n^k\binom{2 n}{n}}$.\\

$$\text{Weight 2 IBES}:\ \  IBES(0;2)=\frac{\pi ^2}{2}$$

$$\text{Weight 3 IBES}:\ \  IBES(0;3)=\pi ^2 \log (2)-\frac{7 \zeta (3)}{2},\  IBES(1;2)=\frac{7 \zeta (3)}{2}+\pi ^2 \log (2)$$
$$$$
\begin{center}
\begin{tabular}{|c|c|c|c|c|c|}
\hline
Weight 4 IBES & $\pi^4$ & $\pi^2 \log^2(2)$ & $\log^4(2)$ & $\zeta(3) \log(2)$ & $\text{Li}_4(\frac{1}{2})$ \\ \hline
$(0;4)$ & $-\frac{19}{360}$ & $\frac{2}{3}$ & $\frac{1}{3}$ & $0$ & $8$ \\ \hline
$(1;3)$ & $\frac{1}{90}$ & $\frac{4}{3}$ & $-\frac{1}{3}$ & $0$ & $-8$ \\ \hline
$(2;2)$ & $-\frac{1}{90}$ & $\frac{2}{3}$ & $\frac{1}{3}$ & $0$ & $8$ \\ \hline
$(11;2)$ & $\frac{9}{20}$ & $2$ & $-1$ & $0$ & $-24$ \\ \hline
\end{tabular}
\end{center}
$$$$
{\small
\begin{center}
\begin{tabular}{|c|c|c|c|c|c|c|c|c|}
\hline
Weight 5 IBES & $\pi^4 \log(2)$ & $\pi^2 \log^3(2)$ & $\log^5(2)$ & $\pi^2 \zeta(3)$ & $\zeta(3) \log^2(2)$ & $\zeta(5)$ & $\text{Li}_4(\frac{1}{2}) \log(2)$ & $\text{Li}_5(\frac{1}{2})$ \\ \hline
$(0;5)$ & $-\frac{19}{180}$ & $\frac{4}{9}$ & $\frac{2}{15}$ & $1$ & $0$ & $\frac{31}{8}$ & $0$ & $-16$ \\ \hline
$(1;4)$ & $\frac{1}{45}$ & $\frac{8}{9}$ & $-\frac{2}{15}$ & $\frac{1}{2}$ & $0$ & $-\frac{31}{2}$ & $0$ & $16$ \\ \hline
$(2;3)$ & $-\frac{1}{45}$ & $\frac{4}{9}$ & $\frac{2}{15}$ & $\frac{1}{4}$ & $0$ & $\frac{31}{4}$ & $0$ & $-16$ \\ \hline
$(11;3)$ & $\frac{9}{10}$ & $\frac{4}{3}$ & $-\frac{2}{5}$ & $-\frac{7}{4}$ & $0$ & $-62$ & $0$ & $48$ \\ \hline
$(3;2)$ & $-\frac{19}{180}$ & $\frac{4}{9}$ & $\frac{2}{15}$ & $-\frac{3}{4}$ & $0$ & $\frac{217}{8}$ & $0$ & $-16$ \\ \hline
$(12;2)$ & $\frac{19}{180}$ & $\frac{8}{9}$ & $-\frac{2}{15}$ & $\frac{3}{2}$ & $0$ & $-\frac{155}{8}$ & $0$ & $16$ \\ \hline
$(111;2)$ & $\frac{91}{36}$ & $\frac{16}{9}$ & $-\frac{2}{3}$ & $3$ & $0$ & $-\frac{155}{8}$ & $0$ & $80$ \\ \hline
\end{tabular}
\end{center}}
$$$$

\subsection{QBES, IQBES}

Tabulated in subsection 8-3-4 already.\\

\clearpage

\section{Miscellaneous}
We tabulate some typical results without proof. All can be derived using methods above.\\

\subsection{LSI}
$$\int_0^{\frac{\pi }{2}} x^2 \cot (x) \log (2 \sin (x)) \log (2 \cos (x)) \, dx$$ $$=\frac{\pi ^2 \zeta (3)}{12}-\frac{31 \zeta (5)}{64}-\frac{7}{8} \zeta (3) \log ^2(2)+\frac{1}{12} \pi ^2 \log ^3(2)-\frac{1}{192} \pi ^4 \log (2)$$

$$\int_0^{\frac{\pi }{2}} x \log ^2(2 \sin (x)) \log (2 \cos (x)) \, dx$$ $$=\text{Li}_5\left(\frac{1}{2}\right)+\text{Li}_4\left(\frac{1}{2}\right) \log (2)\\-\frac{155 \zeta (5)}{128}-\frac{\pi ^2 \zeta (3)}{192}+\frac{7}{16} \zeta (3) \log ^2(2)+\frac{\log ^5(2)}{30}-\frac{1}{36} \pi ^2 \log ^3(2)$$

$$\int_0^{\frac{\pi }{2}} x^2 \log (2 \sin (x)) \log ^2(2 \cos (x)) \, dx$$ $$=-\pi  \text{Li}_5\left(\frac{1}{2}\right)-\pi  \text{Li}_4\left(\frac{1}{2}\right) \log (2)\\+\frac{3 \pi ^3 \zeta (3)}{64}+\frac{121 \pi  \zeta (5)}{128}-\frac{7}{16} \pi  \zeta (3) \log ^2(2)-\frac{1}{30} \pi  \log ^5(2)+\frac{1}{36} \pi ^3 \log ^3(2)$$

$$\int_0^{\frac{\pi }{2}} x \log ^3(2 \sin (x)) \log (2 \cos (x)) \, dx$$ $$=\frac{9}{4} ES(1;-5)-\frac{1}{8} \pi ^2 \text{Li}_4\left(\frac{1}{2}\right)+6 \text{Li}_6\left(\frac{1}{2}\right)+3 \text{Li}_4\left(\frac{1}{2}\right) \log ^2(2)+6 \text{Li}_5\left(\frac{1}{2}\right) \log (2)+\frac{3 \zeta (3)^2}{2}$$ $$-\frac{1579 \pi ^6}{161280}+\frac{7}{8} \zeta (3) \log ^3(2)+\frac{\log ^6(2)}{12}\\-\frac{7}{64} \pi ^2 \zeta (3) \log (2)-\frac{13}{192} \pi ^2 \log ^4(2)+\frac{1}{192} \pi ^4 \log ^2(2)$$

$$\int_0^{\frac{\pi }{2}} x^4 \log ^2(2 \sin (x)) \, dx$$ $$=\frac{3}{2} \pi  ES(1;-5)+\frac{1}{2} \pi ^3 \text{Li}_4\left(\frac{1}{2}\right)\\+\frac{7}{16} \pi ^3 \zeta (3) \log (2)-\frac{269 \pi ^7}{40320}+\frac{1}{48} \pi ^3 \log ^4(2)-\frac{1}{48} \pi ^5 \log ^2(2)$$\\

\subsection{PLSI}
$$\int_0^{\frac{\pi }{2}} x \text{Ti}_2(\tan (x)) \log (2 \sin (x)) \, dx$$ $$=\frac{7}{4} \pi  \text{Li}_4\left(\frac{1}{2}\right)+\frac{49}{32} \pi  \zeta (3) \log (2)-\frac{163 \pi ^5}{11520}+\frac{7}{96} \pi  \log ^4(2)-\frac{7}{96}  \pi ^3 \log ^2(2)$$

$$\int_0^{\frac{\pi }{2}} \text{Cl}_2(2 x) \log (2 \sin (x)) \log (2 \cos (x)) \, dx$$ $$=2 \text{Li}_5\left(\frac{1}{2}\right)+2 \text{Li}_4\left(\frac{1}{2}\right) \log (2)\\-\frac{279 \zeta (5)}{128}+\frac{7}{8} \zeta (3) \log ^2(2)+\frac{\log ^5(2)}{15}-\frac{1}{18} \pi ^2 \log ^3(2)$$

$$\int_0^{\frac{\pi }{2}} \text{Cl}_2(4 x) \log ^2(2 \cos (x)) \, dx$$ $$=16 \text{Li}_5\left(\frac{1}{2}\right)+16 \text{Li}_4\left(\frac{1}{2}\right) \log (2)\\-\frac{465 \zeta (5)}{32}-\frac{7 \pi ^2 \zeta (3)}{48}+7 \zeta (3) \log ^2(2)+\frac{8 \log ^5(2)}{15}-\frac{4}{9} \pi ^2 \log ^3(2)$$

$$\int_0^{\frac{\pi }{2}} x \text{Cl}_2(4 x) \cot (x) \log (2 \sin (2 x)) \, dx$$ $$=16 \text{Li}_4\left(\frac{1}{2}\right) \log(2)-\frac{23 \pi ^2 \zeta (3)}{48}\\+\frac{31 \zeta (5)}{4}+14 \zeta (3) \log ^2(2)+\frac{2 \log ^5(2)}{3}-\frac{2}{3} \pi ^2 \log ^3(2)-\frac{151}{720} \pi ^4 \log (2)$$\\

\subsection{ITF}
$$\int_0^{\infty } \frac{\log ^2\left(x^2+1\right) \tan ^{-1}(x)^2}{x^4} \, dx$$ $$=-\frac{\pi  \zeta (3)}{6}+\frac{\pi ^3}{18}-\frac{4}{9} \pi  \log ^3(2)+\frac{4}{3} \pi  \log ^2(2)-\frac{2}{9} \pi ^3 \log (2)+\frac{4}{3} \pi  \log (2)$$

$$\int_0^{\infty } \frac{\log ^2\left(x^2+1\right) \tan ^{-1}(x)^3}{x^2} \, dx$$ $$=24 \text{Li}_5\left(\frac{1}{2}\right)+\frac{3 \pi ^2 \zeta (3)}{8}\\-\frac{93 \zeta (5)}{2}-\frac{1}{5} \log ^5(2)+\frac{4}{3} \pi ^2 \log ^3(2)+\frac{47}{60} \pi ^4 \log (2)$$

$$\int_0^{\infty } \frac{\log ^2\left(x^2+1\right) \tan ^{-1}(x)^2}{x(x^2+1)} \, dx$$ $$=8 \text{Li}_5\left(\frac{1}{2}\right)-\frac{3 \pi ^2 \zeta (3)}{8}\\-\frac{217 \zeta (5)}{16}-\frac{1}{15} \log ^5(2)+\frac{4}{9} \pi ^2 \log ^3(2)+\frac{79}{360} \pi ^4 \log (2)$$\\

\subsection{QLI $(W=5)$}
{\small$$\int_0^1 \frac{x \log ^2\left(x^2+1\right) \tan ^{-1}(x)^2}{x^2+1} \, dx$$ $$=\frac{1}{4} \pi  C \log ^2(2)+2 \pi  \Im(\text{Li}_4(1+i))+\pi  \log (2) \Im(\text{Li}_3(1+i))+8 \Re\left(\text{Li}_5\left(\frac{1}{2}+\frac{i}{2}\right)\right)-8 \text{Li}_5\left(\frac{1}{2}\right)-\frac{5}{4} \text{Li}_4\left(\frac{1}{2}\right) \log (2)$$ $$\ \ \ -\frac{\pi ^2 \zeta (3)}{8}+\frac{465 \zeta (5)}{512}+\frac{35}{64} \zeta (3) \log ^2(2)+\frac{\log ^5(2)}{80}-\frac{119}{576} \pi ^2 \log ^3(2)-\frac{2473 \pi ^4 \log (2)}{23040}-\frac{\pi  \left(\psi ^{(3)}\left(\frac{1}{4}\right)-\psi ^{(3)}\left(\frac{3}{4}\right)\right)}{3072}$$

$$\int_0^1 \frac{\log ^2\left(x^2+1\right) \tan ^{-1}(x)^2}{x^2+1} \, dx$$ $$=-8 \sqrt{2} \, _6F_5\left(\frac{1}{2},\frac{1}{2},\frac{1}{2},\frac{1}{2},\frac{1}{2},\frac{1}{2};\frac{3}{2},\frac{3}{2},\frac{3}{2},\frac{3}{2},\frac{3}{2};\frac{1}{2}\right)-\frac{1}{8} \pi ^2 C \log (2)-\frac{1}{4} \pi ^2 \Im(\text{Li}_3(1+i))-8 \Im\left(\text{Li}_5\left(\frac{1}{2}+\frac{i}{2}\right)\right)$$ $$\ +\frac{5}{8} \pi  \text{Li}_4\left(\frac{1}{2}\right)+\frac{61}{64} \pi  \zeta (3) \log (2)+\frac{709 \pi ^5}{15360}+\frac{1}{24} \pi  \log ^4(2)+\frac{9}{128} \pi ^3 \log ^2(2)+\frac{\log (2) \left(\psi ^{(3)}\left(\frac{1}{4}\right)-\psi ^{(3)}\left(\frac{3}{4}\right)\right)}{1536}$$}\\

\subsection{QLSI}
$$\int_0^{\frac{\pi }{4}} x \log ^2(\sin (x)) \, dx$$ $$=\frac{1}{8} \pi  C \log (2)+\frac{1}{4} \pi  \Im(\text{Li}_3(1+i))+\frac{5 \text{Li}_4\left(\frac{1}{2}\right)}{16}-\frac{35}{128} \zeta (3) \log (2)-\frac{11 \pi ^4}{1152}+\frac{5 \log ^4(2)}{384}+\frac{1}{384} \pi ^2 \log ^2(2)$$

$$\int_0^{\frac{\pi }{4}} \log ^3(2 \cos (x)) \, dx$$ $$=3 \Im(\text{Li}_4(1+i))-\frac{3}{2} \log (2) \Im(\text{Li}_3(1+i))\\ +\frac{3}{8} C \log ^2(2)-\frac{3 \pi  \zeta (3)}{8}+\frac{1}{16} \pi  \log ^3(2)-\frac{\psi ^{(3)}\left(\frac{1}{4}\right)-\psi ^{(3)}\left(\frac{3}{4}\right)}{2048}$$

$$\int_0^{\frac{\pi }{4}} \log (2 \sin (x)) \log ^2(2 \cos (x)) \, dx$$ $$=\frac{1}{8} C \log ^2(2)+\Im(\text{Li}_4(1+i))+\frac{\pi  \zeta (3)}{16}\\-\frac{1}{2} \log (2) \Im(\text{Li}_3(1+i))+\frac{1}{48} \pi  \log ^3(2)-\frac{5 \left(\psi ^{(3)}\left(\frac{1}{4}\right)-\psi ^{(3)}\left(\frac{3}{4}\right)\right)}{6144}$$\\

\subsection{NLI}

{\small $$\int_0^1 \log ^3(1-x) \log ^2(x+1) \, dx$$ $$=24 \text{Li}_4\left(\frac{1}{2}\right)+24 \text{Li}_5\left(\frac{1}{2}\right)-2 \pi ^2 \zeta (3)+45 \zeta (3)+\frac{3 \zeta (5)}{4}+12 \zeta (3) \log ^2(2)-24 \zeta (3) \log (2)+\frac{\pi ^4}{30}+6 \pi ^2-120$$ $$+\frac{3 \log ^5(2)}{5}-3 \log ^4(2)-\frac{2}{3} \pi ^2 \log ^3(2)+16 \log ^3(2)+2 \pi ^2 \log ^2(2)-48 \log ^2(2)-\frac{1}{30} \pi ^4 \log (2)-6 \pi ^2 \log (2)+96 \log (2)$$

$$\int_0^1 \frac{\log ^2(1-x) \log ^2(x) \log ^2(x+1)}{x^2} \, dx$$ $$=-8 ES(1;-5)-\frac{8}{3} \pi ^2 \text{Li}_4\left(\frac{1}{2}\right)+8 \text{Li}_4\left(\frac{1}{2}\right)+24 \text{Li}_5\left(\frac{1}{2}\right)+16 \text{Li}_6\left(\frac{1}{2}\right)+8 \text{Li}_4\left(\frac{1}{2}\right) \log ^2(2)$$ $$+24 \text{Li}_4\left(\frac{1}{2}\right) \log (2)+16 \text{Li}_5\left(\frac{1}{2}\right) \log (2)+\frac{19 \pi ^2 \zeta (3)}{6}-\frac{227 \zeta (3)^2}{8}-62 \zeta (5)-7 \zeta (3) \log ^3(2)$$ $$-\frac{35}{2} \zeta (3) \log ^2(2)-9 \zeta (3) \log (2)-\frac{217}{2} \zeta (5) \log (2)+\frac{35}{6} \pi ^2 \zeta (3) \log (2)+\frac{\pi ^4}{90}+\frac{25 \pi ^6}{378}+\frac{2 \log ^6(2)}{9}$$ $$+\frac{4 \log ^5(2)}{5}-\frac{5 \log ^4(2)}{3}-\frac{5}{18} \pi ^2 \log ^4(2)+\frac{2}{3} \pi ^2 \log ^3(2)+\pi ^2 \log ^2(2)+\frac{13}{36} \pi ^4 \log ^2(2)+\frac{1}{6} \pi ^4 \log (2)$$

$$\int_0^1 \log ^2(1-x) \log ^2(x) \log ^2(x+1) \, dx$$ $$=-4 ES(1;-5)-\frac{16}{3} \pi ^2 \text{Li}_4\left(\frac{1}{2}\right)-32 \text{Li}_4\left(\frac{1}{2}\right)+24 \text{Li}_5\left(\frac{1}{2}\right)-16 \text{Li}_6\left(\frac{1}{2}\right)-8 \text{Li}_4\left(\frac{1}{2}\right) \log ^2(2)+24 \text{Li}_4\left(\frac{1}{2}\right) \log (2)$$ $$-16 \text{Li}_5\left(\frac{1}{2}\right) \log (2)+\frac{275 \zeta (3)^2}{8}+\frac{34 \pi ^2 \zeta (3)}{3}-237 \zeta (3)-\frac{341 \zeta (5)}{2}+7 \zeta (3) \log ^3(2)-\frac{91}{2} \zeta (3) \log ^2(2)- \frac{77}{6} \pi ^2 \zeta (3) \log (2)$$ $$+156 \zeta (3) \log (2)+\frac{217}{2} \zeta (5) \log (2)+\frac{73 \pi ^6}{3024}-40 \pi ^2-\frac{29 \pi ^4}{45}+720-\frac{2}{9}  \log ^6(2)+\frac{4 \log ^5(2)}{5}-\frac{1}{18} \pi ^2 \log ^4(2)+\frac{2 \log ^4(2)}{3}$$ $$+\frac{2}{3} \pi ^2 \log ^3(2)-24 \log ^3(2)-\frac{1}{36} \pi ^4 \log ^2(2)-8 \pi ^2 \log ^2(2)+144 \log ^2(2)+\frac{2}{3} \pi ^4 \log (2)+32 \pi ^2 \log (2)-480 \log (2)$$}

\subsection{NPLI}

{\footnotesize $$\int_0^1 \text{Li}_2\left(\frac{1-x}{2}\right) \text{Li}_2\left(\frac{x-1}{x+1}\right) \log (x+1) \, dx$$ $$=6 \text{Li}_4\left(\frac{1}{2}\right)-24 \text{Li}_5\left(\frac{1}{2}\right)-18 \text{Li}_4\left(\frac{1}{2}\right) \log (2)+\frac{\pi ^2 \zeta (3)}{12}-\frac{3 \zeta (3)}{2}+\frac{189 \zeta (5)}{8}-\frac{21}{4} \zeta (3) \log ^2(2)+\frac{23}{4} \zeta (3) \log (2)$$ $$-\frac{43 \pi ^4}{720}-\frac{\pi ^2}{4}-\frac{3 \log ^5(2)}{5}+\frac{\log ^4(2)}{3}+\log ^3(2)+\frac{17}{36} \pi ^2 \log ^3(2)+3 \log ^2(2)-\frac{3}{8} \pi ^2 \log ^2(2)+\frac{1}{3} \pi ^2 \log (2)-\frac{29}{360} \pi ^4 \log (2)$$

$$\int_0^1 \text{Li}_3\left(\frac{x+1}{2}\right) \log (x) \log (1-x) \, dx$$ $$=-7 \text{Li}_4\left(\frac{1}{2}\right)+11 \text{Li}_5\left(\frac{1}{2}\right)+4 \text{Li}_4\left(\frac{1}{2}\right) \log (2)+\frac{49 \pi ^2 \zeta (3)}{96}-10 \zeta (3)-\frac{463 \zeta (5)}{32}+\frac{21}{16} \zeta (3) \log ^2(2)-\frac{11}{4} \zeta (3) \log (2)+\frac{121 \pi ^4}{1440}$$ $$-\frac{13 \pi ^2}{6}+20+\frac{3 \log ^5(2)}{40}-\frac{5 \log ^4(2)}{24}-\frac{1}{24} \pi ^2 \log ^3(2)+\frac{\log ^3(2)}{3}+\frac{5}{24} \pi ^2 \log ^2(2)-3 \log ^2(2)+\frac{1}{288} \pi ^4 \log (2)+\frac{2}{3} \pi ^2 \log (2)+12 \log (2)$$}\\

\subsection{NQLI}

$$\int_0^1 \log (1-x) \log (x+1) \log \left(x^2+1\right) \tan ^{-1}(x) \, dx$$ $$=C^2-\frac{5 \pi ^2 C}{24}+\frac{\pi  C}{2}-4 C-\frac{1}{2} \pi  C \log (2)-C \log (2)+\frac{1}{2} \pi  \Im(\text{Li}_3(1+i))+10 \Im(\text{Li}_3(1+i))$$ $$-24 \Im(\text{Li}_4(1+i))+3 \log (2) \Im(\text{Li}_3(1+i))-\frac{13 \text{Li}_4\left(\frac{1}{2}\right)}{4}+\frac{105 \pi  \zeta (3)}{128}-\frac{19 \zeta (3)}{32}-\frac{163}{64} \zeta (3) \log (2)$$ $$+\frac{247 \pi ^4}{23040}+\frac{7 \pi ^2}{24}-\frac{11 \pi ^3}{48}-3 \pi -\frac{1}{6} \log ^4(2)+\frac{\log ^3(2)}{4}+\frac{5}{16} \pi  \log ^3(2)+\frac{9}{64} \pi ^2 \log ^2(2)-\frac{3 \log ^2(2)}{2}$$ $$-\frac{13}{8} \pi  \log ^2(2)+\frac{5}{24} \pi ^3 \log (2)-\frac{7}{48} \pi ^2 \log (2)+\frac{5}{2} \pi  \log (2)+6 \log (2)+\frac{3}{256} \left(\psi ^{(3)}\left(\frac{1}{4}\right)-\psi ^{(3)}\left(\frac{3}{4}\right)\right)$$\\

Using $\int _0^1\int _0^1\frac{f(x)-f(y)}{x-y}dxdy=2 \int_0^1 f(t) (\log (t)-\log (1-t)) \, dt$ for certain continuous $f$:\\

{\small $$\int _0^1\int _0^1\frac{\log (x+1) \log \left(x^2+1\right) \tan ^{-1}(x)-\log (y+1) \log \left(y^2+1\right) \tan ^{-1}(y)}{x-y}dxdy$$ $$=\frac{5 \pi ^2 C}{24}+\pi  C+8 C+\frac{1}{2} C \log ^2(2)-\frac{1}{2} \pi  C \log (2)-2 C \log (2)-16 \Im(\text{Li}_3(1+i))+40 \Im(\text{Li}_4(1+i))-4 \log (2) \Im(\text{Li}_3(1+i))$$ $$+2 \text{Li}_4\left(\frac{1}{2}\right)-\frac{49 \pi  \zeta (3)}{32}-\frac{33 \zeta (3)}{8}+\frac{57}{16} \zeta (3) \log (2)+\frac{5 \pi ^3}{12}-\frac{\pi ^2}{4}-\frac{209 \pi ^4}{11520}+\frac{7 \log ^4(2)}{48}-\frac{\log ^3(2)}{3}-\frac{2}{3} \pi  \log ^3(2)$$ $$-\frac{3}{16} \pi ^2 \log ^2(2)+\frac{9}{4} \pi  \log ^2(2)+\log ^2(2)-\frac{37}{96} \pi ^3 \log (2)+\frac{3}{8} \pi ^2 \log (2)-\pi  \log (2)-\frac{7}{384}\left( \psi ^{(3)}\left(\frac{1}{4}\right)-\psi ^{(3)}\left(\frac{3}{4}\right)\right)$$}\\

Now recall subsection 8-1-2 and 8-4-9 for higher weight:\\

$$\int_0^1 \log ^5\left(x^2+1\right) \, dx$$ $$=-1920 \sqrt{2} \, _6F_5\left(\frac{1}{2},\frac{1}{2},\frac{1}{2},\frac{1}{2},\frac{1}{2},\frac{1}{2};\frac{3}{2},\frac{3}{2},\frac{3}{2},\frac{3}{2},\frac{3}{2};\frac{1}{2}\right)+1920 C-40 C \log ^3(2)+240 C \log ^2(2)$$ $$-960 C \log (2)-1920 \Im(\text{Li}_3(1+i))+1920 \Im(\text{Li}_4(1+i))-240 \log ^2(2) \Im(\text{Li}_3(1+i))+960 \log (2) \Im(\text{Li}_3(1+i))$$ $$-960 \log (2) \Im(\text{Li}_4(1+i))-240 \pi  \zeta (3)+360 \pi  \zeta (3) \log (2)+\frac{19 \pi ^5}{3}+70 \pi ^3+960 \pi -3840+\log ^5(2)+\frac{135}{2} \pi  \log ^4(2)$$ $$-10 \log ^4(2)-240 \pi  \log ^3(2)+80 \log ^3(2)+\frac{145}{4} \pi ^3 \log ^2(2)+600 \pi  \log ^2(2)-480 \log ^2(2)-70 \pi ^3 \log (2)$$ $$-960 \pi  \log (2)+1920 \log (2)-\frac{5 \psi ^{(3)}\left(\frac{1}{4}\right)}{16}+\frac{5 \psi ^{(3)}\left(\frac{3}{4}\right)}{16}+\frac{5}{32} \log (2) \psi ^{(3)}\left(\frac{1}{4}\right)-\frac{5}{32} \log (2) \psi ^{(3)}\left(\frac{3}{4}\right)$$

\subsection{NQPLI}

$$\int_0^1 \text{Ti}_2(x) \log (x+1) \log \left(x^2+1\right) \, dx$$ $$=-\frac{\pi ^2 C}{24}+\pi  C+4 C+\frac{1}{2} C \log ^2(2)-\frac{3}{4} \pi  C \log (2)-3 C \log (2)+\pi  \Im(\text{Li}_3(1+i))$$ $$-4 \Im(\text{Li}_3(1+i))-2 \Im(\text{Li}_4(1+i))+3 \log (2) \Im(\text{Li}_3(1+i))-2 \text{Li}_4\left(\frac{1}{2}\right)+\frac{63 \pi  \zeta (3)}{128}$$ $$-\frac{97 \zeta (3)}{32}-\frac{7}{8} \zeta (3) \log (2)+\frac{13 \pi ^3}{96}+\frac{\pi ^2}{12}-\frac{29 \pi ^4}{5760}-3 \pi -\frac{1}{12} \log ^4(2)+\frac{\log ^3(2)}{12}-\frac{1}{6} \pi  \log ^3(2)$$ $$+\frac{1}{48} \pi ^2 \log ^2(2)-\frac{1}{8} \pi  \log ^2(2)-\log ^2(2)-\frac{1}{16} \pi ^3 \log (2)-\frac{1}{48} \pi ^2 \log (2)+2 \pi  \log (2)+6 \log (2)$$

{\small $$\int_0^1 \text{Li}_2\left(-x^2\right) \log (x+1) \log \left(x^2+1\right) \, dx$$ $$=-4 C^2-\frac{\pi ^2 C}{12}-4 \pi  C+16 C+7 C \log ^2(2)+4 \pi  C \log (2)-24 C \log (2)-6 \pi  \Im(\text{Li}_3(1+i))$$ $$+40 \Im(\text{Li}_3(1+i))+40 \Im(\text{Li}_4(1+i))-24 \log (2) \Im(\text{Li}_3(1+i))-2 \text{Li}_4\left(\frac{1}{2}\right)+\frac{73 \zeta (3)}{8}-\frac{7 \pi  \zeta (3)}{4}$$ $$-\frac{63}{16} \zeta (3) \log (2)+\frac{2 \pi ^4}{9}-\frac{13 \pi ^3}{12}-\frac{\pi ^2}{12}+18 \pi -96-\frac{1}{12} \log ^4(2)+\frac{10 \log ^3(2)}{3}+\frac{13}{12} \pi  \log ^3(2)+\frac{5}{16} \pi ^2 \log ^2(2)$$ $$-\frac{1}{2} \pi  \log ^2(2)-25 \log ^2(2)+\frac{5}{48} \pi ^3 \log (2)+\frac{1}{3} \pi ^2 \log (2)-11 \pi  \log (2)+84 \log (2)-\frac{1}{64} \left(\psi ^{(3)}\left(\frac{1}{4}\right)-\psi ^{(3)}\left(\frac{3}{4}\right)\right)$$}\\

\subsection{NS, MLI}
$$\sum _{m=1}^{\infty } \frac{(-1)^{m-1} H_m}{m} \sum _{n=1}^m \frac{(-1)^{n-1} H_n^{(2)}}{n}=\int_0^1 \int_0^1 \frac{\log (x) \log (1-y) \log \left(\frac{1-x y}{1-y}\right)}{(1-x) y (y+1)} \, dxdy$$
$$=6 \text{Li}_5\left(\frac{1}{2}\right)\\+\frac{11 \pi ^2 \zeta (3)}{48}-\frac{49 \zeta (5)}{8}-\frac{1}{2} \zeta (3)\log ^2(2)-\frac{1}{20} \log ^5(2)+\frac{1}{36} \pi ^2 \log ^3(2)+\frac{13}{360} \pi ^4 \log (2)$$

$$\sum _{m=1}^{\infty } \frac{(-1)^{m-1} \left(\widetilde{H_m}\right){}^2 }{m}\sum _{n=1}^m \frac{(-1)^{n-1} H_n}{n}=\int _0^1\int _0^1\int _0^1\frac{\log (1-z) \log (\frac{(1-x) (1-y) (1-z) (1-x y z)}{2 (x y+1) (x z+1) (y z+1)})}{(x+1) (y+1) (z+1)} dxdydz$$
$$=22 \text{Li}_5\left(\frac{1}{2}\right)\\+\frac{\pi ^2 \zeta (3)}{8}-\frac{349 \zeta (5)}{16}-\frac{5}{8} \zeta (3)\log ^2(2)-\frac{17 \log ^5(2)}{60}-\frac{7}{72} \pi ^2 \log ^3(2)+\frac{43}{240} \pi ^4 \log (2)$$

$$\sum _{m=1}^{\infty } \frac{(-1)^{m-1} \widetilde{H_m} }{m}\sum _{n=1}^m \frac{(-1)^{n-1} H_n \widetilde{H_n}}{n}$$
$$=\int _0^1\int _0^1\int _0^1\frac{\log (1-y) \log \left(\frac{(1-y) (1-z)}{2 (y z+1)}\right)}{(x+1) (y+1) (z+1)}dxdydz-\int _0^1\int _0^1\int _0^1\frac{x \log (1-y) \log \left(\frac{(1-z) (x y+1)}{2 (1-x y z)}\right)}{(x+1) (z+1) (1-x y)}dxdydz$$
$$=14 \text{Li}_5\left(\frac{1}{2}\right)+\frac{\pi ^2 \zeta (3)}{12}-\frac{115 \zeta (5)}{8}-\frac{3}{8} \zeta (3)\log ^2(2)-\frac{11 \log ^5(2)}{60}-\frac{1}{18} \pi ^2 \log ^3(2)\\+\frac{23}{180} \pi ^4 \log (2)$$\\

Denote $f(t)=\frac{t}{1+t}$, then by symmetry:
$$ES(1,-1,-1,-1;-1)=\int_{(0,1)^4}\frac{\left(f(-w x y z)-3 f(x y z)+3 f(-x y)-f(x)\right)\log (1-x)}{ x (y+1) (z+1)(w+1)}dx dy dz dw$$
$$=2 \text{Li}_5\left(\frac{1}{2}\right)+\frac{13 \pi ^2 \zeta (3)}{96}-\frac{83 \zeta (5)}{32}\\+\frac{9}{16} \zeta (3) \log ^2(2)-\frac{13}{60}  \log ^5(2)+\frac{41}{72} \pi ^2 \log ^3(2)-\frac{1}{180} \pi ^4 \log (2)$$\\

Let $d\mu= dx dy dz dw, p_1=\frac{x}{(z+1) (1-x y)}+\frac{1}{(y+1) (z+1)}, p_2=\frac{x^2 y}{(1-x y) (x y z+1)}, p_3=\frac{y}{(y+1) (1-y z)}, q_1=f(w)+z f(-z w), q_2=f(w)+x y z f(-x y z w), q_3=f(w)-y z f(y z w), f(t)=\frac{\log(2)-\log(1-t)}{1+t}$, then:
{\small$$\sum _{1\leq j\leq k\leq l\leq m\leq n} \frac{(-1)^{j-1} (-1)^{k-1} (-1)^{l-1} (-1)^{m-1} (-1)^{n-1}}{j k l m n}$$
$$=\int_{(0,1)^4} \frac{p_1 q_1-p_2 q_2+p_3 q_3}{x+1} d\mu =\frac{\pi ^2 \zeta (3)}{48}+\frac{3 \zeta (5)}{16}+\frac{1}{8} \zeta (3) \log ^2(2)+\frac{\log ^5(2)}{120}+\frac{1}{72} \pi ^2 \log ^3(2)+\frac{1}{160} \pi ^4 \log (2)$$}\\

\subsection{Miscellaneous}

By $x\to \tan(u)$, IBP twice, contour deformation and brute force:\\
{\small$$\int_0^{\tan (a)} \tan ^{-1}(x)^n \, dx=n a^{n-1} \log (2 \cos (a))+a^n \tan (a)-n (n-1) \Re(X),$$
$$X=(-1)^{n-1} \left(1-2^{1-n}\right) (n-2)! \zeta (n)-(2 i a)^{n-2} \text{Li}_2\left(-e^{2 i a}\right)+\sum _{k=1}^{n-2} \frac{(-1)^{k-1} (n-2)! (2 i a)^{-k+n-2} \text{Li}_{k+2}\left(-e^{2 i a}\right)}{(2 i)^{n-1} (-k+n-2)!}$$}

\noindent By repeated IBP and definition of $\text{Ls}$ \cite{ref10} :\\
{\small $$\int_0^1 \log ^n\left(x^2+1\right) \, dx=\sum _{k=0}^n \frac{(-2)^k n!} {(n-k)!}\left(\log ^{n-k}(2)-(n-k) (-2)^{n-k} \sum _{j=0}^{-k+n-1} (-\log (2))^j \binom{-k+n-1}{j} A(j,k,n)\right),$$
$$A(j,k,n)=\underset{\{a,b\}\to \{0,0\}}{\text{lim}}\frac{\partial ^{-j-k+n-1}}{\partial a^{-j-k+n-1}}\left(2^{a+b-1} B\left(\frac{a+1}{2},\frac{b+1}{2}\right)\right)+\frac{1}{2} \text{Ls}_{n-k-j}\left(\frac{\pi}{2}\right)$$}

\section{Comments}

\textbf{LI/PLI/QLI/QPLI/ES/QESs evaluation are superficies. MZV theory is the essence.}\\

Superficies: This paper summarizes nearly all known methods of generating LI/PLI/QLI/QPLI relations without using MZVs, resulting in closed-forms of LI/PLI with $W\leq5$, QLI/QPLI with $W\leq 4$. Since for higher weights irreducible constants are inevitable in both cases, in the sense of `expressing these integrals in polylog closed-forms' nothing more can be done. Nevertheless, while ESs are evaluated by totally MZV-free methods, we seek for help from level 4 MZV due to huge difficulties on solving the last 10 QESs. This shows the boundary of elementary methods and reveals the necessity of establishing a deeper theory.\\

Essence: MZV theory offers a more essential way of establishing relations between higher weight LI/ES/PLIs \cite{ref6, ref18}. For QLI/QES/QPLIs, though not fully developed, it is clear that QMZV theory (focusing on quadratic MZVs i.e. multiple polylogarithms of level 4) is their essence. \cite{ref24} provides good example for $W\leq5$, which is sufficient for problems discussed here. Furthermore, a successful investigation on poly-MZVs (level $2n$) will be producing much more seemingly incomprehensible but profound results. For instance, the theory of logsine integrals with argument $\frac \pi 3$ and corresponding Apery-binomial series \cite{ref10}, should be natural consequences of solving level 6 MZV system.\\

\clearpage

\noindent \textbf{\Large Acknowledgements}\\

1. Special thanks to following contributors for crucial improvements of corresponding sections:\\

K. C. Au: For introducing me this general topic through establishment of section 2 and level 4 MZV through that of subsection 7-8, as well as offering important techniques that I used to complete subsection 5-1, 8-4-7. \\

Chen: For introducing me the basic FL theory, offering problems of subsection 8-3-2, 8-4-10(4) and the latter's solution.\\

J. D'Aurizio: For further explanation on FL theory, suggestions on hypergeometric representation in subsection 8-3-4 part 2, and derivation of FL expansion of $\frac{f(x)}x$ in part 3, given that of $f(x)$.\\

A. Shather: For suggestions on double integration in subsection 6-10-1, 8-5-1, offering problem of 8-4-9(6) and its solution.\\

C. I. Valean: For his challenging problem book \cite{ref15}, offering insights on QLI/QPLI evaluation from which I developed section 5, 6 later.\\\\

2. I also appreciate following Math. SE users for presenting interesting problems, solutions on related topics, or minor suggestions on this article:\\

Cody, FDP, M.N.C.E., Renascence 5, Tunk-Fey, user 178256, user 1591719, user 97357329, xuce1234, Zacky. Please forgive if I've forgotten someone.\\\\

3. Last but not least, I am grateful to my family and friends for always supporting me during this paper's reformulation.\\\\

Thank you all.\\

\end{document}